%% file: main.tex
\documentclass[12pt]{amsart}

\usepackage{amsmath}
\usepackage{amsfonts}
\usepackage{amssymb}
\usepackage{amscd}
\usepackage{graphicx}
\usepackage[abbrev,alphabetic]{amsrefs}
\RequirePackage[dvipsnames,usenames]{color}
\usepackage{soul,xcolor}
\setstcolor{red}
\usepackage{stmaryrd}
\usepackage{mathtools}
\usepackage{booktabs}
\usepackage{multirow}
\newtagform{tiny}{\tiny(}{)}

\newenvironment{smallertags}{\usetagform{tiny}}{\ignorespacesafterend}

\usepackage{mathtools}
\usepackage{hyperref}
\usepackage[margin=1.25in]{geometry}

\usepackage{amsthm}
\usepackage{comment}
\usepackage[all,cmtip]{xy}
\usepackage{tikz-cd}
\usetikzlibrary{cd}

\usepackage[all]{xy}

\makeatletter
\def\@tocline#1#2#3#4#5#6#7{\relax
  \ifnum #1>\c@tocdepth 
  \else
    \par \addpenalty\@secpenalty\addvspace{#2}%
    \begingroup \hyphenpenalty\@M
    \@ifempty{#4}{%
      \@tempdima\csname r@tocindent\number#1\endcsname\relax
    }{%
      \@tempdima#4\relax
    }%
    \parindent\z@ \leftskip#3\relax \advance\leftskip\@tempdima\relax
    \rightskip\@pnumwidth plus4em \parfillskip-\@pnumwidth
    #5\leavevmode\hskip-\@tempdima
      \ifcase #1
       \or\or \hskip 1em \or \hskip 2em \else \hskip 3em \fi%
      #6\nobreak\relax
    \hfill\hbox to\@pnumwidth{\@tocpagenum{#7}}\par
    \nobreak
    \endgroup
  \fi}
\makeatother

\newsavebox{\pullback}
\sbox\pullback{%
\begin{tikzpicture}%
\draw (0,0) -- (1ex,0ex);%
\draw (1ex,0ex) -- (1ex,1ex);%
\end{tikzpicture}}

\newsavebox{\pullbackdl}
\sbox\pullbackdl{%
\begin{tikzpicture}%
\draw (-1ex,0ex) -- (0ex,0ex);%
\draw (0ex,-1ex) -- (0ex,0ex);%
\end{tikzpicture}}

\newcommand{\cyan}{\color{cyan}}

\newcommand{\cred}{\color{black}}

\newcommand{\rdown}[1]{\lfloor #1 \rfloor}

\renewcommand{\mod}{\ \textrm{mod}\ }
 
\renewcommand{\P}{\mathbb{P}}
\newcommand{\Z}{\mathbb{Z}}
\newcommand{\Q}{\mathbb{Q}}
\newcommand{\R}{\mathbb{R}}
\newcommand{\F}{\mathbb{F}}

\newcommand{\bC}{\mathbb{C}}

\newcommand{\bE}{\mathbb{E}}

\newcommand{\bP}{\mathbb{P}}
\newcommand{\bQ}{\mathbb{Q}}

\newcommand{\bZ}{\mathbb{Z}}

\newcommand{\cO}{\mathcal{O}}

\newcommand{\MO}{\mathcal{O}}
\newcommand{\sO}{\mathcal{O}}

\newcommand{\m}{\mathfrak{m}}

\newcommand{\red}{\mathrm{red}}

\newcommand{\Proj}{\mathrm{Proj}}

\DeclareMathOperator{\Supp}{Supp}
\DeclareMathOperator{\Spec}{Spec}
\DeclareMathOperator{\Diff}{Diff}
\DeclareMathOperator{\codim}{codim}
\DeclareMathOperator{\Hom}{Hom}

\DeclareMathOperator{\Pic}{Pic}
\DeclareMathOperator{\NE}{NE}
\DeclareMathOperator{\Exc}{Exc}
\DeclareMathOperator{\Ex}{Exc}

\DeclareMathOperator{\coeff}{coeff}
\DeclareMathOperator{\mult}{mult}

\newcommand{\mydot}{{{\,\begin{picture}(1,1)(-1,-2)\circle*{2}\end{picture}\ }}}

\theoremstyle{plain}
\newtheorem{theorem}{Theorem}[section]

\newtheorem{proposition}[theorem]{Proposition}
\newtheorem{lemma}[theorem]{Lemma}
\newtheorem{corollary}[theorem]{Corollary}
\newtheorem{conjecture}[theorem]{Conjecture}

\newtheorem{question}[theorem]{Question}
\newtheorem{claim}[theorem]{Claim}
\newtheorem*{claim*}{Claim}
\newtheorem{step}{Step}

\newtheorem{theoremA}{Theorem}

\theoremstyle{definition}
\newtheorem{definition}[theorem]{Definition}

\newtheorem{example}[theorem]{Example}
\newtheorem{notation}[theorem]{Notation}
\newtheorem{nothing}[theorem]{}
\newtheorem*{setup*}{Setup}

\theoremstyle{remark}
\newtheorem{remark}[theorem]{Remark}

\def\todo#1{\textcolor{Mahogany}%
{\footnotesize\newline{\color{Mahogany}\fbox{\parbox{\textwidth-15pt}{\textbf{todo: } #1}}}\newline}}



\numberwithin{equation}{theorem}

\title[Quasi-F-splittings in birational geometry II]
{Quasi-F-splittings in birational geometry II}

\author{Tatsuro Kawakami}
\address{Department of Mathematics, Graduate School of Science, Kyoto University, Kyoto 606-8502, Japan} 
\email{tkawakami@math.kyoto-u.ac.jp}
\author{Teppei Takamatsu}
\address{Department of Mathematics, Graduate School of Science, Kyoto University, Kyoto 606-8502, Japan}
\email{teppeitakamatsu.math@gmail.com}
\author{Hiromu Tanaka} 
\address{Graduate School of Mathematical Sciences, 
The University of Tokyo, 
3-8-1 Komaba, Meguro-ku, Tokyo 153-8914, JAPAN} 
\email{tanaka@ms.u-tokyo.ac.jp}
\author{Jakub Witaszek} 
\address{Department of Mathematics \\
Fine Hall, Washington Road\\
Princeton University\\  
Princeton NJ 08544-1000, USA}
\email{jwitaszek@princeton.edu}
\author{Fuetaro Yobuko}
\address{Graduate School of Mathematics, Nagoya University, Furocho, Chikusaku, Nagoya 464-0813, Japan}
\email{soratobumusasabidesu@gmail.com}
\author{Shou Yoshikawa}
\address{Tokyo Institute of Technology, Tokyo 152-8551, Japan}
\email{yoshikawa.s.al@m.titech.ac.jp}

\begin{document}

\begin{abstract}
Over an algebraically closed field of characteristic $p>41$, 
we prove that three-dimensional $\mathbb Q$-factorial affine klt varieties are quasi-$F$-split. 
Furthermore, we show that the bound on the characteristic is optimal. 
\end{abstract}

\maketitle

\setcounter{tocdepth}{2}

\tableofcontents
\input section1.tex

\input section2.tex
\input section3.tex
\input section4.tex

\input section5.tex

\input section6.tex

\input section7.tex

\input section8.tex
\input section9.tex

\input section10.tex
\input section11.tex

\input main.bbl


\end{document}

%% file: section1.tex
\section{Introduction}

This paper is a continuation of \cite{KTTWYY1}.  We refer thereto for the background and motivation behind studying quasi-$F$-splittings in birational geometry.

Given an algebraic variety $V$ in positive characteristic, 
we say that $V$ is {\em $n$-quasi-$F$-split} if there exist 
a $W_n\MO_V$-module homomorphism $\alpha : F_* W_n\MO_V \to \MO_V$ which completes the following commutative diagram: 
\begin{equation*} 
\begin{tikzcd}
W_n\cO_V \arrow{r}{{F}} \arrow{d}{R^{n-1}} & F_* W_n \arrow[dashed]{ld}{\exists\alpha} \cO_V \\
\cO_V. 
\end{tikzcd}
\end{equation*}
Further, we say that $V$ is \emph{quasi-$F$-split} if it is $n$-quasi-$F$-split for some $n \in \bZ_{>0}$.
In the first part \cite{KTTWYY1}, we have proven the following results. 
{\setlength{\leftmargini}{2.7em}
\begin{enumerate}
\item[(I)] Every affine klt surface $S$ is quasi-$F$-split. 
\item[(II)] Every three-dimensional $\Q$-factorial affine klt variety $X$ is $2$-quasi-$F$-split 
if $p \gg  0$. 
\end{enumerate}
}

\noindent The purpose of this article is to give the optimal bound for when $X$ is quasi-$F$-split as explained by the following results.
\begin{theoremA}[Theorem \ref{t-3dim-klt-QFS}]\label{tA-3klt-QFS}
Let $k$ be 
a perfect field of characteristic $p >41$. 
Let $V$ be a three-dimensional $\Q$-factorial affine klt variety over $k$. 
Then $V$ is quasi-$F$-split. 
In particular, $V$ lifts to $W_2(k)$. 
\end{theoremA}

\begin{theoremA}[Theorem \ref{t-3dim-klt-nonQFS}]\label{tA-3klt-nonQFS}
Let $k$ be an algebraically closed field of characteristic $41$. 
Then there exists a three-dimensional $\Q$-factorial affine klt variety $V$ over $k$ which is not quasi-$F$-split. 
\end{theoremA}

Roughly speaking, Theorem \ref{tA-3klt-QFS} and Theorem \ref{tA-3klt-nonQFS} 
are reduced, by {\cred the} global-to-local correspondence in birational geometry, 
to the following corresponding problems for log del Pezzo pairs with standard coefficients. 
For the definition of quasi-$F$-split pairs $(X, \Delta)$, see Subsection \ref{ss-def-QFS}. 

\begin{theoremA}[Theorem \ref{t-LDP-QFS}]\label{tA-LDP-QFS}
Let $k$ be a 
perfect field of characteristic $p >41$. 
Let $(X, \Delta)$ be a log del Pezzo pair over $k$ with standard coefficients, that is, 
$(X, \Delta)$ is a two-dimensional projective klt pair such that $-(K_X+\Delta)$ is ample and 
all coefficients of $\Delta$ are contained in $\{ 1 - \frac{1}{n} \,|\, n \in \Z_{>0}\}$. 
Then $(X, \Delta)$ is quasi-$F$-split. 
\end{theoremA}

\begin{theoremA}[Theorem \ref{t-P^2-main}]\label{tA-LDP-nonQFS}
Let $k$ be an algebraically closed field of characteristic $41$.  
Then there exists a log del Pezzo pair $(X, \Delta)$ over $k$ with standard coefficients which is not quasi-$F$-split. 
\end{theoremA}

We emphasise that Theorem \ref{tA-3klt-QFS} and Theorem \ref{tA-LDP-QFS} are false with the notion of quasi-$F$-split replaced by $F$-split even if $p$ is arbitrarily large, as shown in \cite{CTW15a}.\\

One of the main difficulties of positive characteristic birational geometry is that klt singularities are, in this setting, little understood. In fact, some of their properties, such as Cohen-Macauliness, break already in dimension three for low characteristics (cf.~\cite{CT17}). Theorem \ref{tA-3klt-QFS} indicates that these singularities should, however, behave very similarly to their characteristic {\cred zero} counterparts when, at least, $p>41$. 

In what follows we provide another example, besides liftability, supporting this statement, by generalising the result in \cite{Kaw4} of the first author about the logarithmic extension theorem from $F$-injective three-dimensional terminal singularities to quasi-$F$-split ones. In particular, this result holds for all three-dimensional terminal singularities when $p > 41$.

\begin{theoremA}[Corollary \ref{c-3dim-1form}]\label{tA-1form}
Let $V$ be a three-dimensional terminal variety over a perfect field 
of characteristic $p>41$.  
Then $V$ satisfies the logarithmic extension theorem for one-forms. 
More explicitly, for any proper birational morphism $f\colon W\to V$ from a normal variety $W$, the natural injective homomorphism 
\[
f_{*}\Omega^{[1]}_W(\log\,E)\hookrightarrow \Omega^{[1]}_V. 
\]
is an isomorphism, 
where $E$ denotes the reduced Weil divisor satisfying $\Supp E = \Ex(f)$. 
\end{theoremA}

\subsection{Overview of the proofs}

We now give an overview of how to prove Theorem \ref{tA-3klt-QFS}--\ref{tA-LDP-nonQFS}. 
Following global-to-local  correspondence in birational geometry, 
we shall establish the following two implications (I) and (II).  
\begin{enumerate}
\item[(I)] Theorem \ref{tA-LDP-QFS} $\Rightarrow$ Theorem \ref{tA-3klt-QFS}. 
\item[(II)] Theorem \ref{tA-LDP-nonQFS} $\Rightarrow$ Theorem \ref{tA-3klt-nonQFS}. 
\end{enumerate}
The implication (I) has already been proven in \cite[{\cred the proof of Theorem 6.19}]{KTTWYY1}. To show the implication (II), we construct an example of a non-quasi-$F$-split three-dimensional klt singularity by taking an orbifold cone over the non-quasi-$F$-split log del Pezzo surface from Theorem \ref{tA-LDP-nonQFS}. The details are contained in Section \ref{s-Qcone-QFS} and Section \ref{s-Qcone-klt}, where we prove the following:
\begin{itemize}
    \item a log del Pezzo surface with standard coefficients is quasi-$F$-split if and only if its orbifold cone is so (Corollary \ref{cor:corresponding});
    \item an orbifold cone over a Picard-rank-one log del Pezzo surface with standard coefficients is $\bQ$-factorial and klt (Theorem \ref{t-klt})\footnote{this result was readily believed to be true by the experts, but we could not find a reference that works beyond the case of characteristic zero}.
\end{itemize}

\subsubsection{Sketch of the proof of Theorem \ref{tA-LDP-nonQFS}}\label{sss-LDP-nonQFS}
The example in the statement of Theorem \ref{tA-LDP-nonQFS} is explicitly given by 
\[
(X, \Delta) =\left( \P^2, \frac{1}{2}L_1 + \frac{2}{3} L_2 + \frac{6}{7} L_3 + \frac{40}{41}L_4\right),
\]
where $L_1, L_2, L_3, L_4$ are lines such that $L_1 + L_2 + L_3 + L_4$ is simple normal crossing. 
It is easy to check that $(X, \Delta)$ is a log del Pezzo pair. 
What is remaining is to prove that $(X, \Delta)$ is not quasi-$F$-split. 
Recall that $(X, \Delta)$ is quasi-$F$-split if and only if there exists $m \in \Z_{>0}$ 
such that $\Phi_{X, \Delta, m}: \MO_X \to Q_{X, \Delta, m}$ is a split injection (cf.\ Subsection \ref{ss-def-QFS}), 
{\cred where $\Phi_{X, \Delta, m}: \MO_X \to Q_{X, \Delta, m}$ is defined by the following pushout diagram of $W_m\MO_X$-modules: 
\[
\begin{tikzcd}
W_m\cO_X(\Delta) \arrow{r}{{F}} \arrow{d}{R^{m-1}} & F_* W_m\cO_X(p\Delta)  \arrow{d} \\
\MO_X \arrow{r}{\Phi_{X, \Delta, m}} & Q_{X, \Delta, m}. 
\end{tikzcd}
\]
}
\noindent 
Fix $m \in \Z_{>0}$ and set $D := K_X+\Delta$. 
By $\Phi_{X, \Delta, m} \otimes \MO_X(K_X) = \Phi_{X, D, m}$, it is enough to show that 
\[
H^2(\Phi_{X, D, m}) : H^2(X, \MO_X(D)) \to H^2(X, Q_{X, D, m})
\]
is not injective. 
Set $\Delta_{\red} := L_1 + L_2 +L_3 + L_4$. 
By the following exact sequence (cf.\ Lemma \ref{lem:Serre's map}): 
\[
0 \to \MO_X(D) \xrightarrow{\Phi_{X, D, n}} Q_{X, D, n} \to B_m\Omega_X^1(\log \Delta_{\red})(p^mD) \to 0, 
\]
it suffices to show that the connecting map 
\[
\delta_m : H^1(X,B_m\Omega_X^1(\log \Delta_{\red})(p^mD)) \to H^2(X,\MO_X(D))
\]
is nonzero. 
If $m=1$, then this follows from $H^1(X, Q_{X, D, 1})=H^1(X, F_*\MO_X(pD)) \simeq H^1(X, \MO_X(\rdown{pD}))= 0$ and $H^2(X, B_m\Omega_X^1(\log \Delta_{\red})(p^mD)) \neq 0$, 
where the latter one can be checked by the Euler sequence and standard exact sequences on log higher Cartier operators. 
For the general case, the strategy is to apply induction on $m$ via  the following commutative diagram: 
\begin{equation*}
\begin{tikzcd}
H^1(B_{m}\Omega_X^1(\log \Delta_{\red})(p^{m}D)) \arrow[d, "C"']\arrow{r}{{\delta_{m}}} & H^2(\MO_X(D)) \arrow[equal]{d}   \\
H^1(B_{m-1}\Omega_X^1(\log \Delta_{\red})(p^{m-1}D)) \arrow{r}{{\delta_{m-1}}} & H^2(\MO_X(D))
\end{tikzcd}
\end{equation*}
Again by the Euler sequence and log higher Cartier operators, 
we can check that $C$ is an isomorphism. 
Therefore, $\delta_m$ is nonzero, because so is $\delta_1$. 
For more details, see Section \ref{s-P^2-cex}.

\medskip

\subsubsection{Sketch of the proof of Theorem \ref{tA-LDP-QFS}}\label{sss-LDP-QFS}
First, let us recall \cite[Conjecture 6.20]{KTTWYY1}.

\begin{conjecture} \label{conj:log-liftability-intro}
Let $(X,\Delta)$ be a log del Pezzo pair over a perfect field of characteristic $p>5$ and such that $\Delta$ has standard coefficients. Then there exists {a} 
log resolution $f \colon Y \to X$ of $(X, \Delta)$ such that $(Y, f^{-1}_*\Delta + \Exc(f))$ lifts to {$W(k)$}. 
\end{conjecture}
\noindent In other words, the statement stipulates that $(X,\Delta)$ is \emph{log liftable} to $W(k)$.

As indicated in \cite[Remark 6.21]{KTTWYY1}, Theorem \ref{tA-LDP-QFS} follows from this conjecture by using  
the method of higher Cartier {\cred operators}. For the convenience of the reader, we summarise this argument. Write $K_Y+\Delta_Y = f^*(K_X+\Delta)$ and set $E := \Supp(f^{-1}_*\Delta + \Exc(f))$. 
By Theorem \ref{t-QFS-criterion}, it is enough to verify that
\begin{enumerate}
    \item[(B)] $H^0(Y, \MO_Y(K_Y+E+p^l(K_Y+\Delta_Y))) =0$ for every $l \geq 1$, and
    \item[(C)] $H^0(Y, \Omega_Y(\log E)(K_Y+B_Y))=0$.
\end{enumerate}
Here, (C) follows by a nef-and-big Akizuki-Nakano vanishing (see \cite[Theorem 2.11]{Kaw3}) contingent upon Conjecture \ref{conj:log-liftability-intro}. Interestingly, (B) is false in general if we do not assume that $p> 42$. The validity of (B) under this assumption follows readily from the ACC statement for the nef threshold for surfaces (this is exactly where the constant $42$ comes into the picture) which in characteristic $0$ is known by Koll\'ar (\cite{Kol94}). We extend this result (alas only for log del Pezzo surfaces) to positive characteristic in Subsection \ref{ss-explicit-ACC} by lifting to characteristic $0$ and applying Koll\'ar's original statement.\\

Therefore, to understand Theorem \ref{tA-LDP-QFS}, we need to address Conjecture \ref{conj:log-liftability-intro}. Note that log liftability of del Pezzo pairs (that is, the case $\Delta=0$) for $p>5$ is already known by  Arvidsson-Bernasconi-Lacini \cite{ABL20}. In trying to generalise their result to pairs with standard coefficients, we first show the following local result which we believe is of independent interest (both in characteristic zero and positive characteristics).

\begin{lemma}[{Corollary \ref{cor:Hiromu-resolution}}] \label{lem:min-resolutions-intro}
Let $(X,\Delta)$ be a two-dimensional klt pair with standard coefficients 
over an algebraically closed field. Then there exists a projective birational morphism $\varphi \colon V \to X$ such that $(V,\Delta_V)$ is simple normal crossing and $\Delta_V$ is effective, where $\Delta_V$ is the $\bQ$-divisor defined by $K_V+\Delta_V = \varphi^*(K_X+\Delta)$.
\end{lemma}

One may ask if a stronger statement is true, that 
\begin{equation} \label{eq:stronger-minimal-resolution}
    (V, \varphi^{-1}_*\Delta + \Exc(\varphi)) 
    \textrm{ is simple normal crossing.} 
\end{equation}
By minimal resolutions, we may assume that (\ref{eq:stronger-minimal-resolution}) is valid when $\Delta=0$. Unfortunately, in general, one cannot construct a resolution for which $\Delta_V\geq 0$ and  (\ref{eq:stronger-minimal-resolution}) holds. A counterexample is as follows:
\begin{itemize}
    \item $x \in X$ a smooth point, and
    \item $\Delta=\frac{1}{2}C$ with $C$ being an irreducible curve with a cusp\footnote{a singularity looking like $(y^2 =x^{2n+1})$} at $x$.
\end{itemize} 
We shall call such a point $x$, a special point of $(X,\Delta)$. On the positive side, we can show that special points are the only possible counterexamples to constructing resolutions satisfying $\Delta_V \geq 0$ and  (\ref{eq:stronger-minimal-resolution}) (see Theorem \ref{t-hiromu-resol}).

Building on the above lemma and \cite{ABL20}, we can prove the following weaker version of Conjecture \ref{conj:log-liftability-intro}.

\begin{theorem}[Corollary \ref{cor:Hiromu-resolution}, 
Theorem \ref{thm:liftability-of-Hiromu-resolution}] \label{thm:liftability-of-log-del-Pezzo-intro}
Let $(X,\Delta)$ be a log del Pezzo pair with standard coefficients over an algebraically closed field of characteristic $p>5$. Then there exists a projective birational morphism $f \colon Y \to X$ such that $(Y, \Delta_Y)$ is log smooth and lifts over $W(k)$, where $K_Y + \Delta_Y = f^*(K_X+\Delta)$ and $\Delta_Y \geq 0$.
\end{theorem}

Now, let us go back to the proof of the Theorem \ref{t-3dim-klt-QFS}. In what follows, we shall call a special point $x$ of $(X,\Delta)$ to be 
{\em of type} $n$, if $n$ is the minimum number of blowups which 
resolve $(X,\Delta)$ at $x$. We make the following observations.
\begin{itemize}
    \item If there are no special points of $(X,\Delta)$, then we can construct a log resolution as in Conjecture  \ref{conj:log-liftability-intro}, and so the sketch of the argument above shows that $(X,\Delta)$ is quasi-$F$-split.
    \item If every special point of $(X,\Delta)$ is of type $n$ for $2n+1 < p$, then we can still recover vanishing (C) by combining the above strategy with ideas of Graf (\cite{graf21}) on extensions of logarithmic forms\footnote{We consider a sequence of birational morphisms $f \colon Y \xrightarrow{\theta} V \xrightarrow{\varphi} X$, where $\varphi$ is as in Lemma \ref{lem:min-resolutions-intro} and $\varphi \circ \theta$ is a log resolution of $(X,\Delta)$. Specifically, we choose them so that $\theta$ is identity over all points $x \in X$ unless $x$ is special. Then, by the inequality $2n+1 < p$, the determinant $(-1)^n(2n+1)$ of a suitable intersection matrix of $\Exc(\theta)$ is not divisible by $p$ (cf.\  Remark \ref{r-cusp-resol}), and so we can apply the argument of Graf to show that $H^0(Y, \Omega^1_Y(\log E)(K_Y+\Delta_Y)) = H^0(Y, \Omega^1_Y(\log \Supp \Delta_Y)(K_Y+\Delta_Y))$, which is then zero by Akizuki-Nakano and our Theorem \ref{thm:liftability-of-log-del-Pezzo-intro}; see Subsection \ref{t-BS-vanishing}.}.  
    Thus, we again get that $(X,\Delta)$ is quasi-$F$-split.  
\end{itemize}

Thus, from now on, in our explanation we may assume that there exists a special point $x$ of $(X,\Delta)$ of type $n$ for $42 \leq p \leq 2n+1$. Since $K_X+\Delta$ is anti-ample, the existence of a curve $\frac{1}{2}C \subseteq \Delta$ with singularity of high multiplicity, imposes strict restriction on the geometry of $(X,\Delta)$. Thus, one may hope to classify all such $(X,\Delta)$. Informally speaking, $(X,\Delta)$ should be birational to a Hirzebruch surface $(\mathbb{F}_n, \Delta_{\mathbb{F}_n})$ with large $n$ and simple $\Supp \Delta_{\mathbb{F}_n}$. Granted that, one should be able to prove Conjecture \ref{conj:log-liftability-intro} for such $(X,\Delta)$ by hand, concluding the proof that it is quasi-$F$-split.\\

In what follows, we provide more details (in practice we avoid the use of explicit classifications by a mix of intersection and deformation theories). But before doing so, let us note that a natural first step would be to take the minimal resolution 
$f \colon Y \to X$ of $X$ and run a $K_Y$-MMP to construct a Hirzebruch surface as above. Unfortunately, if we write $K_Y+\Delta_Y = f^*(K_X+\Delta)$, then $\Delta_Y$ will not have standard coefficients anymore. In fact, its coefficients may be arbitrarily small, and so hard to control. 

To circumvent this problem, we use 
canonical models over $X$ as suggested to us by J.\ Koll\'ar. Specifically:
\begin{enumerate}
    \item we take the canonical model 
    $f \colon Y \to X$ over $X$. Then the singularities of $Y$ are canonical, and hence Gorenstein.
    \item The coefficients of $\Delta_Y$, defined by $K_Y+\Delta_Y =  f^*(K_X+\Delta)$,
    are at least $\frac{1}{3}$.
    \item Then we run a $K_Y$-MMP $g \colon Y \to Z$, so that 
        \begin{itemize}
        \item[(i)] $\rho(Z)=1$ or
        \item[(ii)] $Z$ admits a Mori fibre space $\pi \colon Z \to \mathbb{P}^1$. 
        \end{itemize}
    Set $\Delta_Z := g_*\Delta_Y$ and $C_Z := g_*f^{-1}_*C$, where $C_Z$ is still highly singular.
    \item Now, canonical del Pezzo surfaces are bounded, and given that $\frac{1}{2}C_Z \leq \Delta_Z$ with singularities of high multiplicity, one can check by intersection theory that the case (i) does not happen, and so $Z$ admits a Mori fibre space $\pi \colon  Z \to \mathbb{P}^1$.
    \item Since $(K_Z+\Delta_Z) \cdot F <0$ for a generic fibre $F$ of $\pi$, we get that $0 < \Delta_Z \cdot F < 2$. With a little bit of work, one can then check that $(\Delta_Z)_{\rm red} \cdot F \leq 3$.
    {\cred From this, we can deduce the vanishing of $H^2(Z,T_Z(-\log\,(\Delta_Z)_{\rm red}))$, which contains the obstruction of the log lifting of $(Z,\Delta_Z)$. Thus $(Z,\Delta_Z)$ is log liftable, and so is $(X,\Delta)$ (see Lemma \ref{l-vertical-lift}; thanks to this lemma, we do not need to classify $(Z,\Delta_Z)$ explicitly).}
\end{enumerate}

\medskip
\noindent {\bf Acknowledgements.}
The authors thank J\'anos Koll\'ar and Yuya Matsumoto 
for very valuable conversations related to the content of the paper. 
The authors also thank the referee for 
constructive suggestions and  
reading the manuscript carefully.

\begin{itemize}
    \item Kawakami was supported by JSPS KAKENHI Grant number JP22KJ1771.
    \item Takamatsu was supported by JSPS KAKENHI Grant number JP22KJ1780.
    \item Tanaka was supported by JSPS KAKENHI Grant numbers JP18K13386, JP22H01112, and JP23K03028.
    \item Witaszek was supported by NSF research grant DMS-2101897.
    \item Yobuko was supported by JSPS KAKENHI Grant number JP19K14501.
    \item Yoshikawa was supported by JSPS KAKENHI Grant number JP20J11886 and RIKEN iTHEMS Program. 
\end{itemize}

%% file: section2.tex
\section{Preliminaries}

\subsection{Notation}\label{ss-notation}

\begin{enumerate}
    \item 
    Throughout this paper, we work over an algebraically closed field $k$ of characteristic $p>0$, unless otherwise specified. However, we emphasise that all the results in the preliminaries also hold over perfect fields. Moreover, at the end of the paper we deduce from the algebraically-closed-field case that the main results of our article hold over perfect fields, too.
    
    We shall freely use notation and terminologies from \cite{hartshorne77} and \cite{KTTWYY1}. 
    \item 
    We say that $X$ is a {\em variety} (over field $k$) if 
$X$ is an integral scheme which is separated and of finite type over $k$. 
We say that $X$ is a {\em surface} 
(resp.\ a {\em curve}) 
if $X$ is a variety of dimension two (resp.\ one). 
\item We say that $(X, \Delta)$ is {\em log smooth} if $X$ is a smooth variety and 
$\Delta$ is a simple normal crossing effective $\Q$-divisor. 
\item\label{ss-n-log-lift} 
Let $X$ be a normal surface and let $D$ be a $\Q$-divisor on $X$. 
We say that $(X, D)$ is {\em log liftable} if 
there exists a log resolution $f : Y \to X$ of $(X, D)$ such that $(Y, f_*^{-1}D_{\red} +\Ex(f))$ lifts to $W(k)$. 
If $(X, D)$ is log liftable, then $(Z, g_*^{-1}D_{\red} +\Ex(g))$ lifts to $W(k)$ for any log resolution $g : Z \to X$ of $(X, D)$ (Proposition \ref{p-log-lift-any}). 
\item 
Given a normal variety $X$ and a $\Q$-divisor $D$, 
recall that the subsheaf $\MO_X(D)$ of the constant sheaf $K(X)$ is defined by 
\[
\Gamma(U, \MO_X(D)) := \{ \varphi \in K(X) \,|\, ({\rm div}\,(\varphi) + D)|_U \geq 0\}. 
\]
Clearly, we have $\MO_X(D) = \MO_X(\rdown{D})$. 
Moreover, when $X = \Spec R$, we set $R(D) := \Gamma(X, \MO_X(D))$. 
\item Given a normal variety $X$ and a $\Q$-divisor $\Delta$, 
$\coeff (\Delta)$ denotes its coefficient set. More precisely, 
if $\Delta = \sum_{i \in I} \delta_i \Delta_i$ is the irreducible decomposition 
such that $\delta_i\neq 0$ for every $i \in I$, 
then $\coeff (\Delta) := \{ \delta_i \,|\, i \in I\} \subseteq \Q$. 
Similarly, a coefficient of $\Delta$ is always assumed to be nonzero under our notation.
\item Given a normal variety $X$ and a $\Q$-divisor $D$ with the irreducible decomposition $D = \sum_{i \in I} d_iD_i$, 
we set $D_{\red} := \sum_{i \in I} D_i$. 
{We say that $D$ is {\em simple normal crossing} 
if for every closed point $x \in D_{\red}$, $X$ is smooth at $x$ and there exists a regular system of parameters $x_1,\ldots, x_d$ in the maximal ideal $\m$ of $\cO_{X,x}$ and $1 \leq r \leq d$ such that $D_{\red}$ is defined by $x_1\cdots x_r$ in $\cO_{X,x}$.} 
\item Let $X$ be a normal variety, let $D$ be a reduced divisor on $X$, and let  $j\colon U\to X$ be the inclusion morphism from the log smooth locus $U$ of $(X,D)$. We denote the sheaf of $i$-th reflexive logarithmic differential forms $j_{*}\Omega_U^{i}(\log D|_U)$ by $\Omega_X^{[i]}(\log D)$.
\end{enumerate}

\subsection{Log liftability} 

\begin{lemma}\label{lemma:lift of bl-down}
    Let $Y$ be a smooth projective surface and let $E$ be a $(-1)$-curve on $Y$.
    Let $f\colon Y\to X$ be the contraction of $E$.
    Suppose that there exists a lift $\mathcal{Y}$ over $W(k)$ of $Y$.
    Then there exist a lift $\mathcal{X}$ of $X$ and $\tilde{f}\colon \mathcal{Y}\to \mathcal{X}$ of $f$ over $W(k)$. 
\end{lemma}
\begin{proof}
    Since $H^1(E, N_{E|_{{\cred Y}}})=H^1(E, \sO_E(-1))=0$, there exists a lift $\mathcal{E}$ of $E$ {\cred by \cite[Theorem 6.2(a)]{hartshorne_deformation}.}
    Let $\mathcal{L}$ be an ample Cartier divisor on $\mathcal{Y}$ and $L\coloneqq \mathcal{L}|_{Y}$.
    Take $\lambda\in\Q_{>0}$ so that $(L+\lambda E)\cdot E=0$.
    Then $L+\lambda E-K_Y$ is nef and big over $X$ and the base point free theorem shows that there exists an ample Cartier divisor $A$ on $X$
     and $m\in\Z_{>0}$ such that $m(L+\lambda E)=f^{*}A$. Replacing $A$ with its multiple, we can assume that $H^i(X, \sO_X(A))=0$ for every $i>0$.
     Since $R^if_{*}\sO_Y(f^{*}A)=\sO_X(A)\otimes R^{i}f_{*}\sO_X=0$, we have $H^i(Y, \sO_Y(m(L+\lambda E)))= H^i(Y, \sO_Y(f^{*}A))\cong H^i(X,\sO_X(A))=0$.
     Then it follows from the Grauert theorem that $H^0(\mathcal{Y}, \sO_{\mathcal{Y}}(\mathcal{L}+\lambda \mathcal{E}))$ is a free $W(k)$-module and
     \[
     H^0(\mathcal{Y}, \sO_{\mathcal{Y}}(\mathcal{L}+\lambda \mathcal{E}))\otimes_{W(k)}k \cong H^0(Y, \sO_{Y}(L+\lambda E)).
     \]
     Let $\mathcal{Z}$ be the base locus of $\mathcal{L}+\lambda \mathcal{E}$. Then the image of $X \to \Spec\,W(k)$ is a closed subset since $\mathcal{Y}$ is proper over $\Spec\,W(k)$.
     Since $L+\lambda E$ is base point free, the base locus $\mathcal{Z}$ should be the empty set.
     Now, the morphism associated to $\mathcal{L}+\lambda \mathcal{E}$ gives a lift of the morphism associated to $L+\lambda E$, which is $f$.
\end{proof}

\begin{proposition}\label{p-log-lift-any}
    Let $(X, D)$ be a pair of a normal projective surface and a $\Q$-divisor such that $(X,D)$ is log liftable.
    Then for any log resolution $f\colon Y\to X$ of $(X,D)$, the pair $(Y, f^{-1}_{*}D+\Exc(f))$ lifts to $W(k)$.
\end{proposition}
\begin{proof}
    The assertion follows from Lemma \ref{lemma:lift of bl-down} and the proof of \cite[Lemma 2.9]{Kaw3}.
\end{proof}

\subsection{Singularities in 
the minimal model program}

We recall some notation in the theory of singularities in the minimal model program. 
For more details, we refer the reader to \cite[Section 2.3]{KM98} and \cite[Section 1]{kollar13}.

We say that $(X, \Delta)$ is a \textit{log pair} 
if $X$ is a normal variety over $k$ and $\Delta$ is an effective $\Q$-divisor such that 
$K_X + \Delta$ is $\Q$-Cartier. 
For a proper birational morphism $f: X' \to X$ from a normal variety $X'$ 
and a prime divisor $E$ on $X'$, the \textit{discrepancy} of $(X, \Delta)$ 
at $E$ is defined as
\[
a(E, X, \Delta) := \text{the coefficient of }E\text{ in }K_{X'} - f^* (K_X + \Delta). 
\]

We say that $(X, \Delta)$ is \textit{klt} (resp.\ \textit{lc}) 
if $(X, \Delta)$ is  a log pair such that 
$a (E, X, \Delta) > -1$ (resp.\ $a (E, X, \Delta) \geq -1$)  for any prime divisor $E$ over $X$. 
Note that, in contrast to \cite{KM98} and \cite{kollar13}, 
we always assume $\Delta$ to be effective for a klt (or lc) pair $(X, \Delta)$. 
We say that $X$ is \textit{klt} (resp.\ \textit{lc}) if so is $(X, 0)$. 



\subsection{Log del Pezzo pairs}

\begin{definition}
{We say that}
\begin{enumerate}
\item $(X, \Delta)$ is a {\em log del Pezzo pair} if 
$(X, \Delta)$ is a two-dimensional projective klt pair such that $-(K_X+\Delta)$ is ample; 
\item $(X, \Delta)$ is a {\em weak log del Pezzo pair} if 
$(X, \Delta)$ is a two-dimensional projective klt pair such that $-(K_X+\Delta)$ is nef and big;
\item $X$ is {\em of del Pezzo type} if 
$X$ is a projective normal surface and there exists an effective $\Q$-divisor $\Delta$ such that 
$(X, \Delta)$ is a log del Pezzo pair. 
\end{enumerate}
\end{definition}

\begin{remark}\label{r-log-dP-summary}
{The following properties of del Pezzo pairs will be used later on.}
\begin{enumerate}
\item 
Let $(X, \Delta)$ be a weak log del Pezzo pair. 
Since $-(K_X+\Delta)$ is big, we can write $-(K_X+\Delta) = A+E$ for some ample $\Q$-divisor $A$ and an effective $\Q$-divisor $E$. 
Then $(X, \Delta+ \epsilon E)$ is a log del Pezzo pair for $0 < \epsilon \ll 1$. 
Indeed, $(X, \Delta+ \epsilon E)$ is klt and 
\[
-(K_X+\Delta +\epsilon E)=-(1-\epsilon)(K_X+\Delta) + \epsilon A
\]
is ample.
\item
We have the following implications: 
\begin{align*}
(X, \Delta) \text{ is a log del Pezzo pair} 
&\Rightarrow 
(X, \Delta) \text{ is a weak log del Pezzo pair} \\
&\Rightarrow 
X \text{ is of del Pezzo type.}
\end{align*}
Indeed, the first implication is clear and the second one follows from (1). 
\item 
Let $X$ be a surface of del Pezzo type. 
Fix an effective  $\Q$-divisor $\Delta$ on $X$ such that $(X, \Delta)$ is klt 
and $-(K_X+\Delta)$ is ample. 
Then the following hold. 
\begin{enumerate}
 \item $X$ is $\Q$-factorial \cite[Theorem 5.4]{tanaka12}. 
 \item 
 For any $\Q$-divisor $D$, we may run a $D$-MMP \cite[Theorem 1.1]{tanaka12}. 
 Indeed, we have 
 \begin{equation}\label{e1-log-dP-summary}
 \epsilon D = K_X+\Delta -(K_X+\Delta) +\epsilon D \sim_{\Q} K_X+\Delta'
 \end{equation}
 for some $0 < \epsilon \ll 1$ and effective $\Q$-divisor $\Delta'$ such that $(X, \Delta')$ is klt. 
 \item 
 If $D$ is a nef $\Q$-divisor, then $D$ is semi-ample by (\ref{e1-log-dP-summary}) 
 and \cite[Theorem 1.2]{tanaka12}. 
 \item 
 It holds that 
 \[
 \NE(X) = \overline{\NE}(X) = \sum_{i=1}^n R_i
 \]
 for finitely many extremal rays $R_1, ..., R_n$ of $\NE(X)$  \cite[Theorem 3.13]{tanaka12}. 
 For each $R_i$, there exists $f_i : X \to Y_i$, called the contraction of $R_i$, 
 which is a morphism to a projective normal variety $Y_i$ 
 such that $(f_i)_*\MO_X = \MO_{Y_i}$ and, given a curve $C$ on $X$, 
 $f_i(C)$ is a point if and only if $[C] \in R_i$ \cite[Theorem 3.21]{tanaka12}. 
\end{enumerate}
\end{enumerate}
\end{remark}

\subsection{Quasi-$F$-splitting}\label{ss-def-QFS}
 
In this subsection, we summarise the definition and some basic properties of a quasi-$F$-splitting. 
For details, we refer to \cite{KTTWYY1}, in which we work in a more general setting. 

\begin{definition}\label{d-QFS}
Let $X$ be a normal variety and let $\Delta$ be an effective $\Q$-divisor on $X$ satisfying $\rdown{\Delta} =0$. 
For $n \in \Z_{>0}$, we say that $(X, \Delta)$ is $n$-{\em quasi-$F$-split} if 
there exists 
a $W_n\MO_X$-module homomorphism $\alpha : F_* W_n\MO_X(p\Delta) \to \MO_X$ which completes the following commutative diagram: 
\begin{equation*} \label{diagram:intro-definition}
\begin{tikzcd}
W_n\cO_X(\Delta) \arrow{r}{{F}} \arrow{d}{R^{n-1}} & F_* W_n \arrow[dashed]{ld}{\exists\alpha} \cO_X(p\Delta) \\
\MO_X. 
\end{tikzcd}
\end{equation*}
We say that $(X, \Delta)$ is {\em quasi-$F$-split} if there exists $n \in \Z_{>0}$ such that $(X, \Delta)$ is $n$-quasi-$F$-split. 
\end{definition}

\begin{remark}\label{r-QFS}
Let $X$ be a normal variety and let $\Delta$ be an effective $\Q$-divisor on $X$ satisfying $\rdown{\Delta} =0$. 
Fix $n \in \Z_{>0}$. 
We then define $Q_{X, \Delta, n}$ by the following pushout diagram of $W_n\MO_X$-modules: 
\[
\begin{tikzcd}
W_n\cO_X(\Delta) \arrow{r}{{F}} \arrow{d}{R^{n-1}} & F_* W_n\cO_X(p\Delta)  \arrow{d} \\
\MO_X \arrow{r}{\Phi_{X, \Delta, n}} & Q_{X, \Delta, n}. 
\end{tikzcd}
\]
\begin{enumerate}
\item  
The $W_n\MO_X$-module $Q_{X, \Delta, n}$ is naturally a coherent $\MO_X$-module, and hence 
$\Phi_{X, \Delta, n} : \MO_X \to Q_{X, \Delta, n}$ is an $\MO_X$-module homomorphism 
\cite[Proposition 3.6(2)]{KTTWYY1}. 
In particular, the following {conditions} are equivalent \cite[Proposition 3.7]{KTTWYY1}. 
\begin{itemize}
\item  $(X, \Delta)$ is $n$-quasi-$F$-split. 
\item  
$\Phi_{X, \Delta, n} : \MO_X \to Q_{X, \Delta, n}$ splits as an $\MO_X$-module homomorphism.
\end{itemize}
\item  
Assume that $X$ is projective. 
Then 
$(X, \Delta)$ is $n$-quasi-$F$-split if and only if 
\begin{equation} \label{eq:coh-def-of-qfsplit}
H^{\dim X}(\Phi_{X, K_X+\Delta, n}) : 
H^{\dim X}(X,  
 \MO_X(K_X+\Delta)) \to H^{\dim X}(X, Q_{X, K_X+\Delta, n})
\end{equation}
is injective \cite[Lemma 3.13]{KTTWYY1}. 
\end{enumerate}
\end{remark}

\begin{theorem}\label{t-QFS-criterion}
Let $(X, \Delta)$ be a log del Pezzo pair. 
Assume that there exists a log resolution $f: Y \to X$ of $(X, \Delta)$ and a $\Q$-divisor $B_Y$ on $Y$ such that the following conditions hold: 
\begin{enumerate}
\item[(A)] $\rdown{B_Y} \leq 0$, $f_*B_Y  =\Delta$, and $-(K_Y+B_Y)$ is ample, 
\item[(B)] $H^0(Y, \MO_Y(K_Y + E + p^{\ell}(K_Y+B_Y))) =0$ 
for every $\ell \in \Z_{>0}$, 
\item[(C)] $H^0(Y, \Omega_Y^1(\log E) \otimes \MO_Y(K_Y+B_Y))=0$,
\end{enumerate}
where $E := (B_Y)_{\red}$. Then $(X, \Delta)$ is quasi-$F$-split. 
\end{theorem}

\begin{proof}
The assertion holds by  \cite[Theorem 5.13]{KTTWYY1} 
and the following inclusion: 
\[
H^0(Y, B_1\Omega_Y^2(\log E)(p^{\ell}(K_Y+B_Y))) 
\subseteq 
H^0(Y, \MO_Y(K_Y + E + p^{\ell}(K_Y+B_Y))) =0. 
\]
\end{proof}


\subsection{Higher log Cartier operator for $\Q$-divisors}

For definitions and fundamental properties of higher log Cartier operators, 
we refer to \cite[Subsection 5.2]{KTTWYY1}.

\begin{lemma}\label{l-Cartier-op}
Let $(Y,E)$ be a log smooth pair, where $E$ is a reduced divisor. 
Let \vspace{0.05em} $D$ be a $\Q$-divisor satisfying $\Supp\, \{D\}\subseteq E$. 
Then we have the following exact \vspace{0.1em} sequences for all $i \in \Z_{\geq 0}$ and $m, r\in\Z_{>0}$ satisfying $r<m$.
\begin{smallertags}
{\small \begin{align}
&\!\!0 \to B_m\Omega^{i}_Y(\log\,E)(p^mD)\to Z_m\Omega^{i}_Y(\log\,E)(p^mD)
\overset{{C^{m}}}{\to} 
\Omega^{i}_Y(\log\,E)(D)\to 0,\label{exact1}\\[0.5em]
&\!\!0 \to F_{*}^{m{-}r}\!B_{r}\Omega^{i}_Y(\log E)(p^{m}\!D) \!\to\! B_{m}\Omega^{i}_Y(\log E)(p^m\!D) \overset{{C^{r}}}{\to} B_{m-r}\Omega^{i}_Y(\log E)(p^{m-r}\!D) \!\to\! 0,\label{exact2}\\[0.5em]
&\!\!0 \to F_{*}^{m{-}r}B_{r}\Omega^{i}_Y(\log E)(p^{m}\!D) \!\to\! Z_{m}\Omega^{i}_Y(\log E)(p^m\!D) \overset{{C^{r}}}{\to} Z_{m-r}\Omega^{i}_Y(\log E)(p^{m-r}\!D) \!\to\! 0,\label{exact3}\\[0.5em]
&\!\!0 \to Z_m\Omega^{i}_Y(\log\,E)(p^{m}D) \to F_{*}Z_{m-1}\Omega^{i}_Y(\log\,E)(p^{m}D) 
\xrightarrow{\psi}
B_1\Omega^{i+1}_Y(\log\,E)(pD) \to 0\label{exact4},
\end{align}}
\end{smallertags}
\!\!where $\psi := F_{*}d\circ C^{m{-}1}$.
\end{lemma}
\begin{proof}
Here, (\ref{exact1}) is exact by \cite[(5.7.1)]{KTTWYY1}. 
The exactness of (\ref{exact2}) and (\ref{exact3}) follow from the constructions of $B_m\Omega^{i}_Y(\log\,E)(p^mD)$ and $Z_m\Omega^{i}_Y(\log\,E)(p^mD)$ (cf.\  \cite[Definition 5.6]{KTTWYY1}).  
As for the last one (\ref{exact4}), see \cite[Lemma 5.8]{KTTWYY1}.
\end{proof}

\begin{lemma}\label{lem:Serre's map}
Let $(Y,E)$ be a log smooth pair, where $E$ is a reduced divisor. 
Let $D$ be a $\Q$-divisor satisfying $\Supp \{D\}\subseteq E$.
Then  there exists the following exact sequence of $W_m\MO_Y$-modules 
\[
0\to W_m\MO_Y(D) \xrightarrow{F} F_{*}W_m\MO_Y(pD) \xrightarrow{s} B_m\Omega_Y^{1}(\log\,E)(p^mD)\to 0, 
\]
where $s$ is defined by 
\[
s(F_*(f_0,\ldots,f_{n-1}))\coloneqq 
F_*^n(f_{0}^{p^{n-1}-1}df_{0} + f_1^{p^{n-2}-1}df_1+\ldots+df_{n-1}).
\]
Furthermore, ${\rm Coker}(\Phi_{Y, D, m} : \MO_Y(D) \to Q_{Y, D, m})  \simeq B_m\Omega_Y^{1}(\log\,E)(p^mD)$. 
\end{lemma}

\begin{proof}
See \cite[Lemma 5.9]{KTTWYY1}. 
\end{proof}

%% file: section3.tex
\section{Log resolutions of klt surfaces with effective boundaries}\label{s-klt-resol}


Throughout this section, we work over an algebraically closed field $k$. 
The main question of this section is as follows.

\begin{question}\label{q-klt-resol}
Let $(X, \Delta)$ be a two-dimensional klt pair with standard coefficients. 
Does there exist a log resolution $\varphi : V \to X$ of $(X, \Delta)$ 
such that $\Delta_V$ is effective for the $\Q$-divisor $\Delta_V$ defined by $K_V+\Delta_V = \varphi^*(K_X+\Delta)$?  
\end{question}

If $\Delta=0$, then it is well known that the answer is affirmative, 
since the minimal resolution of $X$ is a log resolution of $(X, \Delta=0)$ \cite[Theorem 4.7]{KM98}. 
Unfortunately, the answer is negative in general. 
For example, we can not find such a log resolution when $(X, \Delta) = (\mathbb A^2, \frac{1}{2} C)$, where $C$ is a cusp $\{y^2 = x^{2n+1}\}$, and $n \in \Z_{>0}$. 
The purpose of this section is to prove that 
if $(X,\Delta)$ is a counterexample to the above question, then up to localising at some point $x\in X$ we have that $X$ is smooth, $\Delta=\frac{1}{2}C$ for some prime divisor $C$, and $(X,\Delta)$ has the same dual graph as the example $(\mathbb A^2, \frac{1}{2} C)$ as above. 
 (cf.\ Theorem \ref{t-hiromu-resol}, Figure \ref{figure-V}). 

In Subsection \ref{ss-klt-resol-sm}, 
we treat the case when $X$ is smooth and the answer to Question \ref{q-klt-resol} is affirmative. 
In Subsection \ref{ss-klt-resol-sing}, 
we prove that Question \ref{q-klt-resol} is affirmative 
for the case when $X$ is singular. 
Based on these subsections, we shall prove the main theorem of this subsection in Subsection \ref{ss-klt-resol-general}.

\subsection{The smooth case}\label{ss-klt-resol-sm}

In this subsection, we consider Question \ref{q-klt-resol} 
for the case when $X$ is smooth. 
The main objective is to show that the answer to Question \ref{q-klt-resol} is affirmative for many cases (Lemma \ref{l-klt-resol-r=3}, Lemma \ref{l-klt-resol-r=2}, Lemma \ref{l-klt-resol-r=1}). 
We start by introducing the following notation.

\begin{notation}\label{n-klt-resol-sm}
Let $(X, \Delta)$ be a two-dimensional klt pair such that $X$ is smooth. 
Assume that 
there exists a closed point $x$ of $X$ such that 
$\Delta_{\red}|_{X \setminus x}$ is smooth 
and $\Delta$ is not simple normal crossing at $x$. 
Let $\Delta = \sum_{i=1}^r c_i C_i$ be the irreducible decomposition, 
where $c_i \geq 1/2$ holds for any $1 \leq i \leq r$. 
Assume that all of $C_1, ..., C_r$ pass through $x$. 
\end{notation}

\begin{lemma}\label{l-klt-resol-r=3}
We use Notation \ref{n-klt-resol-sm}. 
Assume $r\geq 3$. 
Then $r=3$ and there exists a log resolution $\varphi : V \to X$ of $(X, \Delta)$ such that 
$\varphi(\Ex(\varphi)) = \{x\}$ and $\Delta_V$ is effective 
for the $\Q$-divisor $\Delta_V$ defined by $K_V + \Delta_V = \varphi^*(K_X+\Delta)$. 
\end{lemma}

\begin{proof}
Let $f: Y \to X$ be the blowup at $x$. Set $E:=\Ex(f)$. We have 
\[
K_Y + \Delta_Y = K_Y + f_*^{-1}\Delta + bE = f^*(K_X+\Delta)\quad \text{for}\quad \Delta_Y := f_*^{-1}\Delta + bE
\]
\[
\text{and}\quad 
b := -1 + \mult_x \Delta = -1 + \sum_{i=1}^r c_i \mult_x C_i \geq  -1 + \frac{1}{2}\sum_{i=1}^r \mult_x C_i. 
\]
If $r \geq 4$, then we obtain $b \geq 1$, which contradicts the assumption that $(X, \Delta)$ is klt. 
Hence $r=3$. 
By decreasing the coefficients of $\Delta,$ we may assume that all the coefficients of $\Delta$ are $\frac{1}{2}$, that is,  $\Delta = \frac{1}{2}(C_1+C_2+C_3)$. 
Then 
\[
\mult_x \Delta = \frac{3}{2}\qquad\text{and}\qquad 
\mult_x C_1 = \mult_x C_2 = \mult_x C_3 = 1. 
\]
Hence each $C_i$ is smooth at $x$. In particular,  we obtain 
\[
\Delta_Y = \frac{1}{2}(E+C_{1, Y} +  C_{2, Y} + C_{3, Y}) \qquad 
{\text{for}\qquad C_{i, Y} := f_*^{-1}C_i}. 
\]
After permuting $C_1, C_2, C_3$ if necessary, we have the following three cases. 
\begin{enumerate}
    \item Each of $C_1+C_2, C_2+C_3, C_3+C_1$ is simple normal crossing. 
    \item None of $C_1+C_2, C_2+C_3, C_3+C_1$ is simple normal crossing. 
    \item $C_1+C_2$ is not simple normal crossing, but $C_1+C_3$ is simple normal crossing. 
\end{enumerate}
For $1 \leq i < j \leq 3$, 
$C_i + C_j$ is simple normal crossing at $x$ if and only if $C_{i, Y} \cap C_{j, Y} = \emptyset$. 

Assume (1). 
Then the blowup $Y \to X$ is a log resolution of $(X, \Delta)$, 
because $C_{i, Y} \cap C_{j, Y} = \emptyset$ for $1 \leq i < j \leq 3$. 
We are done. 

Assume (2). 
We have $C_{1, Y} \cap E = C_{2, Y} \cap E = C_{3, Y} \cap E =: y$, which is a point. 
Then it holds that 
\[
\mult_y \Delta_Y  = \frac{1}{2}\mult_y (E+C_{1, Y} +  C_{2, Y} + C_{3, Y})=2. 
\]
However, this implies that $(Y, \Delta_Y)$ is not klt, which contradicts the assumption that $(X, \Delta)$   is klt. Hence the case (2) does not occur. 

Assume (3).  In this case, we have 
\[
y := C_{1, Y} \cap E = C_{2, Y} \cap E \neq C_{3, Y} \cap E. 
\]
Note that $\Delta_Y$ is simple normal crossing around $E \cap C_{3, Y}$. 
On the other hand, three curves $E, C_{1, Y}, C_{2, Y}$ intersects at a single point $y$. 
We again take the blowup $g: Z \to Y$ at $y$. 
Then the same argument as above can be applied after replacing $(Y, \Delta_Y, y)$ by $(X, \Delta, x)$, 
{\cred so that} we can repeat the same procedure. 
{\cred This procedure will terminate after finitely many times, because the intersection multiplicity strictly drops: 
$(C_1 \cdot C_2)_x >  (C_{1, Y} \cdot C_{2, Y})_y$ (for the definition of intersection multiplicities, see \cite[Ch. V, Section 1]{hartshorne77}). 
{Note that this inequality $(C_1 \cdot C_2)_x >  (C_{1, Y} \cdot C_{2, Y})_y$ is proven as follows: by taking a compactification of $X$ such that $(C_1 \cap C_2)_{\red} = \{x\}$, we may assume that 
$X$ is proper and $C_1 \cdot C_2 = (C_1 \cdot C_2)_x$; then 
$(C_1 \cdot C_2)_x = C_1 \cdot C_2 = C_{1, Y} \cdot C_{2, Y} +1 >C_{1, Y}\cdot C_{2, Y}  \geq (C_{1, Y} \cdot C_{2, Y})_y$.}} 
As an alternative proof, 
this termination is assured by the fact that 
a log resolution is obtained by taking successive blowups of points.  
\qedhere
\end{proof}

\begin{lemma}\label{l-klt-resol-r=2}
We use Notation \ref{n-klt-resol-sm}. 
Assume $r=2$. 
Then there exists a log resolution $\varphi: V \to X$ of $(X, \Delta)$ such that 
$\varphi(\Ex(\varphi)) = \{x\}$ and $\Delta_V$ is effective 
for the $\Q$-divisor $\Delta_V$ defined by $K_V + \Delta_V = \varphi^*(K_X+\Delta)$. 
\end{lemma}
\begin{proof}
By decreasing the coefficients of $\Delta$, we may assume that all the coefficients of $\Delta$ are equal to $\frac{1}{2}$, 
that is, $\Delta = \frac{1}{2}(C_1 +C_2)$.  
Since $(X, \Delta)$ is klt, we have $\mult_x \Delta <2$. 
As each $C_i$ passes through $x$, we get $\mult_x \Delta \geq 1$. 
By $\mult_x \Delta \in \frac{1}{2} \Z$, we have the following two cases. 
\begin{enumerate}
    \item[(I)] $\mult_x \Delta = 1$. 
    \item[(II)] $\mult_x \Delta =3/2$. 
\end{enumerate}

(I) 
Assume $\mult_x \Delta =1$. 
Then $\mult_x C_1 = \mult_x C_2 = 1$, that is, both $C_1$ and $C_2$ are smooth at $x$. 
Let $f: Y \to X$ be the blowup at $x$. 
Then 
\[
K_Y + \Delta_Y = f^*(K_X+\Delta) \qquad \text{for}\qquad 
\]
\[
\Delta_Y = \frac{1}{2}(C_{1, Y} + C_{2, Y}), \qquad C_{1, Y} := f_*^{-1}C_1,\qquad C_{2, Y} := f_*^{-1}C_2. 
\]
For any point $y \in C_{1, Y} \cap C_{2, Y}$, 
we have an inquality of intersection multiplicities 
$(C_{1, Y} \cdot C_{2, Y})_y \leq (C_1 \cdot C_2)_x -1$ 
(e.g. if $X$ is projective, then $C_{1, Y} \cdot C_{2, Y} = C_1 \cdot C_2 -1$). 
We take blowups until (the proper transforms of) $C_1$ and $C_2$ will be disjoint. 
Then we obtain 
\[
K_{V} + \Delta_V = \varphi^*(K_X+\Delta)
\]
for $\Delta_V = \frac{1}{2}(C_{1, V} + C_{2, V}), 
C_{1, V} := \varphi_*^{-1}C_1, 
C_{2, V} := \varphi_*^{-1}C_2$, that is,  
the coefficient of any $\varphi$-exceptional prime divisor is zero. 
Furthermore, it is easy to see that 
$\varphi: V \to X$ is a log resolution of $(X, \Delta)$.

(II) Assume   $\mult_x \Delta =3/2$. Then 
\[
K_Y + \frac{1}{2}(E+ C_{1, Y} + C_{2, Y}) = f^*(K_X+\Delta). 
\]
After permuting $C_1$ and $C_2$ if necessary, 
we may assume that $\mult_x C_1 =1$ and $\mult_x C_2 =2$. 
Hence, $C_1$ is smooth at $x$ and $C_2$ is not smooth at $x$. 
We have $C_{1, Y} \cdot E =1$ and $C_{2, Y} \cdot E =2$. 
In particular, $C_{1, Y} + E$ is simple normal crossing and $y_1 := C_{1, Y} \cap E$ is a point. 


We now treat the case when $C_{2, Y}+E$ is simple normal crossing, that is,  
$C_{2, Y} \cap E = \{y_2, y'_2\}$ with $y_2 \neq y'_2$. 
If $y_1 \neq y_2$ and $y_1 \neq y'_2$, then $\Delta_Y$ is simple normal crossing. 
We are done. 
Then we may assume that $y_1 = y_2$. 
In this case, we can apply 
{Lemma \ref{l-klt-resol-r=3}.}
This completes the case when $C_{2, Y}+E$ is simple normal crossing. 

We may assume that $C_{2, Y}+E$ is not simple normal crossing, that is,  
$y_2 :=(C_{2, Y} \cap E)_{\red}$ is one point. 
If $y_1 = y_2$, then we can apply 
{Lemma \ref{l-klt-resol-r=3}.}
The problem is reduced to the case when $y_1 \neq y_2$. 
Then $(Y, \Delta_Y)$ is simple normal crossing at $y_1$. 
If $\mult_{y_2}\Delta_Y = 1$, then we are done by (I). 
We obtain $\mult_{y_2}\Delta_Y = 3/2$. 
After replacing $(Y, \Delta_Y, y_2)$ by $(X, \Delta, x)$, 
we may apply the same argument as above. 
In other words, we again take the blowup $g: Z \to Y$ at $y_2$. 
Applying this argument repeatedly, 
we obtain a log resolution $\varphi : V \to X$ such that $\Delta_V$ is effective 
for $K_V + \Delta_V = \varphi^*(K_X+\Delta)$. 
\end{proof}

\begin{lemma}\label{l-klt-resol-r=1}
We use Notation \ref{n-klt-resol-sm}. 
Assume $\mult_x \Delta \geq \frac{3}{2}$. 
Then there exists a log resolution $\varphi : V \to X$ of $(X, \Delta)$ such that 
$\varphi(\Ex(\varphi)) = \{x\}$ and $\Delta_V$ is effective 
for the $\Q$-divisor $\Delta_V$ defined by $K_V + \Delta_V = \varphi^*(K_X+\Delta)$. 
\end{lemma}

\begin{proof}
Let $f: Y \to X$ be the blowup at $x$. 
We obtain 
\[
K_Y + \Delta_Y =f^*(K_X+\Delta)
\]
\[
\text{for}\qquad E := \Ex(f)\qquad\text{and}\qquad \Delta_Y := f_*^{-1}\Delta + (\mult_x \Delta -1) E. 
\]
By $\mult_x \Delta \geq \frac{3}{2}$, any coefficient of $\Delta_Y$ is $\geq \frac{1}{2}$. 
Then any non-log smooth point of $(Y, \Delta_Y)$ is contained in at least two irreducible components of $\Delta_Y$. 
Therefore, we can find a log resolution $\varphi_Y: V \to Y$ of $(Y, \Delta_Y)$ 
such that $\Delta_V$ is effective for $K_V + \Delta_V = \varphi_Y^*(K_Y+\Delta_Y)$ 
by Lemma \ref{l-klt-resol-r=3} and Lemma \ref{l-klt-resol-r=2}. 
Then the composition $\varphi : V \xrightarrow{\varphi_Y} Y \xrightarrow{f} X$ is a required log resolution 
of $(X, \Delta)$. 
\end{proof}

\subsection{The singular case}\label{ss-klt-resol-sing}

The purpose of this subsection is to show that 
Question \ref{q-klt-resol} is affirmative for the case when $X$ is singular and $(X, \Delta)$ is log smooth outside the singular points of $X$ (Proposition \ref{p-klt-sing}). 
We first treat two typical cases in 
 Lemma \ref{l-klt-sing1} and 
Lemma \ref{l-klt-sing2}. 
The general case will be reduced to these cases (cf.\ the proof of Proposition \ref{p-klt-sing}).

\begin{lemma}\label{l-klt-sing1}
Let $(X, \Delta)$ be a two-dimensional klt pair such that all the coefficients of $\Delta$ are 
$\geq \frac{1}{2}$. 
Let $x$ be a singular point of $X$ 
such that $(X \setminus \{x\}, \Delta|_{X \setminus \{x\}})$ is log smooth. 
For the minimal resolution $f: Y \to X$ of $X$ and 
the $\Q$-divisor $\Delta_Y$ defined by $K_Y+\Delta_Y = f^*(K_X+\Delta)$, assume that $E:= \Ex(f)$ is a prime divisor and $ \Delta_Y$ is not simple normal crossing. 
Then $E \subseteq \Supp \Delta_Y$ and 
the coefficient of $E$ in $\Delta_Y$ is $\geq \frac{1}{2}$. 
\end{lemma}

\begin{proof}
If $x \not \in \Supp \Delta$, then $\Delta_Y$ must be simple normal crossing which contradicts our assumption. 
Thus $x \in \Supp \Delta$, and so $\Supp \Delta_Y = E \cup \Supp f_*^{-1}\Delta$. By decreasing the coefficients of $\Delta,$ we may assume that every coefficient of $\Delta$ is equal to $\frac{1}{2}$. 
Set $m:=-E^2 \geq 2$. 
We have 
\[
K_Y + f_*^{-1}\Delta + E =f^*(K_X+\Delta) + bE 
\]
for the rational number $b$ satisfying the following: 
\[
b = \frac{1}{E^2} ( -2 + f_*^{-1}\Delta \cdot E) = \frac{1}{m}
(2 -  f_*^{-1}\Delta \cdot E). 
\]
Since $2f_*^{-1}\Delta$ is Cartier and 
$\Supp \Delta_Y = \Supp (f_*^{-1} \Delta + E)$ is not simple normal crossing, 
we obtain $f_*^{-1}\Delta \cdot E \geq 1$. 
Hence 
\[
b  = \frac{1}{m}
(2 -  f_*^{-1}\Delta \cdot E) \leq  \frac{1}{m}
(2 -  1) \leq \frac{1}{2}.  
\]
Therefore, any coefficient of $\Delta_Y =  f_*^{-1}\Delta + (1-b)E$ is $\geq \frac{1}{2}$. 
\end{proof}

\begin{lemma}\label{l-klt-sing2}
Let $(X, \Delta)$ be a two-dimensional klt pair such that all the coefficients of $\Delta$ are $\geq \frac{1}{2}$. 
Let $x$ be a singular point of $X$ 
such that $(X \setminus \{x\}, \Delta|_{X \setminus \{x\}})$ is  log smooth.  
For the minimal resolution $f: Y \to X$ of $X$, assume that the following hold. 
\begin{enumerate}
    \item 
    $\Ex(f) = E_1 \cup E_2$ is the irreducible decomposition, 
    where $E_1$ and $E_2$ are distinct prime divisors on $Y$. 
    \item 
    $f_*^{-1}\Delta$ contains $E_1 \cap E_2$. 
\end{enumerate}
Let $\Delta_Y$ be the $\Q$-divisor defined by $K_Y+ \Delta_Y = f^*(K_X+\Delta)$. 
Then $\Ex(f) \subseteq \Supp \Delta_Y$ and 
all the coefficients of $\Delta_Y$ are $\geq \frac{1}{2}$. 
\end{lemma}

\begin{proof}
Note that $E_1 \cdot E_2 =1$. 
In other words, $E_1 +E_2$ is simple normal crossing and $y := E_1 \cap E_2$ is a closed point of $Y$. 
Set $m_1 := -E_1^2$ and $m_2 := -E_2^2$. Of course, we have 
\begin{equation}\label{e1-klt-sing1}
m_1 \geq 2\qquad \text{and}\qquad m_2 \geq 2. 
\end{equation}
By (2), 
we obtain $x \in \Supp \Delta$ and $\Supp \Delta_Y = \Exc(f) \cup \Supp f_*^{-1}\Delta$. 
By decreasing the coefficients of $\Delta,$ we may assume that 
any coefficient of $\Delta$ is equal to $\frac{1}{2}$. 
Let $b_1, b_2 \in \Q$  be the rational numbers defined by 
\[
K_Y + f_*^{-1}\Delta + E_1+E_2 =f^*(K_X+\Delta) + b_1 E_1+b_2E_2.  
\]
Since $(X, \Delta)$ is klt, we have $b_1 >0$ and $b_2 >0$. 
We  obtain 
\[
\Delta_ Y = f_*^{-1}\Delta + (1-b_1)E_1+ (1-b_2)E_2. 
\]
It follows from (2) that 
\begin{equation}\label{e2-klt-sing1}
f_*^{-1}\Delta \cdot E_1 >0 \qquad \text{and} \qquad f_*^{-1}\Delta \cdot E_2 >0. 
\end{equation}
For the intersection matrix 
\[
A := (E_i \cdot E_j) = 
\begin{pmatrix}
-m_1 & 1\\
1 & -m_2 
\end{pmatrix}, 
\]
its inverse matrix $A^{-1}$ is given as follows: 
\[
-A^{-1} = \frac{-1}{m_1m_2-1}
\begin{pmatrix}
-m_2 & -1\\
-1 & -m_1
\end{pmatrix}
=\frac{1}{m_1m_2-1}
\begin{pmatrix}
m_2 & 1\\
1 & m_1
\end{pmatrix}. 
\]
Therefore, 
\[
\begin{pmatrix}
b_1\\
b_2
\end{pmatrix}
=
-A^{-1}
\begin{pmatrix}
2 - (f_*^{-1}\Delta + E_2) \cdot E_1\\
2 - (f_*^{-1}\Delta + E_1) \cdot E_2
\end{pmatrix} = 
\frac{1}{m_1m_2-1}
\begin{pmatrix}
m_2 & 1\\
1 & m_1
\end{pmatrix}
\begin{pmatrix}
1 - f_*^{-1}\Delta \cdot E_1\\
1 - f_*^{-1}\Delta \cdot E_2
\end{pmatrix}. 
\]
By (\ref{e2-klt-sing1}) and $f_*^{-1}\Delta \cdot E_i \in \frac{1}{2}\Z$, 
it holds that 
\begin{equation}\label{e3-klt-sing1}
1 - f_*^{-1}\Delta \cdot E_1 \leq \frac{1}{2} \qquad \text{and} \qquad 
1 - f_*^{-1}\Delta \cdot E_2 \leq \frac{1}{2}. 
\end{equation}

By (\ref{e1-klt-sing1}) and (\ref{e3-klt-sing1}), the following holds: 
\begin{align*}
b_1 &= \frac{1}{m_1m_2-1}( m_2(1 - f_*^{-1}\Delta \cdot E_1) + (1 - f_*^{-1}\Delta \cdot E_2))\\
&\leq \frac{m_2 +1}{2(m_1m_2-1)} \leq \frac{m_2 +1}{2(2m_2-1)} 
\leq 
\frac{1}{2}, 
\end{align*}
where the last inequality holds by 
$\frac{m_2 +1}{2m_2-1}  = \frac{1}{2}+ \frac{3/2}{2m_2-1}\leq 
  \frac{1}{2}+ \frac{3/2}{2 \cdot 2-1} =1$. 
By symmetry, we get $b_2 \leq \frac{1}{2}$. 
Hence any coefficient of $\Delta_Y = f_*^{-1}\Delta + (1-b_1)E_1+ (1-b_2)E_2$ is $\geq 1/2$. 
\end{proof}

\begin{proposition}\label{p-klt-sing}
Let $(X, \Delta)$ be a two-dimensional klt pair such that all the coefficients of $\Delta$ are 
$\geq \frac{1}{2}$. 
Let $x$ be a singular point of $X$ 
such that $(X \setminus \{x\}, \Delta|_{X \setminus \{x\}})$ is log smooth. 
Then there exists a log resolution $\varphi : V \to X$ of $(X, \Delta)$ such that $\varphi(\Ex(\varphi)) = \{x\}$ and  $\Delta_V$ is effective 
for the $\Q$-divisor $\Delta_V$ defined by $K_V + \Delta_V = \varphi^*(K_X+\Delta)$. 
\end{proposition}

\begin{proof}
Let 
\[
f : Y \to X
\]
be the minimal resolution of $X$. 
Let $\Ex(f) = E_1 \cup \cdots \cup E_n$ be the irreducible decomposition. 
We define an effective $\Q$-divisor $\Delta_Y$ on $Y$ by 
\[
K_Y + \Delta_Y = f^*(K_X+\Delta).
\]
If $x \not\in \Supp \Delta$, then 
it is enough to set $V:=Y$. 
In what follows, we assume that $x \in \Supp \Delta$. 
In particular, 
\begin{equation}\label{e1-klt-sing}
\Supp \Delta_Y = \Ex(f) \cup \Supp f_*^{-1}\Delta = 
E_1 \cup \cdots \cup E_n \cup \Supp f_*^{-1}\Delta. 
\end{equation}
We emphasise that $\Ex(f)$ is simple normal crossing (but $\Supp \Delta_Y$ need not be simple normal crossing).

\begin{claim}\label{c1-klt-sing}
Fix $1 \leq i \leq n$. 
If there exists a closed point $y \in E_i$ 
at which $\Delta_Y$ is not simple normal crossing, 
then the coefficient of $E_i$ in $\Delta_Y$ is $\geq \frac{1}{2}$. 
\end{claim}

\begin{proof}[Proof of Claim \ref{c1-klt-sing}]
Since $E_1 \cup \cdots \cup E_n$ is simple normal crossing, there are the following two possible cases:\ (1) and (2). 
\begin{enumerate}
    \item $y \not\in E_{j}$ for every $j$ satisfying $j \neq i$. 
    \item There exists $1 \leq i' \leq n$ 
    such that $i' \neq i$, $y \in E_{i} \cap E_{i'}$, 
    and $y \not\in E_j$ for any $j \neq i, i'$.
\end{enumerate}

Assume (1).  
Let $f' : Y \to X'$ be the projective birational morphism 
to a normal surface $X'$ such that $\Ex(f')=E_i$. 
Then $f'$ factors through $f$: 
\[
f : Y \xrightarrow{f'} X' \xrightarrow{g} X. 
\]
By $E_i^2 \leq -2$, 
$f': Y \to X'$ is the minimal resolution of $X'$. 
Set $\Delta' := g_*^{-1} \Delta$. 
We have 
\[
K_{X'} + \Delta' \leq g^*(K_X+\Delta).
\]
We define a $\Q$-divisor $\Delta'_Y$ by 
$K_Y+\Delta'_Y = f'^*(K_{X'}+\Delta')$. 
Then $\Delta'_Y$ is effective, since $f' : Y \to X'$ is the minimal resolution of $X'$. 
It holds that 
\[
K_Y + \Delta'_Y = f'^*(K_{X'}+\Delta')\leq 
f'^*(g^*(K_X+\Delta)) = f^*(K_X+\Delta) = K_Y+\Delta_Y.
\]
Hence it is enough to show that the coefficient 
of $E_i$ in $\Delta'_Y$ is $\geq \frac{1}{2}$. 
For $x' := f'(E_i)$ and a suitable open neighbourhood ${\widetilde X'}$ of $x' \in X'$, 
we may apply Lemma \ref{l-klt-sing1} for 
$({\widetilde X'}, \Delta'|_{{\widetilde X'}})$ and the minimal resolution 
$f'|_{f'^{-1}({ \widetilde X'})} : f'^{-1}({ \widetilde X'}) \to { \widetilde X'}$ of ${\widetilde X'}$, 
because $\Delta'_Y$ is not simple normal crossing at $y$. 
This completes the proof of Claim \ref{c1-klt-sing} for the case when (1) holds. 

Assume (2). 
Let $f' : Y \to X'$ be the projective birational morphism 
to a normal surface $X'$ such that $\Ex(f')=E_i \cup E_{i'}$. 
Then the same argument as that of (1) works by using Lemma \ref{l-klt-sing2} instead of Lemma \ref{l-klt-sing1}.  
This completes the proof of Claim \ref{c1-klt-sing}. 
\end{proof}

Fix a closed point $y \in Y$ at which $(Y, \Delta_Y)$ is not log smooth. 
By (\ref{e1-klt-sing}), 
it is enough to find an open neighbourhood $\widetilde Y$ of $y \in Y$ and a log resolution $\psi : \widetilde V \to \widetilde Y$ of $(\widetilde Y, \Delta_Y|_{\widetilde Y})$ 
such that $\Delta_{\widetilde V}$ is effective for the $\Q$-divisor defined by $K_{\widetilde V} + \Delta_{\widetilde V} =
\psi^*(K_{\widetilde Y}+ (\Delta_Y|_{\widetilde Y}))$. 
By Claim \ref{c1-klt-sing}, 
there exists an open neighbourhood $\widetilde Y$ of $y \in Y$  
such that all the coefficients of $\Delta_Y|_{\widetilde Y}$ are $\geq \frac{1}{2}$. 
Moreover, $y$ is contained in at least two irreducible components of $\Delta_Y|_{\widetilde Y}$: 
one of them is some $f$-exceptional prime divisor $E_i$ and 
we can find another one from $f_*^{-1}\Delta$. 
Therefore, we can find a required log resolution $\psi : \widetilde V \to \widetilde Y$ of $(\widetilde Y, \Delta_Y|_{\widetilde Y})$ by Lemma \ref{l-klt-resol-r=3} 
and Lemma \ref{l-klt-resol-r=2}. 
\qedhere
\end{proof}

\subsection{The general case}\label{ss-klt-resol-general}

We are ready to prove the main theorem of this section.

\begin{figure}
\centering
\begin{tikzpicture}
\draw (1.5, 2) node{$\frac{1}{2}C_V$};
\draw (-2.5, -2) node{$0E_n$};
\draw (-0.8, -1.1) node{$0E_{n-1}$};
\draw (1.7, -1.7) node{$0E_{n-2}$};
\draw (4.5, -1.3) node{$0E_{2}$};
\draw (7.5, -1.3) node{$0E_{1}$};

\draw (-2.4, 0) node{$-1$};  
\draw (-1.2, -1.8) node{$-2$}; 
\draw (1.5, -0.9) node{$-2$}; 
\draw (4.7, -2.0) node{$-2$}; 
\draw (7.0, -2.0) node{$-2$}; 

\draw[domain=-2:1, samples=500] plot (\x, {1+sqrt(0.333*\x+0.666)});
\draw[domain=-2:1, samples=500] plot (\x, {1-sqrt(0.333*\x+0.666)});

\draw(-2,2.2)--(-2,-2.2); 

\draw(-2.2,-0.95)--(0.2,-2.05);
\draw(-0.2,-2.05)--(2.2,-0.95);
\draw[shift={(3,-1.5)}] node{$\hdots$}; 
\draw(4,-2)--(6.2,-0.95);
\draw(5.8,-0.95)--(8,-2);

\end{tikzpicture}
\caption{The coefficients of $\Delta_V$ and self-intersection numbers on $V$}\label{figure-V}
\end{figure}
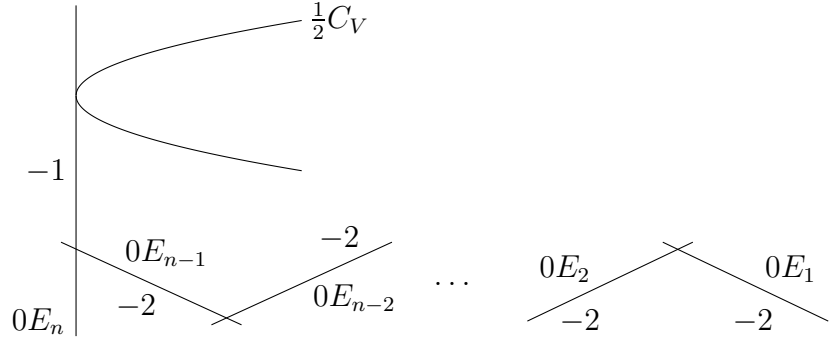

\begin{theorem}\label{t-hiromu-resol}
Let $(X, \Delta)$ be a two-dimensional klt pair with standard coefficients. 
Fix a closed point $x$ of $X$. 
Assume that $(X \setminus \{x\}, \Delta|_{X \setminus \{x\}})$ 
is log smooth, 
$(X, \Delta)$ is not log smooth, and  
all the irreducible components of $\Delta$ pass through $x$. 
Then one of the following holds. 
\begin{enumerate}
\item 
There exists a log resolution $\varphi : V \to X$ of $(X, \Delta)$ such that 
$\varphi(\Ex(\varphi)) = \{x \}$ and 
$\Delta_V$ is effective 
for the $\Q$-divisor $\Delta_V$ defined by $K_V + \Delta_V = \varphi^*(K_X+\Delta)$. 
\item 
$X$ is smooth, $\Delta = \frac{1}{2} C$ for a prime divisor $C$,  
and there exists a projective birational morphism $\varphi \colon V \to X$ from a smooth surface $V$ 
such that, for the $\Q$-divisor $\Delta_V$ defined by $K_V+\Delta_V =\varphi^*(K_X+\Delta)$, the dual graph of $\Delta_V$ is given by Figure \ref{figure-V}.\footnote{In Figure \ref{figure-V}, 
$C_V$ is the strict transform of $C$, 
$E_1, \ldots, E_n$ are the $\varphi$-exceptional prime divisors, 
the rational number in front of each prime divisor is 
its coefficient in $\Delta_V$, and the other numbers are the self-intersection numbers}

Specifically, there exists a sequence of projective birational morphisms of smooth surfaces 
\[
\varphi: V :=X_n \xrightarrow{\varphi_n} X_{n-1} \xrightarrow{\varphi_{n-1}} \cdots \xrightarrow{\varphi_1} X_0 =X
\]
satisfying the following property 
for 
the proper transforms $C_i$ of $C$ on $X_i$. 
\begin{enumerate}
    \item 
    For each $0 \leq i \leq n-1$, $C_i$ has a unique singular point $x_i$. 
%
Furthermore, $\mult_{x_i} C_i=2$. 
    \item 
    For each $0 \leq i \leq n-1$, $\varphi_{i+1} : X_{i+1} \to X_i$ is the blowup at $x_i$,  and $K_{X_{i+1}} + \frac{1}{2}C_{i+1} = \varphi^*_{i+1}(K_{X_{i}} + \frac{1}{2}C_{i})$. 
    \item 
    $C_V := C_n$ is a smooth prime divisor and $\Delta_V = \frac{1}{2} C_V$. 
    \item 
    The dual graph of $C_V \cup \Ex(\varphi)$ is given by 
    \[
    C_V = E_n - E_{n-1} - \cdots - E_1,
    \]
    where each $E_i$ denotes the proper transform of $\Ex(\varphi_i)$ on $V$ and 
    $C_V=E_n$ means that $(C_{V} \cap E_n)_{\red}$ is a single point with $C_V \cdot E_n =2$ (cf.\ Figure \ref{figure-V}). 
\end{enumerate}
\end{enumerate}
\end{theorem}

A closed point $x$ satisfying the above property (2) is called 
a {\em special point} {of} $(X, \Delta)$.

\begin{proof}
If $x$ is a singular point of $X$, then (1) holds by Proposition \ref{p-klt-sing}. 
Hence the problem is reduced to the case when $X$ is smooth. 
If $\Delta$ is not irreducible, then (1) holds by 
Lemma \ref{l-klt-resol-r=3} and 
Lemma \ref{l-klt-resol-r=2}. 
Therefore, we may assume that $\Delta$ is irreducible. 
We can write $\Delta = \frac{m-1}{m} C$ for an integer $m \geq 2$ 
and a prime divisor $C$. 
Since $(X, \Delta = \frac{m-1}{m} C)$ is not log smooth at $x$,  $C$ is singular at $x$, that is, $\mult_x C \geq 2$. 
If $m \geq 4$, then we have 
$\mult_x \Delta = \frac{m-1}{m} \mult_x C \geq \frac{3}{4} \cdot 2 = \frac{3}{2}$. 
In this case, (1) holds by Lemma \ref{l-klt-resol-r=1}. 
Therefore, we obtain $m=2$ or $m=3$
. 
It follows from a similar argument that $\mult_x C = 2$.


We start with the case when $m=2$. 
Let $\varphi_1: X_1 \to X$ be the blowup at $x$. 
Then we have 
\[
K_{X_1} + \frac{1}{2} C_{1} 
=\varphi_1^*\left(K_X+ \frac{1}{2}C\right)
\]
for 
$C_{1} := (\varphi_1)_*^{-1}C$.  
We have $C_{1} \cdot E =2$ for $E := \Ex(\varphi_1)$. 
There are the following three cases. 
\begin{enumerate}
    \item[(I)] $C_{1} \cap E$ consists of two points. 
    \item[(II)] $x_1 :=(C_{1} \cap E)_{\red}$ is one point, and $C_{1}$ is smooth at $x_1$   {and tangent to $E$}. 
    \item[(III)] $x_1:=(C_{1} \cap E)_{\red}$ is one point and $C_{1}$ is not smooth at $x_1$. 
\end{enumerate}
Assume (I). Then $C_{1} + E$ is simple normal crossing. 
In this case, (1) holds. 
Assume (II). 
In this case, we obtain (2). 
Assume (III). 
In this case, we take the blowup at $x_1$. 
Then, after replacing $(X_1, \Delta_{1} := \frac{1}{2}C_{1}, x_1)$ by $(X, \Delta, x)$, 
we can apply the same argument as above. 
Note that, for the blowup $\varphi_2 : X_2 \to X_1$ at $x_1$,  
the strict transforms $C_2$ and $(\varphi_2)_*^{-1}E$ of $C_{1}$ and $E$ are disjoint, because the following holds for $F := \Ex(\varphi_2)$: 
\[
C_2 \cdot (\varphi_2)_*^{-1}E 
= (\varphi_2^*C_1-2F) \cdot (\varphi_2)_*^{-1}E 
= C_1 \cdot E -2 F \cdot (\varphi_2)_*^{-1}E =2-2 \cdot 1=0.
\]
Repeat this argument, and then this procedure will terminate after finitely many times \cite[Corollary V.3.7]{hartshorne77}. 
Then (1) or (2) holds when $m=2$.

Let us treat the remaining case, that is, $m = 3$. We use the same notation as above. 
We then get 
\[
K_{X_1} + \frac{2}{3} C_{1} + \frac{1}{3}E 
=\varphi_1^*\left(K_X+ \frac{2}{3}C\right). 
\]
If (I) holds, then we get (1) as above. 
If (III) holds, then $\mult_{x_1}(\frac{2}{3} C_{1} + \frac{1}{3}E) \geq \frac{5}{3}$, 
so that we again obtain (1) by Lemma \ref{l-klt-resol-r=1}. 
Assume (II). 
We then get 
\[
K_{X_2} + \frac{2}{3} C_2 + \frac{1}{3}\varphi_2^{-1}E
=\varphi_2^*\left( K_{X_1} + \frac{2}{3} C_{1} + \frac{1}{3}E \right)
=\varphi_2^*\varphi_1^*\left(K_X+ \frac{2}{3}C\right). 
\]
Repeat this procedure, which will terminate after finitely many times. 
Hence (1) holds when $m=3$. 

\end{proof}

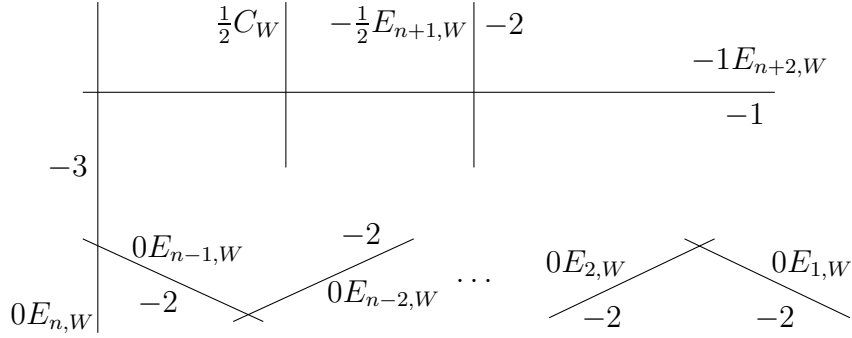
\begin{figure}
\centering
\begin{tikzpicture}
\draw (0.0, 1.9) node{$\frac{1}{2}C_W$};
\draw (2.0, 1.9) node{$-\frac{1}{2}E_{n+1, W}$};
\draw (6.8, 1.4) node{$-1E_{n+2, W}$};

\draw (-2.6, -2) node{$0E_{n, W}$};
\draw (-0.8, -1.1) node{$0E_{n-1, W}$};
\draw (1.8, -1.7) node{$0E_{n-2, W}$};
\draw (4.5, -1.3) node{$0E_{2, W}$};
\draw (7.5, -1.3) node{$0E_{1, W}$};

\draw (3.4, 1.9) node{$-2$}; 
\draw (6.6, 0.7) node{$-1$}; 
\draw (-2.4, 0) node{$-3$};  
\draw (-1.2, -1.8) node{$-2$}; 
\draw (1.5, -0.9) node{$-2$}; 
\draw (4.7, -2.0) node{$-2$}; 
\draw (7.0, -2.0) node{$-2$}; 


\draw(-2.2, 1)--(7, 1); 
\draw(0.5, 2.2)--(0.5, 0); 
\draw(3, 2.2)--(3, 0); 

\draw(-2,2.2)--(-2,-2.2); 

\draw(-2.2,-0.95)--(0.2,-2.05);
\draw(-0.2,-2.05)--(2.2,-0.95);
\draw[shift={(3,-1.5)}] node{$\hdots$}; 
\draw(4,-2)--(6.2,-0.95);
\draw(5.8,-0.95)--(8,-2);

\end{tikzpicture}
\caption{The coefficients of $\Delta_W$ and self-intersection numbers on $W$}\label{figure-W}
\end{figure}


\begin{remark} \label{r-cusp-resol}
We here summarise some properties of special points $x$ of $(X, \Delta)$. {This will be applied at the end of the proof of Theorem \ref{t-LDP-QFS}, and so the reader might postpone reading this remark until reaching Theorem \ref{t-LDP-QFS}.}
Let $(X, \Delta)$ be a two-dimensional klt pair with standard coefficients and 
let $x$ be a special point of $(X, \Delta)$, 
that is, Theorem \ref{t-hiromu-resol}(2) holds around $x$. 
Taking the blowups $\theta \colon W \to V$ at the non-log smooth point twice, the composition 
\[
\psi : W\xrightarrow{\theta} V \xrightarrow{\varphi} X
\]
is a log resolution of $(X, \Delta)$ (cf.\ Figure \ref{figure-V}, Figure \ref{figure-W}). 
Set $C_W := \theta_*^{-1}C_V$ and $E_{i, W} := \theta_*^{-1}E_i$. 
Let $E_{n+1, W}$ and $E_{n+2, W}$ be the $\theta$-exceptional prime divisors obtained by the first and second blowups, respectively. 


Set $K_W + \Delta_W = \psi^*(K_X+\Delta)$. 
By Figure \ref{figure-V} and 
$K_V +\frac{1}{2}C_V= K_V +\Delta_V=\varphi^*(K_X+\Delta)$, it is easy to see that   
\[
\Delta_W = \frac{1}{2}C_W -\frac{1}{2}E_{n+1, W} - E_{n+2, W}. 
\]
In Figure \ref{figure-W}, the rational number in front of each prime divisor is 
its coefficient in $\Delta_W$, and the other numbers indicate the self-intersection numbers. 
Set 
\[
\mathbb E_W := E_{1, W} \cup \cdots \cup E_{n, W}. 
\]
Given a $\psi$-exceptional prime divisor $E_{i, W}$ with $1 \leq i \leq n+2$, 
it follows from Figure \ref{figure-W} that {$E_{i, W}$ has an integer coefficient in $K_W+\Delta_W$} and 
 $\rdown{K_W+\Delta_W} \cdot E_{i, W} = 0$ 
  if and only if $1 \leq i \leq n$ (indeed, $\rdown{K_W+\Delta_W} \cdot E_{W, i} = -\frac{1}{2}(C_W +E_{n+1, W}) \cdot E_{W, i}$). 
The determinant $d(\mathbb E_W)$ of the intersection matrix $(E_{i, W} \cdot E_{j, W})_{1 \leq i, j\leq n}$ of $\mathbb E_W$ can be computed as follows: 
\begin{align*}
d(\mathbb E_W) = 
\det &\left[\begin{matrix}
-3 & 1 & 0 & 0 & \ldots & \ldots \\
1 & -2 & 1 & 0 & \ldots & \ldots \\
0 & 1 & -2 & 1 & \ldots & \ldots  \\
0 & 0 & 1 & -2 & \ldots & \ldots &  \\
\ldots & \ldots & \ldots & \ldots & \ddots & \ldots   \\
\ldots & \ldots & \ldots & \ldots & \ldots & \ddots 
\end{matrix} \right] = -3\det A_{n-1} - \det A_{n-2}, \textrm{ where} \\
A_n := &\left[\begin{matrix}
-2 & 1 & 0 & 0 & \ldots & \ldots \\
1 & -2 & 1 & 0 & \ldots & \ldots \\
0 & 1 & -2 & 1 & \ldots & \ldots  \\
0 & 0 & 1 & -2 & \ldots & \ldots &  \\
\ldots & \ldots & \ldots & \ldots & \ddots & \ldots   \\
\ldots & \ldots & \ldots & \ldots & \ldots & \ddots 
\end{matrix} \right] \textrm{ is an } n \times n \textrm{  matrix.}
\end{align*}
We see that
\[
\det A_n = -2 \det A_{n-1} - \det A_{n-2}, 
\]
which by induction yields $\det A_n = (-1)^n(n+1)$ for $n\geq 1$. In particular, 
\[
d(\mathbb E_W) = 
(-1)^n3n + (-1)^{n-1}(n-1) = (-1)^{n}(2n+1). 
\]
For its absolute value $a := |d(\mathbb E_W)| = 2n +1$, 
a point $x$ is called  \emph{of type $a$}. 
To summarise, the following hold. 
\begin{itemize}
    \item $n$ is the minimum number of blowups to resolve the singularity of $C$. 
    \item $n+2$ is the minimum number of blowups to reach a log resolution of $(X, \Delta = \frac{1}{2}C)$. 
    \item $d(\mathbb E_W) = \det (E_{i, W} \cdot E_{j, W})_{1 \leq i, j\leq n} = (-1)^n(2n+1)$. 
    \item $x$ is a special point of type $2n+1$. 
\end{itemize}



\end{remark}

We shall later need the following two consequences. 

\begin{corollary} \label{cor:Hiromu-resolution} 
Let $(X,\Delta)$ be a two-dimensional klt pair with standard coefficients. 
Then there exists a projective birational morphism $\varphi \colon V \to X$ 
such that {$(V, \Delta_V)$ is log smooth and $\Delta_V$ is effective},  where $\Delta_V$ is the $\Q$-divisor defined by $K_V + \Delta_V = \varphi^*(K_X+\Delta)$.
\end{corollary}

\begin{proof}
The assertion immediately follows from Theorem \ref{t-hiromu-resol}. 
\end{proof}

\begin{corollary}\label{c-klt-resol-r=1'}
We use Notation \ref{n-klt-resol-sm}. 
Assume $\mult_x \Delta \geq \frac{4}{3}$. 
Then there exists a log resolution $\varphi : V \to X$ of $(X, \Delta)$ such that $\Delta_V$ is effective 
for the $\Q$-divisor $\Delta_V$ defined by $K_V + \Delta_V = \varphi^*(K_X+\Delta)$. 
\end{corollary}

\begin{proof}
By Lemma \ref{l-klt-resol-r=3} and Lemma \ref{l-klt-resol-r=2},  we may assume that $\Delta$ is irreducible. 
We can write $\Delta = cC$ for some $c \geq 1/2$ and 
prime divisor $C$.  
By Lemma \ref{l-klt-resol-r=1}, we may assume that $\mult_x C=2$. 
It follows from $\mult_x \Delta = c \mult_x C =2c$ and 
the assumption $\mult_x \Delta \geq \frac{4}{3}$ that  $c \geq 2/3$. 
By decreasing the coefficient $c$, we may assume $c = 2/3$. 
Then the assertion follows from Theorem \ref{t-hiromu-resol}. 
\end{proof}

%% file: section4.tex
\section{Log liftablity of log del Pezzo surfaces}\label{s-log-lift}


Throughout this section, we work over an algebraically closed field $k$ of characteristic $p>0$. 
We start by stating the following conjecture. 

\begin{conjecture}\label{c-LDP-loglift}
Assume $p>5$. 
Let $(X, \Delta)$ be a log del Pezzo pair with standard coefficients. 
Then $(X, \Delta)$ is log liftable. 
\end{conjecture}

This conjecture is known to hold when $\Delta =0$ \cite[Theorem 1.2]{ABL20}. 
{Assuming} Conjecture \ref{c-LDP-loglift}, we can show that 
a log del Pezzo pair $(X, \Delta)$ with standard coefficients is quasi-$F$-split {when $p>41$}
(cf.\ Lemma \ref{l-log-lift-enough}). 
However, Conjecture \ref{c-LDP-loglift} is open as far as the authors know. 
As replacements of Conjecture \ref{c-LDP-loglift}, we {establish the following results:}
\begin{itemize}
\item Theorem \ref{t-bdd-or-liftable} stating that either $(X, \Delta)$ is log liftable or special points of $(X, \Delta)$ are bounded with respect to the arithmetic genus (see Subsection \ref{ss-bdd-or-lift1} and \ref{ss-bdd-or-lift2}),
\item Theorem \ref{thm:liftability-of-Hiromu-resolution} assuring the existence of a projective birational map $f \colon Y \to X$ such that $(Y,\Delta_Y)$ is  log smooth and lifts to $W(k)$, 
 where $K_Y+\Delta_Y = f^*(K_X+\Delta)$ (see Subsection \ref{ss-weak-log-lift}).
\end{itemize}

\subsection{Bounds for Mori fibre spaces}\label{ss-bdd-or-lift1}

In this subsection, we shall establish two auxiliary results: 
Proposition \ref{p-bound-rho1} and Proposition \ref{p-bound-rho2}, 
which will be 
used { to establish} the main theorem of this section: Theorem \ref{t-bdd-or-liftable}. 
In {this theorem and} its proof, 
the situation is as follows: {for a log del Pezzo pair}
$(X, \Delta)$ with standard coefficients, 
we take  the canonical model $f: Y \to X$ of $X$ (cf.\ Definition \ref{d-cano-model}), 
and we run a $K_Y$-MMP, whose end result is denoted by $Z$. 
In particular, $Z$ admits a Mori fibre space structure $\pi : Z \to T$. 
Proposition \ref{p-bound-rho1} and Proposition \ref{p-bound-rho2} treat the cases when  $\dim T=0$ and $\dim T=1$, respectively. 
The purpose of this subsection is to give an explicit upper bound for $(K_Z+C_Z) \cdot C_Z$ under suitable assumptions.  

\begin{proposition}\label{p-bound-rho1}
Let $Z$ be a canonical del Pezzo surface with $\rho(Z)=1$ and 
let $C_Z$ is a prime divisor such that $-(K_Z+\frac{1}{2}C_Z)$ is ample. 
Then $(K_Z + C_Z) \cdot C_Z <18$. 
\end{proposition}

\begin{proof}
It is easy to show that
\[
 {\rm Cl}\,(Z) /\! \equiv \hspace{2mm} \simeq\, \mathbb Z, 
\]
where $\equiv$ denotes the numerical equivalence. 
Let $A$ be an ample Weil divisor on $Z$ which generates ${\rm Cl}\,(Z) / \equiv$. 
We have 
\[
-K_Z \equiv rA \qquad \text{and} \qquad C_Z \equiv mA
\]
for some $r \in \Z_{>0}$ and $m \in \Z_{>0}$. 
Then $-(K_Z+ \frac{1}{2} C_Z)$ 
is numerically equivalent to $(r -\frac{m}{2})A$. 
Since $-(K_Z+ \frac{1}{2} C_Z)$ is ample, we obtain $m < 2r$. 
We then get 
\[
(K_Z + C_Z) \cdot C_Z
= (-r +m) \cdot m A^2 
< (-r + 2r) \cdot 2r A^2 =2r^2 A^2 =2K_Z^2 \leq 18, 
\]
where the last inequality holds by the fact that $Z$ is a projective rational surface with canonical singularities (indeed, we have $K_Z^2 =K_{Z'}^2 \leq 9$ for the minimal resolution $\mu : Z' \to Z$). 
\end{proof}

\begin{proposition}\label{p-bound-rho2}
Let $(Z, \Delta_Z = \frac{1}{2} C_Z + B_Z)$ be a weak log del Pezzo pair, 
where $Z$ is canonical, $\rho(Z)=2$, $C_Z$ is a non-smooth prime divisor, and $B_Z$ is an effective $\Q$-divisor satisfying $C_Z \not\subseteq \Supp\,B_Z$. 
Assume that ${\rm (a)}$--${\rm (c)}$ hold. 
\begin{enumerate}
\item[(a)] All the coefficients of $B_Z$ are $\geq \frac{1}{3}$. 
    \item[(b)] There exists a surjective morphism $\pi : Z \to T$ to a smooth projective curve $T$ with $\pi_*\MO_Z = \MO_T$. 
    \item[(c)] There does not exist a prime divisor $D$ on $Z$ such that $K_Z \cdot D < 0$ and 
    $D^2 <0$. 
\end{enumerate}
Then one of the following holds. 
\begin{enumerate}
    \item $(\Delta_Z)_{\red} \cdot F \leq 3$ for a general fibre $F$ of $\pi : Z \to T$. 
    \item $(K_Z+ C_Z) \cdot C_Z <36$. 
\end{enumerate}
\end{proposition}

{Assumption (a) is natural when one considers minimal Gorenstein (canonical) models of klt surfaces (see Lemma \ref{l-coeff-1/3}), 
 whilst Assumption (c) says that there is no projective birational contraction $Z \to Z'$ induced by a $K_Z$-MMP.
\begin{remark}
{Let us briefly explain the motivation for Proposition \ref{p-bound-rho2}. Assume that $Z$ is smooth, and so $Z \simeq \mathrm{Proj}_{\mathbb P^1}(\mathcal{O}_{\mathbb P^1} \oplus \mathcal{O}_{\mathbb P^1}(n))$ is the $n$-th Hirzebruch surface. In what follows we suppose that $n\geq 6$ and show that (1) holds. Since $-(K_Z+\frac{1}{2}C_Z)$ is big, it is easy to see that if $n < 6$, then the arithmetic genus of $C_Z$ 
 must be bounded from above by a constant (see the proof of Proposition \ref{p-bound-rho2} for the exact calculations), and so the validity of (2) in this case should not come as a surprise. 

Let $\Gamma \simeq \mathbb P^1$ be the negative section such that $\Gamma^2 = -n$ and let $F$ be a general fibre. As $C_Z$ is singular, $C_Z \neq \Gamma$ and $C_Z \cdot F \geq 2$. 
Write $B_Z = a\Gamma + B'_Z$ for $a \in \mathbb{Q}_{\geq 0}$, where $B'_Z \geq 0$ and $\Gamma \not \subseteq \Supp B'_Z$. Then $a \geq \frac{n-2}{n}$ by the following calculation:
\[
0 \geq (K_Z+\Delta_Z) \cdot \Gamma \geq (K_Z+a\Gamma)\cdot \Gamma = (K_Z+\Gamma)\cdot \Gamma + (a-1)\Gamma^2 = -2 - (a-1)n.
\]
Since $n \geq 6$, we get $a \geq \frac{2}{3}$.

To show (1), we first observe that $0 > (K_Z+\Delta_Z) \cdot F = -2 +  \Delta_Z \cdot F$, and so $\Delta_Z \cdot F \leq 2$. In particular,
\begin{equation} \label{eq:hirzebruch-bounds}
2 > \Delta_Z \cdot F = \frac{1}{2}C_Z \cdot F + a\Gamma \cdot F + B'_Z \cdot F \geq 1 + \frac{2}{3} + B'_Z \cdot F, 
\end{equation}
that is $B'_Z \cdot F < \frac{1}{3}$. Since the coefficients in $B'_Z$ are at least $\frac{1}{3}$, we get that $B'_Z \cdot F = 0$. Moreover, the calculation in (\ref{eq:hirzebruch-bounds}) also shows that $C_Z \cdot F = 2$. This concludes the proof of (1).}
\end{remark}
\begin{proof}[Proof of Proposition \ref{p-bound-rho2}]
Since $Z$ is of del Pezzo type and  $\rho(Z)=2$, 
$\NE(Z)$ has exactly two extremal rays and there exist the corresponding contraction morphisms (cf.\ Remark \ref{r-log-dP-summary}). 
One of them corresponds to $\pi : Z \to T = \P^1$. 
Let $R$ be the other extremal ray. 
Fix a curve $\Gamma$ with $[\Gamma] \in R$. 
By $\rho(Z)=2$, we have $\Gamma^2 = 0$ or $\Gamma^2<0$ 
(cf.\ \cite[Proof of the case where $C^2>0$ of Theorem 3.21]{tanaka12}). 
Set $\gamma := -\Gamma^2$ and let $F$ be a general fibre of $\pi : Z \to T$. 
By adjunction, we have $F \simeq \P^1$ and $K_Z \cdot F = -2$. In particular, 
\[
\NE(Z) = \R_{\geq 0}[F] + \R_{\geq 0}[\Gamma]. 
\]

\setcounter{step}{0}
\begin{step}\label{s1-bound-rho2}
Assume $\Gamma^2=0$. Then $(K_Z+C_Z) \cdot C_Z <  16$.
\end{step}

\begin{proof}[Proof of Step \ref{s1-bound-rho2}]
By $\Gamma^2=0$, the corresponding morphism is a Mori fibre space to a curve: 
$\pi' : Z \to T' =\P^1$. 
Let $F$ and $F'$ be general fibres of $\pi : Z \to T$ and $\pi' : Z \to T'$, respectively. 
In particular, $F$ and $F'$ are effective Cartier divisors with 
$F \simeq F' \simeq \P^1$ and $-K_Z \cdot F = -K_Z \cdot F' = 2$. 
Set $d := F \cdot F' \in \Z_{>0}$. 
By $\NE(Z) = \R_{\geq 0} [F] + \R_{\geq 0} [F']$
We can write 
\begin{equation}\label{e1-bound-rho2}
C_Z \equiv a F + a' F'
\end{equation}
for some $a, a' \in \R_{\geq 0}$. 
By Kleimann's criterion, the big $\Q$-divisor $-(K_Z+ \frac{1}{2} C_Z)$ is automatically ample. 
It follows from (\ref{e1-bound-rho2}) that 
\[
0> \left(K_Z+ \frac{1}{2} C_Z\right) \cdot F = -2  + \frac 1 2 a'd, 
\]
and hence $a'd < 4$. By symmetry, we also obtain $ad < 4$. 
{Then} 
\[
C_Z^2 = 2aa' F \cdot F' = \frac{2(ad)(a'd)}{d} < \frac{32}{d} \leq 32, 
\]
which implies 
\[
(K_Z+ C_Z) \cdot C_Z = 
(K_Z+ \frac{1}{2} C_Z + \frac{1}{2} C_Z ) \cdot C_Z < \frac{1}{2}  C_Z^2 < 16. 
\]
This completes the proof of Step \ref{s1-bound-rho2}. 
\end{proof}

\begin{step}\label{s2-bound-rho2}
If $\Gamma^2 <0$, then 
\begin{enumerate}
\renewcommand{\labelenumi}{(\roman{enumi})}
\item $\Gamma \simeq \P^1$, 
\item $\Gamma \neq C_Z$, and 
\item the coefficient of $\Gamma$ in $\Delta_Z$ is nonzero. 
\end{enumerate}
\end{step}

\begin{proof}[Proof of Step \ref{s2-bound-rho2}]
Assume $\Gamma^2<0$. 
We have 
\[
(K_Z + \Gamma) \cdot \Gamma < (K_Z+\Delta_Z) \cdot \Gamma \leq 0, 
\]
which implies $\Gamma \simeq \P^1$ \cite[Theorem 3.19]{tanaka12}. 
Thus (i) holds. Since $C_Z$ is singular, (i) implies (ii). 

Let us show (iii). 
Suppose that the coefficient of $\Gamma$ in $\Delta_Z$ is zero. 
We then get $\Delta_Z \cdot \Gamma \geq 0$, which implies  
\[
0> (K_Z+ \Delta_Z) \cdot \Gamma  \geq  K_Z \cdot \Gamma. 
\]
This contradicts (c). 
This completes the proof of Step \ref{s2-bound-rho2}. 
\end{proof}

\begin{step}\label{s3-bound-rho2}
Assume that $\Gamma^2<0$ and $\pi|_{\Gamma}: \Gamma \to T$ is not birational. 
Then $(K_Z+C_Z) \cdot C_Z <  4$.
\end{step}

\begin{proof}[Proof of Step \ref{s3-bound-rho2}]
By $\NE(Z) = \R_{\geq 0}[F] + \R_{\geq 0}[\Gamma]$, 
$\Gamma$ is the unique curve on $Z$ whose self-intersection number is negative. Hence $-(K_Z+\frac{1}{2}C_Z + c\Gamma)$ is nef and big, where $c$ denotes the coefficient of $\Gamma$ in $\Delta_Z$. 
By Step \ref{s2-bound-rho2}(iii) and (a), we get $c \geq \frac{1}{3}$. 
Since $F$ is nef and $\Gamma \neq C_Z$ (Step \ref{s2-bound-rho2}(ii)), we have 
\begin{equation}\label{e2-bound-rho2}
0> (K_Z+\Delta_Z) \cdot F 
\geq \left( K_Z+ \frac{1}{2}C_Z +c\Gamma \right) \cdot F 
\geq -2 + \frac{1}{2}\cdot 2 +c\Gamma \cdot F, 
\end{equation}
where the last inequality follows from the assumption that $C_Z$ is singular, 
and hence $\pi|_T : C_Z \to T$ is of degree $\geq 2$. 
We then get $c \Gamma \cdot F <1$. 
By $c \geq \frac{1}{3}$, we have $\Gamma \cdot F < 3$. 
Since $\pi|_{\Gamma}: \Gamma \to T$ is not birational, 
it holds that $\Gamma \cdot F \geq 2$, and so $\Gamma \cdot F =2$.  
By (\ref{e2-bound-rho2}), we also have $C_Z \cdot F =2$. 
It follows from $\NE(Z) = \R_{\geq 0}[F] + \R_{\geq 0}[\Gamma]$ that we can write 
$C_Z \equiv a\Gamma + b F$ for some $a, b \in \R_{\geq 0}$. 
By taking the intersection number with $F$, we get $a=1$. 
In particular, 
\begin{equation}\label{e3-bound-rho2}
C_Z \equiv \Gamma + b F\qquad \text{and}\qquad 
C_Z \cdot \Gamma = -\gamma + 2b. 
\end{equation}
Again by (\ref{e2-bound-rho2}), we obtain 
\[
0> (K_Z+\Delta_Z) \cdot F \geq 
\left(K_Z+ \frac{1}{2} C_Z + c \Gamma\right) \cdot F = 
-2 + 1 + 2c,  
\]
which implies $c < \frac{1}{2}$. 
As $-(K_Z+\Delta_Z)$ is nef, 
the following holds: 
\begin{align*}
0 &\geq (K_Z+\Delta_Z) \cdot \Gamma \\
&\geq \left(K_Z+ \frac{1}{2} C_Z + c \Gamma\right) \cdot \Gamma \\
&= \left( (K_Z+\Gamma)+ \frac{1}{2} C_Z + (c-1)\Gamma\right) \cdot \Gamma \\
&\geq -2 + \frac{1}{2} (-\gamma +2b) + \gamma(1-c)\\
&>     -2 + \frac{1}{2} (-\gamma +2b) + \gamma\left( 1-\frac{1}{2}\right)\\
&=-2 +b, 
\end{align*}
where the third inequality holds by (\ref{e3-bound-rho2}) and $(K_Z+\Gamma) \cdot \Gamma \geq -2$.  
We then get $b <2$, which implies 
\[
C_Z^2 =(\Gamma +bF)^2 = \Gamma^2 +4b < 0+4 \cdot 2=8. 
\]
Then 
\[
(K_Z+C_Z) \cdot C_Z = \Big(K_Z+\Delta_Z - \big(\Delta_Z - \frac{1}{2}C_Z\big) + \frac{1}{2}C_Z\Big) \cdot C_Z
\leq \frac{1}{2}C_Z^2 <4. 
\]
This completes the proof of Step \ref{s3-bound-rho2}. 
\qedhere

\end{proof}

\begin{step}\label{s4-bound-rho2}
Assume that $\Gamma^2<0$ and $\pi|_{\Gamma} : \Gamma \to T$ is birational. 
Then (1) or (2) in the statement of Proposition \ref{p-bound-rho2} holds. 
\end{step}

\begin{proof}
Since $\pi|_{\Gamma}: \Gamma \to T = \P^1$ is a finite birational morphism from a curve $\Gamma$ to a smooth curve $T$, $\pi|_{\Gamma}$ is an isomorphism, that is,  $\Gamma$ is a section of $\pi : Z \to T$. 
We can write  
\[
\Delta_Z = \frac{1}{2} C_Z + c \Gamma + \Delta'_Z, 
\]
where $c$ is the coefficient of $\Gamma$ in $\Delta_Z$ and 
$\Delta'_Z$ is a nef effective $\Q$-divisor. 
If $c \geq \frac{2}{3}$, then $\Delta'_Z$ is $\pi$-vertical (that is, $\pi(\Supp \Delta'_Z) \subsetneq T$) {and} $C_Z \cdot F =2$. 
Then (1) holds. 
Hence we may assume that 
\[
\frac{1}{3} \leq c < \frac{2}{3}.
\]
By $\NE(Z) = \R_{\geq 0}[\Gamma] + \R_{\geq 0}[F]$, we can write  
\[
C_Z = a \Gamma + b F 
\]
for some $a, b \in \R_{\geq 0}$. 
By $a = C_Z \cdot F$, we have $a \in \Z_{\geq 0}$. 
Moreover, we have $a = C_Z \cdot F \geq 2$, because $C_Z$ is a singular prime divisor. 
Since $-(K_Z+ \Delta_Z)$ is nef and big, we obtain 
\begin{enumerate}
\item[(I)] $-(K_Z+\Delta_Z) \cdot F >0$. 
\item[(II)] $-(K_Z+\Delta_Z) \cdot \Gamma \geq 0$. 
\end{enumerate}
The inequality (I) implies 
\[
0> (K_Z+\Delta_Z) \cdot F \geq 
\left(K_Z+ \frac{1}{2} C_Z + c \Gamma\right) \cdot F = 
-2 + \frac{1}{2} a + c,  
\]
which yields $\frac{1}{2}a + c < 2$ and $a \in \{ 2, 3\}$.    
The inequality (II) implies 
\begin{align*}
0 &\geq (K_Z+\Delta_Z) \cdot \Gamma \\
&\geq \left(K_Z+ \frac{1}{2} C_Z + c \Gamma\right) \cdot \Gamma \\
&= \left( (K_Z+\Gamma) + \frac{1}{2} C_Z + (c-1) \Gamma\right) \cdot \Gamma \\
&\geq -2 + \frac{1}{2}C_Z \cdot \Gamma + \gamma(1-c)\\
&= -2 + \frac{1}{2} (-a\gamma +b) + \gamma(1-c)\\
&=-2 + \frac{1}{2}b + \gamma -\gamma\left( \frac{1}{2}a +c\right)\\
&> -2 + \frac{1}{2}b + \gamma -2\gamma.
\end{align*}
Hence $b < 2\gamma + 4$. 
The above displayed equation, together with $C_Z \cdot \Gamma \geq 0$, also implies the following {(see the fourth line)}: 
\[
\gamma \leq \frac{2}{1-c}. 
\]
We further have 
\[
C_Z^2 = (a\Gamma + bF)^2 = a^2(-\gamma) +2ab <2ab, 
\]
which implies 
\[
(K_Z+C_Z) \cdot C_Z =   \Big(K_Z+\Delta_Z - \big(\Delta_Z - \frac{1}{2}C_Z\big) + \frac{1}{2}C_Z\Big) \cdot C_Z
\leq \frac{1}{2}C_Z^2 < ab. 
\]
In what follows, we separately treat the two cases: $a=2$ and $a=3$. 

Assume $a=2$. 
Then we obtain 
\[
\gamma \leq \frac{2}{1-c} < \frac{2}{1- \frac{2}{3}} = 6
\]
and $b <2\gamma + 4 <16$.
Therefore, $(K_Z+C_Z) \cdot C_Z < ab <32$. 

Assume $a=3$. 
In this case, we have $c < 2 -\frac{1}{2}a =\frac{1}{2}$  as $\Delta_Z \cdot F < 2$ (see the first paragraph of the proof of Step 4). 
Then we obtain 
\[
\gamma \leq \frac{2}{1-c} < \frac{2}{1- \frac{1}{2}} = 4
\]
and $b <2\gamma + 4 <12$. 
Therefore, $(K_Z+C_Z) \cdot C_Z < ab <36$. 
In any case, (2) holds. 
This completes the proof of Step \ref{s4-bound-rho2}. 
\end{proof}
Step \ref{s1-bound-rho2}, 
Step \ref{s3-bound-rho2}, and 
Step \ref{s4-bound-rho2} complete the proof of 
Proposition \ref{p-bound-rho2}. 
\end{proof}

\subsection{Log liftablity for unbounded log del Pezzo pairs}\label{ss-bdd-or-lift2}

The purpose of this subsection is to prove one of the main results of this section: Theorem \ref{t-bdd-or-liftable}. 
We start by recalling the definition of canonical models. 

\begin{definition}\label{d-cano-model}
Let $X$ be a normal surface. 
We say that $f: Y \to X$ is the {\em canonical model} over $X$ 
if 
\begin{enumerate}
    \item $f$ is a projective birational morphism from a canonical surface $Y$, and  
    \item $K_Y$ is $f$-ample. 
\end{enumerate}
In this case, also $Y$ is called the {\em canonical model} over $X$. 
\end{definition}

\begin{remark}\label{r-cano-model}
Let $X$ be a normal surface. 
\begin{enumerate}
    \item It is well known that canonical models over $X$ {are Gorenstein and} are unique up to isomorphisms over $X$. 
    \item Given a normal surface $X$, 
    the canonical model $Y$ over $X$ is obtained as follows.  
    We first take the minimal resolution $\mu : X' \to X$ of $X$. 
    Then the canonical model $Y$ is obtained by contracting all the $(-2)$-curves contained in $\Ex(\mu)$. 
    In particular, we have the induced morphisms: $\mu : X' \to Y \xrightarrow{f} X$. 
\end{enumerate}
\end{remark}

The following lemma {highlights} the key property of canonical models { that we shall use later on}. 
We emphasise that our proof of Theorem \ref{t-bdd-or-liftable} does not work 
if we use the minimal resolution of $X$ instead of the canonical model over $X$.


\begin{lemma}\label{l-coeff-1/3}
Let $X$ be an 
lc surface and let $f: Y \to X$ be the canonical model over $X$. Let $B$ be the effective $\Q$-divisor defined by 
\[
K_{Y} + B = f^*K_X. 
\]
Then $\Ex(f) = \Supp B$ and every coefficient of $B$ is $\geq \frac{1}{3}$. 
\end{lemma}

\begin{proof}
Let $\mu : X' \to X$ be the minimal resolution of $X$, 
so that we have the induced projective birational morphisms (Remark \ref{r-cano-model}): 
\[
\mu: X' \xrightarrow{\mu_Y} Y \xrightarrow{f} X. 
\]
We define an effective $\Q$-divisor $B'$ on $X'$  by 
\[
K_{X'} + B' =\mu^*K_X. 
\]
Fix a $\mu$-exceptional prime divisor $E$ on $X'$ satisfying $E^2 \leq -3$. 
We can write 
\[
B' = c E + B'', 
\]
where $c$ is the coefficient of $E$ in $B'$ and $B''$ is an effective $\Q$-divisor satisfying $E \not\subseteq \Supp B''$. 
By $(\mu_Y)_*B' = B$, it suffices to show that $c \geq 1/3$. 
This holds by the following computation: 
\[
0 = \mu^*K_X \cdot E= (K_{X'} + B' ) \cdot E = 
(K_{X'} + E + (c-1)E +B'' ) \cdot E 
\]
\[
\geq -2 + (1-c) (-E^2) \geq -2 +(1-c) 3, 
\]
where the last inequality follows from $1 - c \geq 0$, 
which is assured by the assumption that $(X, 0)$ is lc. 
\end{proof}

\begin{lemma}\label{l-vertical-lift}
Assume $p>3$. 
Let $Z$ be a projective normal surface and 
let $D$ be a reduced Weil divisor.   
Let $\pi : Z \to T$ be a morphism to a smooth projective curve $T$ such that $\pi_*\MO_Z = \MO_T$. 
Let $F$ be a general fibre of  $\pi: Z \to T$. 
Assume that the following hold. 
\begin{enumerate}
    \item $K_Z \cdot F <0$. 
    \item $D \cdot F \leq 3$.
\end{enumerate}
Then $(Z, D)$ is log liftable. 
\end{lemma}

\begin{proof}
Let $f\colon W\to Z$ be a log resolution of $(Z, D)$. 
Set $D_W: = f^{-1}_{*}D +\Exc(f)$. 
We aim to lift $(W, D_W)$ to $W(k)$. 
{\cred  
In order to obtain such a lift, it suffices to show the following \cite[Theorem 2.3]{Kawakami-Nagaoka22}: 
\[
H^2(W, T_W(-\log D_W))=0\qquad \text{and}\qquad H^2(W, \sO_W)=0, 
\]
where $T_W(-\log D_W)\coloneqq \mathcal Hom_{\MO_{W}}(\Omega^{1}_W(\log D_W), \sO_W)$
denotes the logarithmic tangent sheaf.}
The latter vanishing follows from (1) and Serre duality.  
Again by Serre duality, we have 
\[
H^2(W, T_W(-\log D_W))
\simeq H^0(W, \Omega^1_W(\log D_W)\otimes \MO_W(K_W))^*, 
\]
where $(-)^*$ denotes the dual vector space. 
{\cred 
Since $f_{*}\Omega^1_W(\log D_W)$ is torsion-free, 
we obtain the following injective $\MO_W$-module homomorphism by applying the double dual $(-)^{**} := \mathcal Hom_{\MO_Z}( \mathcal Hom_{\MO_Z}( -, \MO_Z), \MO_Z)$: 
\[
f_{*}\Omega^1_W(\log D_W)\otimes \MO_W(K_W)\hookrightarrow (f_{*}\Omega^1_W(\log D_W)\otimes \MO_W(K_W))^{**}=(\Omega^{{[1]}}_Z(\log D)\otimes \MO_Z(K_Z))^{**}, 
\]
where the equality at last holds by the fact that 
the second and third terms are reflexive sheaves which coincide 
on $Z\setminus f(\Exc(f))$ (cf.\ \cite[Lemma 4.1]{Kaw3}). 
Thus it suffices to show that 
\[
H^0(Z, (\Omega^{{[1]}}_Z(\log D)\otimes \MO_Z(K_Z))^{**})=0.
\]
Suppose by contradiction that \[H^0(Z, (\Omega^{{[1]}}_Z(\log D)\otimes \MO_Z(K_Z))^{**})\neq 0.\]
Let $j\colon U\hookrightarrow Z$ be the open immersion from the log smooth locus $U$ of $(Z,D)$. 
It follows from  $(\Omega^{{[1]}}_Z(\log D)\otimes \MO_Z(K_Z))^{**}=j_{*}(\Omega^1_U(\log\,D)\otimes \sO_U(K_U))$ 
that 
\begin{align*}
    H^0(Z, (\Omega^{{[1]}}_Z(\log D)\otimes \MO_Z(K_Z))^{**}) &= H^0(U, \Omega^{1}_U(\log D)\otimes \MO_U(K_U)))\\
    &\simeq \Hom_{\sO_U}(\sO_U(-K_U), \Omega^{1}_U(\log D))\\
    &= \Hom_{\sO_Z}(\sO_Z(-K_Z), \Omega^{[1]}_Z(\log D)).
\end{align*}
Thus there exists an injective $\MO_Z$-module homomorphism 
$s\colon \MO_Z(-K_Z) \to \Omega^{{[1]}}_Z(\log D)$.}


Note that every $\pi$-horizontal prime divisor $C \subseteq \Supp D$ is generically \'etale over $T$, as otherwise the assumption (2) leads a contradiction: $3 \geq D \cdot F \geq C \cdot F \geq p >3$.
{Thus} 
we can find a non-empty open subset $T'$ of $T$ such that $(Z', D')$ is simple normal crossing over $T'$, 
where $Z' := \pi^{-1}(T')$ and $D' := D|_{Z'}$ 
(indeed, we may assume that $\pi|_{Z'} : Z' \to T'$ is smooth and $\pi|_{D'} : D' \to T'$ is \'etale). 
We then 
obtain the following horizontal exact sequence:
\begin{equation*}
    \begin{tikzcd}
        & & \MO_{Z'}(-K_{Z'}) \arrow[dotted]{ld}\arrow{d}{s}\arrow{rd}{t} && \\
       0 \arrow{r} & \sO_{Z'}((\pi|_{Z'})^{*}K_{T'}) \arrow{r} & \Omega^{1}_{Z'}(\log D') \arrow{r}{\rho} & \Omega^{1}_{Z'/T'}(\log D') \arrow{r} & 0,
    \end{tikzcd}
\end{equation*}
where $s$ denotes $s|_{Z'}$ by abuse of notation.
The construction of this exact sequence is as follows. 
When $D=0$, this is the usual relative differential sequence \cite[Proposition II.8.11]{hartshorne77}. 
When $D\neq 0$, we define $\rho$ by $d(f^{*}{\tau})\mapsto 0, dz/z\mapsto d\tau/\tau$, where $\tau$ is a coordinate on $T$ and $z$ is a local equation of $D$. Note that $f^{*}{\tau}$ and $z$ form a coordinate on $Z$ since $D$ is simple normal crossing over $T$. 

Set $t := \rho \circ s$. 
Suppose that $t$ is nonzero. 
Then, by restricting $t$ to $F$, we have an injective $\sO_F$-module homomorphism $t|_{F}\colon \sO_F(-K_F) \hookrightarrow 
\Omega^1_{F}(\log (D|_F))=\sO_F(K_F+D|_{F})$, where the injectivity holds since $F$ is chosen to be general. 
{This,} together with (2), 
implies that 
\[
2=\deg (-K_F)\leq\deg(K_F+(D|_F))\leq 1,
\]
which is a contradiction.
Thus $t$ is zero and the homomorphism $s$ induces an injection  $\sO_Z(-K_Z)\hookrightarrow \sO_Z(\pi^{*}K_T)$.
By taking the restriction to $F$, we get
\[
2=\deg (-K_F)\leq\deg(\pi^*K_T|_F)=0,
\]
which is again a contradiction. 
Therefore, we obtain the required vanishing.
\end{proof}

We are ready to prove the main result of this subsection. 

\begin{theorem}\label{t-bdd-or-liftable}
Assume $p>3$. 
Let $(X, \Delta)$ be a log del Pezzo pair with standard coefficients. 
{\cred Assume that there exists  a prime divisor $\Gamma$ on $X$ such that $\Gamma \subseteq \Supp\,\Delta$ and $\Gamma$ has a singular point at which $X$ is smooth.} 
Then one of the following holds. 
\begin{enumerate}
    \item $(X, \Delta)$ is log liftable. 
    \item If $C$ is a prime divisor on $X$ such that $C \subseteq  \Supp\,\Delta$ and  $C$ has a singular point at which $X$ is smooth, then 
    \[
    p_a(C_{X'})  \leq 18
    \]
    for the minimal resolution $\mu : X' \to X$ of $X$ and the proper transform $C_{X'} := \mu_*^{-1}C$ of $C$ on $X'$. 
\end{enumerate}
\end{theorem}

\begin{proof}
{\cred Fix a  prime divisor $C$ such that $C \subseteq \Supp\,\Delta$ and that $C$ has a singular point at which $X$ is smooth.} 
We can write 
\[
\Delta = a C + B, 
\]
where $a$ is the coefficient of $C$ in $\Delta$ and 
$B$ is an effective $\Q$-divisor satisfying $C\not\subseteq \Supp B$. 

Let $f : Y \to X$ be the canonical model over $X$ (cf.\ Definition \ref{d-cano-model}, Remark \ref{r-cano-model}). 
We run a $K_Y$-MMP: 
\[
g : Y  =:Y_0 \xrightarrow{g_0} Y_1 \xrightarrow{g_1} Y_2 \xrightarrow{g_2} \cdots 
\xrightarrow{g_{\ell-1}} Y_{\ell} =: Z,  
\]
where $Z$ is its end result and each $g_i: Y_i \to Y_{i+1}$ is a {projective} birational morphism such that 
$E_i := \Ex(g_i)$ is a prime divisor with $K_{Y_i} \cdot E_i < 0$. 
{We may assume that none of the $K_Z$-negative extremal rays induces a projective birational morphism (in other words,  every induced contraction is a Mori fibre space).}
We define effective $\Q$-divisors $\Delta_Y$ and $\Delta_Z$ by 
\[
K_Y +\Delta_Y = f^*(K_X+\Delta)\qquad\text{and} \qquad \Delta_Z := g_*\Delta_Y. 
\]
Set $C_Y := f_*^{-1}C$ and $C_Z := g_*C_Y$. 
Then $C_Y$ and $C_Z$ are singular prime divisors 
(note that $g(C_Y)$ is not a point, 
because the image of $C_Y$ on $Y_i$ is singular, 
whilst $E_i \simeq \P^1$ by \cite[Theorem 3.19(1)]{tanaka12}). 
Then $(Y, \Delta_Y)$ and $(Z, \Delta_Z)$ are weak del Pezzo pairs. 
Since $-K_Y$ is big, the end result $Z$ has a $K_Z$-Mori fibre space structure $\pi: Z \to T$, that is, $\pi$ is a morphism to a normal projective variety $T$, 
$-K_Z$ is $\pi$-ample, $\pi_*\MO_Z= \MO_T$, and $\dim Z > \dim T$. 
In particular, $\dim T =0$ or $\dim T=1$.

Let $\mu : X' \to X$ and $\mu_Z : Z' \to Z$ be 
the minimal resolutions of $X$ and $Z$, respectively. 
Then $\mu : X' \to X$ factors through $f: Y \to X$ (Remark \ref{r-cano-model}): $\mu: X' \xrightarrow{\mu_Y} Y \xrightarrow{f} X$. 
Furthermore, the induced resolution $g \circ \mu_Y : X' \to Z$  factors through $\mu_Z: Z' \to Z$, 
so that we get the following commutative diagram: 
\[
\begin{tikzcd}
& X' \arrow[d, "\mu_Y"] \arrow[ld, "h"'] \arrow[ddr, bend left, "\mu"]\\
Z' \arrow[d, "\mu_Z"'] & Y\arrow[ld, "g"'] \arrow[rd, "f"]\\
Z & & X.
\end{tikzcd}
\]
Recall that $h : X' \to Z'$ is a composition of blowups of points. 
Hence we have $p_a(C_{X'}) \leq p_a(C_{Z'})$ \cite[Corollary V.3.7]{hartshorne77}, 
where $C_{Z'} := h_*C_{X'}$. 
Then the problem is reduced to showing $p_a(C_{Z'}) \leq 18$, which is equivalent to 
$(K_{Z'} + C_{Z'}) \cdot C_{Z'} \leq 34$, {where $(K_{Z'} + C_{Z'}) \cdot C_{Z'} = 2p_a(C_{Z'}) - 2$}. 
Since the equality $K_{Z'} + C_{Z'} + \Gamma = \mu_Z^*(K_Z+C_Z)$ 
holds 
and some $\mu_Z$-exceptional effective $\Q$-divisor $\Gamma$, 
we obtain 
\[ 
(K_{Z'}+C_{Z'}) \cdot C_{Z'} \leq (K_{Z'} + C_{Z'} +\Gamma) \cdot C_{Z'} 
= (K_Z+C_Z) \cdot C_Z. 
\]
To summarise, in order to show $p_a(C_{X'})  \leq 18$, 
it suffices to prove $(K_Z+C_Z) \cdot C_Z <36$.

Assume 
{$\dim T =0$.} 
It follows from Proposition \ref{p-bound-rho1} that $(K_Z+C_Z) \cdot C_Z < 18 <36$. 
Thus (2) holds.

Assume 
{$\dim T =1$.} 
{To apply Proposition \ref{p-bound-rho2}, 
we need to check that the assumptions (a)--(c) of  Proposition \ref{p-bound-rho2} hold. The assumption (a) holds by Lemma \ref{l-coeff-1/3}, and the assumption (b) is clear. Now suppose that  the assumption (c) does not hold, that is, 
there exists  a prime divisor $D$ on $Z$ satisfying $K_Z \cdot D <0$ and $D^2 <0$. Then $D$ is a generator of a $K_Z$-negative extremal ray, the contraction of which is birational. This is impossible by our assumption on $Z$, and so  (c) holds.}   

By Proposition \ref{p-bound-rho2}, we get that 
$(K_Z+C_Z) \cdot C_Z <36$ or 
$(\Delta_Z)_{\red} \cdot F \leq 3$ for a general fibre $F$ of $\pi : Z \to T$. 
In the former case, (2) holds. 
Hence we may assume that 
 $(\Delta_Z)_{\red} \cdot F \leq 3$. 
 {Then} $(Z, \Delta_Z)$ is log liftable by Lemma \ref{l-vertical-lift}, {and so,}
 automatically, $(Y, \Delta_Y)$ is log liftable (cf.\ Subsection \ref{ss-notation}(\ref{ss-n-log-lift})). 
  Since $\Ex(f) \subseteq \Supp \Delta_Y$ (see Lemma \ref{l-coeff-1/3}),   
  $(X, \Delta)$ is log liftable. Thus (1) holds. 
\end{proof}

\subsection{Weak log liftability for log del Pezzo pairs}\label{ss-weak-log-lift}

The purpose of this subsection is to show the following theorem. 

\begin{theorem} \label{thm:liftability-of-Hiromu-resolution} 
Assume $p>5$. 
Let $(X, \Delta)$ be a log del Pezzo pair. 
Suppose that there exist a birational morphism $f \colon Y \to X$ 
from a smooth projective surface $Y$ 
and an effective simple normal crossing $\Q$-divisor $B_Y$ on $Y$ 
such that $f_*B_Y =\Delta$ and $-(K_Y+B_Y)$ is nef and big. 
Then $(Y, \Supp B_Y)$ lifts to $W(k)$. 
\end{theorem}

\begin{proof}
By Remark \ref{r-log-dP-summary}, we may run a $(-K_Y)$-MMP: $g': Y \to Z'$. 
Then also $Z'$ is of del Pezzo type and $-K_{Z'}$ is nef and big, because 
$-K_Y = -(K_Y+B_Y) + B_Y$ is big. 
As $-K_{Z'}$ is semi-ample (Remark \ref{r-log-dP-summary}), there exists a birational morphism $h : Z' \to Z$ 
to a projective normal surface $Z$ such that $h_*\MO_{Z'} = \MO_Z$, $-K_Z$ is ample, and $-K_{Z'} = h^*(-K_Z)$. 
We then have the induced morphisms: 
\[
g: Y \xrightarrow{g'} Z' \xrightarrow{h} Z. 
\]
We can write 
\[
K_Y+
{\cred \Gamma} 
= (g')^*K_{Z'} = g^*K_Z, 
\]
for some  effective $g'$-exceptional $\Q$-divisor ${\cred \Gamma}$. 
It holds that 
\begin{equation} \label{eq:Delta-in-BY}
\Supp {\cred \Gamma} \subseteq \Supp B_Y, 
\end{equation}
because the MMP $g' : Y \to Z'$ is a $B_Y$-MMP by  $-K_Y = -(K_Y + B_Y)+B_Y$, and hence 
any prime divisor contracted by $g'$ must be an irreducible component of $B_Y$.

Set $E$ to be the maximum reduced divisor such that $\Supp E \subseteq \Ex(g)$ and $\Supp E \subseteq \Supp B_Y$. 
We then obtain 
\[
\kappa(Y, K_Y+E) \leq \kappa(Z, g_*(K_Y+E)) = \kappa(Z, K_Z) = -\infty. 
\]
Hence, as $(Y,E)$ is log smooth, 
this pair lifts to $W(k)$ by $p>5$ and \cite[Theorem 1.3]{Kaw3}.
We denote this lift by $(\mathcal{Y}, \mathcal{E})$.

\begin{claim}\label{claim:liftability-of-Hiromu-resolution} 
Every prime divisor $C \subseteq \Supp B_Y$  lifts to $\mathcal{Y}$. 
\end{claim}

\begin{proof}[Proof of Claim \ref{claim:liftability-of-Hiromu-resolution}] 
Fix a prime divisor $C$ with $C \subseteq \Supp B_Y$. 
If $C \subseteq E$, then there is nothing to prove. 
So we may assume that $C \not \subseteq \Supp E$, that is, $C$ is not $g$-exceptional. 
In particular, $C \not\subseteq \Supp {\cred \Gamma}$. 

First, since $H^2(Y, \mathcal{O}_Y)=0$, it follows from \cite[Corollary 8.5.6 (a)]{fga2005} 
that the line bundle $L=\mathcal{O}_Y(C)$ lifts to a line bundle $\mathcal{L}$ on $\mathcal Y$. 
Second, it holds that $H^i(Y, L) =0$ for $i>0$ by Kawamata-Viehweg vanishing (which holds for every boundary by \cite[Theorem 1.1]{ABL20}), because we can write 
\begin{align*}
C &= K_Y + C + {\cred \Gamma}- (K_Y + {\cred \Gamma}) \\
&= \underbrace{K_Y + (1-\delta)C + ({\cred \Gamma}+\epsilon A)}_{\mathrm{klt}} 
+ 
\underbrace{(g^*(-K_Z) - \epsilon A) + \delta C}_{\mathrm{ample}}, 
\end{align*}
where $A$ is a $g$-anti-ample $g$-exceptional effective $\Q$-divisor and $0 < \delta \ll \epsilon \ll 1$. Here, $(Y, (1-\delta)C + ({\cred \Gamma}+\epsilon A))$ is klt, because 
\begin{itemize}
\item $(Y, \Supp (C + {\cred \Gamma}))$ is simple normal crossing as $\Supp (C+{\cred \Gamma}) \subseteq \Supp B_Y$ (cf.\ (\ref{eq:Delta-in-BY})), therefore
\item $(Y, C + {\cred \Gamma})$ is plt given that $\rdown{{\cred \Gamma}}=0$, and so
\item $(Y, C + {\cred \Gamma} + \epsilon A)$ is plt, because $C \not \subseteq \Supp A$ and a small perturbation of a plt pair by any effective $\bQ$-divisor not containing $C$ is plt.
\end{itemize}


By upper semicontinuity (\cite[Chapter III, Theorem 12.18]{hartshorne77}), we have $H^i(\mathcal{Y}, \mathcal{L})=0$ for $i>0$, and so by Grauert's theorem (\cite[Chapter III, Corollary 2.19]{hartshorne77}), the restriction map 
\[
H^0(\mathcal{Y}, \mathcal{L}) \to H^0(Y, \mathcal L|_Y) = H^0(Y, L) 
\]
is surjective. 
Then there exists an effective Cartier divisor $\mathcal C$ on $\mathcal Y$ such that 
$\mathcal C|_Y =C$ and $\MO_{\mathcal Y}(\mathcal C) \simeq \mathcal L$. 
Then $\mathcal C$ is flat over $W(k)$, because 
the fibre $Y$ over the closed point $pW(k) \in \Spec W(k)$ is irreducible 
and $\mathcal C$ does not contain $Y$. 
This complete the proof of Claim \ref{claim:liftability-of-Hiromu-resolution}. 
\end{proof}

Apply Claim \ref{claim:liftability-of-Hiromu-resolution} to every prime divisor $C \subseteq \Supp B_Y$. 
By the fact that log-smoothness deforms (see \cite[Remark 2.7]{Kaw3}), 
we get that $(Y, \Supp B_Y)$  lifts to $W(k)$. 
This completes the proof of Theorem \ref{thm:liftability-of-Hiromu-resolution}.
\end{proof}

%% file: section5.tex
\section{Log del Pezzo pairs in characteristic $p>41$}\label{s-LDP-QFS42}

Throughout this section, we work over an algebraically closed field of characteristic $p>0$. 
The purpose of this section is to show that $(X, \Delta)$ is quasi-$F$-split 
if $p>41$ and $(X, \Delta)$ is a log del Pezzo pair with standard coefficients (Theorem \ref{t-LDP-QFS}). 

\subsection{An explicit ACC for log Calabi-Yau surfaces}\label{ss-explicit-ACC}

In this subsection, we establish the following vanishing theorem: 
\[
H^0(X, \cO_X(K_X + \Delta_{\red} + \rdown{p^{\ell}(K_X+\Delta)})) =0
\]
for a log del Pezzo pair $(X, \Delta)$ with standard coefficients in characteristic $p>41$  (Theorem \ref{t-42-vanishing}) . 
To this end, we prove an explicit version of ACC for two-dimensional log Calabi--Yau pairs (Theorem \ref{t-ACC2}). 
We start by treating the following essential case.

\begin{theorem} \label{t-ACC1}
Assume $p>5$. 
Let $(X,\Delta = aC + B)$ be a log del Pezzo pair with standard coefficients, 
where $a \geq 0$, $C$ is a prime divisor, and $B$ is an effective $\Q$-divisor 
with $C \not\subseteq \Supp B$. 

Suppose that $(X, tC+B)$ is lc and $K_X + tC + B \equiv 0$ for a real number $t$ with 
$0 \leq t < 1$. Then $t \leq \frac{41}{42}$. 
\end{theorem}
{In characteristic zero, this theorem is a special case of \cite[Definition 1.1, (5.1.2), and Theorem 5.3]{Kol94}. In what follows, we deduce it in positive characteristic by constructing an appropriate lift.}

\begin{proof}
We have $t \in \Q$, because the equality $(K_X + tC + B) \cdot H =0$ holds for an ample Cartier divisor $H$. 
It is clear that $a < t < 1$. 
Moreover, we may assume that $t>\frac{5}{6}$, as otherwise there is nothing to show. 
Note that $-(K_X+sC+B)$ is ample for any $a \leq s <t$, 
because $K_X + tC + B \equiv 0$ and $-(K_X+aC+B)$ is ample. 

\begin{claim}\label{c-ACC1}
There exists a resolution $f : Y \to X$ of $X$ such that: 
\begin{enumerate}
\item[(i)] 
$\Supp f^*C \subseteq \Supp (tC_Y + B_Y)$ and 
$B_Y$ is effective and simple normal crossing, where 
$C_Y := f^{-1}_*C$ and $B_Y$ is the $\bQ$-divisor defined by
\[
K_Y + tC_Y + B_Y = f^*(K_X + tC +B).
\]
\item[(ii)] $f: Y \to X$ is a log resolution of $(X, tC+B)$ over some open subset $X'$ containing $C$. 
\end{enumerate}
\end{claim}


\begin{proof}[Proof of Claim \ref{c-ACC1}]
Since the problem is local on $X$, we fix a closed point $x$ of $X$ around which we shall work. 
It suffices to find an open neighbourhood $\widetilde X$ of $x \in X$ 
and a resolution $\widetilde f: \widetilde{Y} \to \widetilde{X}$ of $\widetilde X$ which satisfy the corresponding properties (i) and (ii). 
If $x \not\in C$, then Corollary \ref{cor:Hiromu-resolution} enables us to find 
such a resolution $\widetilde f: \widetilde{Y} \to \widetilde{X}$ 
for $\widetilde X := X \setminus C$, because $(t C +B)|_{\widetilde X} =B|_{\widetilde X}$ has standard coefficients. 
Then we may assume that $x \in C$ and $(X, tC+B)$ is not log smooth at $x$. 
If $x$ is a singular point of $X$, then we may apply Proposition \ref{p-klt-sing} to $(X,(t-\epsilon_1)C+B)$ and $0 < {\epsilon_1}  \ll 1$ (here we take this $\epsilon_1$ to ensure that $\Supp f^*C \subseteq \Supp (tC_Y + B_Y)$). 
Assume that $X$ is smooth at $x$.  
For $0 < \epsilon_2 \ll 1$, we obtain 
\[
\mult_x ((t-\epsilon_2)C + B) \geq 
\mult_x \left( \frac{5}{6}C + B\right) 
\geq 
\min\left\{ 2 \cdot \frac{5}{6}, \hspace{2mm}\frac{5}{6} + \frac{1}{2} \right\} = \frac{4}{3}. 
\]
Therefore, we can apply Corollary \ref{c-klt-resol-r=1'} to $(X,(t-\epsilon_2)C+B)$. 
This completes the proof of Claim \ref{c-ACC1}.
\end{proof}

Pick $\epsilon \in \bQ$ such that $0 < \epsilon \ll 1$. 
We have 
\[
K_Y+tC_Y + B_Y - \epsilon f^*C = f^*(K_Y+(t-\epsilon)C + B). 
\]
Note that $tC_Y + B_Y - \epsilon f^*C$ is effective and simple normal crossing by Claim \ref{c-ACC1}.
Since $(X, (t-\epsilon)C+B)$ is a log del Pezzo pair, we may apply Theorem \ref{thm:liftability-of-Hiromu-resolution} 
to $f : Y \to X$ and $tC_Y + B_Y - \epsilon f^*C$, 
so that $(Y, \Supp (tC_Y + B_Y - \epsilon f^*C))$ lifts to $W(k)$. As $\Supp (f^*C) \subseteq \Supp (tC_Y + B_Y)$, we immediately get that $(Y,tC_Y +B_Y)$  admits a lift 
$(\mathcal Y, t \mathcal C_{\mathcal Y} + \mathcal B_{\mathcal Y})$ to $W(k)$.

For $\Gamma := tC_Y + B_Y - \epsilon f^*C$, 
$B^{\text{st}}_{Y} := f_*^{-1}B$, and $B^{\text{neg}}_Y := B_Y - B^{\text{st}}_{Y}$, the following holds on $Y$. 
\begin{enumerate}
\item $(Y, \Gamma)$ is a log smooth weak del Pezzo pair. 
\item $K_Y +t C_Y +B_Y \equiv 0$. 
\item 
$B_Y = B^{\text{st}}_{Y} + B^{\text{neg}}_Y$, 
where $B^{\text{st}}_{Y}$ has standard coefficients and 
$B^{\text{neg}}_Y$ is an effective $\Q$-divisor 
which is negative definite, that is, the intersection matrix of $\Supp B^{\text{neg}}_Y$ is negative definite. 
\end{enumerate}
Via the lift $(\mathcal Y, t \mathcal C_{\mathcal Y} + \mathcal B_{\mathcal Y})$, 
the geometric generic fibre 
$Y_{\overline K}$ of $\mathcal Y \to \Spec W(k)$ satisfies the following corresponding properties (1)'--(3)', 
where $\overline K$ denotes the algebraic closure of ${\rm Frac}\,W(k)$ and 
$\Gamma_{\overline K}, C_{Y_{\overline K}}, B_{Y_{\overline K}}$ are the $\Q$-divisors corresponding to 
$\Gamma, C_Y, B_Y$, respectively. 
\begin{enumerate}
\item[(1)'] $(Y_{\overline K}, \Gamma_{\overline K})$ is a log smooth weak del Pezzo pair. 
\item[(2)'] $K_{Y_{\overline K}} + tC_{Y_{\overline K}} + B_{Y_{\overline K}} \equiv 0$. 
\item[(3)'] 
$B_{Y_{\overline K}} = B^{\text{st}}_{Y_{\overline K}} + B^{\text{neg}}_{Y_{\overline K}}$, 
where $B^{\text{st}}_{Y_{\overline K}}$ has standard coefficients and 
$B^{\text{neg}}_{Y_{\overline K}}$ is an effective $\Q$-divisor 
which is negative definite. 
\end{enumerate}
Since $Y_{\overline K}$ is of del Pezzo type by (1)' and Remark \ref{r-log-dP-summary}, 
we can run a $B_{Y_{\overline K}}^{\text{neg}}$-MMP (Remark \ref{r-log-dP-summary}): 
\[
\psi : Y_{\overline K} \to V, 
\]
where $V$ denotes the end result. 
Since $B^{\text{neg}}_{Y_{\overline K}}$ is negative definite, 
its push-forward $\psi_*B^{\text{neg}}_{Y_{\overline K}}$ is either zero or negative definite. 
{Indeed, a $\Q$-divisor $D$ on a projective $\Q$-factorial surface is negative definite if and only if $D'^2 <0$ for every nonzero $\Q$-divisor $D'$ satisfying $\Supp D' \subseteq \Supp D$.} 


As $\psi_*B^{\text{neg}}_{Y_{\overline K}}$ is nef, 
we obtain $\psi_*B^{\text{neg}}_{Y_{\overline K}} =0$. 
By (3)', the $\Q$-divisor $B_V := \psi_*B_{Y_{\overline K}} = \psi_*B^{\text{st}}_{Y_{\overline K}}$ has standard coefficients. It follows from (2)' that 
$K_V + tC_V +B_V \equiv 0$ for $C_V := \psi_*C_{Y_{\overline K}}$. 
Note that $C_V$ is still a prime divisor, because 
the above $B^{\text{neg}}_{Y_{\overline K}}$-MMP only contracts the curves contained in $\Supp B^{\text{neg}}_{Y_{\overline K}}$. 
It follows from \cite[Definition 1.1, (5.1.2), and Theorem 5.3]{Kol94} that $t \leq \frac{41}{42}$. 
\qedhere

\end{proof}

\begin{lemma}\label{l-curve-ACC}
Let $(C, \Delta)$ be a one-dimensional projective lc pair such that $K_C+\Delta \equiv 0$ and 
\[
\coeff(\Delta) \subseteq \left\{ \frac{n-1}{n} \, \middle|\, n \in \Z_{>0} \right\} \cup \left[\frac{5}{6}, 1\right]. 
\]
Then 
\[
\coeff(\Delta) \subseteq 
\left\{ \frac{n-1}{n} \, \middle|\, n \in \Z_{>0}, 1 \leq n \leq 6\right\} \cup  \{1\}. 
\]
\end{lemma}

\begin{proof}
Let $\Delta = \sum_{i=1}^r a_iP_i$ be the irreducible decomposition. 
We may assume that $\Delta \neq 0$. 
It follows from $K_C+ \Delta \equiv 0$ that $C =\P^1$ and $\sum_{i=1}^r a_i =2$. 
By contradiction, we suppose that $\frac{5}{6} < a_1 <1$. 
We then get $r=3$, which implies $a_1 + a_2 +a_3 =2$. 
It is easy to see that there does not exist such a solution. 
\end{proof}


\begin{theorem} \label{t-ACC2}
Assume $p>5$. 
Let $(X,\Delta)$ be a two-dimensional projective lc pair 
such that $K_X + \Delta \equiv 0$ and 
\[
\coeff(\Delta) \subseteq \left\{ \frac{n-1}{n} \, \middle|\, n \in \Z_{>0} \right\} \cup \left[\frac{41}{42}, 1\right]. 
\]
Then 
\[
\coeff(\Delta) \subseteq \left\{ \frac{n-1}{n} \, \middle|\, n \in \Z_{>0}, 1 \leq n \leq 42\right\} \cup  \{1\}. 
\]
\end{theorem}

\begin{proof}
Replacing $(Y, \Delta_Y)$ by $(X, \Delta)$ for a dlt blowup $f: Y \to X$ of $(X, \Delta)$ (see, for example, \cite[Theorem 4.7]{tanaka16_excellent}), 
we may assume that $(X, \Delta)$ is dlt. 
Take the irreducible decomposition of $\Delta$ and define effective $\Q$-divisors $D, E, F$ as follows: 
\[
\Delta =D+E+F = \sum_{i=1}^r d_i D_i + \sum_{j=1}^s e_j E_j + \sum_{k=1}^t F_k, 
\]
\[
D:= \sum_{i=1}^r d_i D_i, \quad E:=\sum_{j=1}^s e_j E_j,\quad F:=\sum_{k=1}^t F_k. 
\]
where $D_i, E_j, F_k$ are prime divisors, 
\[
\frac{1}{2} \leq d_i \leq \frac{41}{42}, \qquad \text{and}\qquad 
\frac{41}{42} < e_j <1
\]
for all $i, j, k$. 
It suffices to show $E =0$. 
Suppose $E \neq 0$. 
Let us derive a contradiction. 

We run a $(K_X+D+F)$-MMP $g: X \to Z$, 
which is a $(-E)$-MMP by 
$-E \equiv K_X+D+F$. 
Then the end result $Z$ has a Mori fibre space structure $\pi: Z \to T$ such that 
$g_*E$ is $\pi$-ample. 
In particular, $g_*E \neq 0$. 
Furthermore, $(Z, g_*D)$ is klt, because 
$(X, D+(1-\epsilon)F)$ is klt and this MMP can be considered as a $(K_X+D+(1-\epsilon)F)$-MMP for some $0 < \epsilon \ll 1$. 
Replacing $(Z, g_*\Delta)$ by $(X, \Delta)$, we may assume, 
in addition to the original assumptions, that (1)--(3) hold. 
\begin{enumerate}
    \item $(X, D)$ is klt. 
    \item $X$ has a $(K_X+D+F)$-Mori fibre space structure $\pi: X \to T$. 
    \item $E$ is $\pi$-ample. 
\end{enumerate}

Assume $\dim T =1$. 
By $p\geq 5$, any $\pi$-horizontal irreducible component of $\Delta$ 
is generically \'etale over $T$. 
Hence, 
we obtain  $\coeff(\Delta|_W) \subseteq \coeff(\Delta)$ 
for a general fibre $W$ of $\pi: Z \to T$. 
By (3), $E|_F \neq 0$. 
This contradicts Lemma \ref{l-curve-ACC}. 
\medskip

Assume $\dim T =0$. {In particular, $\rho(X)=1$.}
Since the bound for ACC for log canonical thresholds is $\frac{5}{6}$, 
$(X, D + \sum_{j=1}^s E_j +\sum_{k=1}^t F_k)$ is lc
\footnote{{\cred this fact can be proven as follows (cf.\ \cite[Section 5]{HMX14}). 
If $(X, \Delta =D+E+F)$ is not klt, then 
extract a prime divisor $C$ whose discrepancy $a(C, X, \Delta)$ is $-1$, and use the adjunction restricting to $C$, so that 
Lemma \ref{l-curve-ACC} is applicable. 
If $(X, \Delta =D+E+F)$ is klt, then 
we can apply the proof for the non-klt case after increasing the coefficients $e_1, ..., e_s$ of $E$ one by one}}. 
For a sufficiently large integer $\ell \gg 42$, the pair 
\[
\left(X, D+  \frac{\ell -1}{\ell}\sum_{j=1}^s E_j + \frac{\ell -1 }{\ell}\sum_{k=1}^t F_k\right)
\]
is klt and $K_X +D + \frac{\ell -1}{\ell}\sum_{j=1}^s E_j + \frac{\ell -1 }{\ell}\sum_{k=1}^t F_k$ is ample. 
On the other hand, 
\[
\left(X, D+  \frac{41}{42}\sum_{j=1}^s E_j + \frac{\ell -1 }{\ell}\sum_{k=1}^t F_k\right)
\]
is klt and $-\left(K_X +D +\frac{41}{42}\sum_{j=1}^s E_j + \frac{\ell -1 }{\ell}\sum_{k=1}^t F_k\right)$ is ample.
Enlarge the coefficient of $E_1$ from $\frac{41}{42}$ to $\frac{\ell-1}{\ell}$. 
Then 
\[
K_X + D+ \frac{\ell -1 }{\ell} E_1 + \frac{41}{42}\sum_{j=2}^s E_j + \frac{\ell -1 }{\ell}\sum_{k=1}^t F_k 
\]
is either nef or anti-ample. 
If this $\Q$-divisor is nef, then there exists $u \in \Q$ 
such that $\frac{41}{42} < u \leq \frac{\ell -1 }{\ell}$ and 
$K_X + D+ u E_1 + \frac{41}{42}\sum_{j=2}^s E_j + \frac{\ell -1}{\ell}\sum_{k=1}^t F_k \equiv 0$. 
However, this contradicts Theorem \ref{t-ACC1}. 
Hence $-(K_X + D+ \frac{\ell -1 }{\ell} E_1 + \frac{41}{42}\sum_{j=2}^s E_j + \frac{\ell -1 }{\ell}\sum_{k=1}^t F_k)$ is ample. 
By enlarging the coefficient of $E_2$ from $\frac{41}{42}$ to $\frac{\ell-1}{\ell}$, the same argument deduces that 
\[
K_X + D+ \frac{\ell -1 }{\ell} (E_1+E_2) + \frac{41}{42}\sum_{j=3}^s E_j + \frac{\ell -1 }{\ell}\sum_{k=1}^t F_k 
\]
is anti-ample. By repeating this procedure finitely many times, 
we get that 
\[
K_X + D+ \frac{\ell -1 }{\ell}\sum_{j=1}^s E_j + \frac{\ell -1 }{\ell}\sum_{k=1}^t F_k 
\]
is anti-ample. 
However, this is again a contradiction. 
Hence $E=0$, as required. 
\end{proof}

\begin{theorem} \label{t-42-vanishing}
Assume $p>41$. 
Let $(X,\Delta)$ be a log del Pezzo pair with standard coefficients. Then
\[
H^0(X, \cO_X(K_X +  \Delta_{\red} + \rdown{p^{\ell}(K_X+\Delta)})) =0
\]
for 
every positive integer $\ell$.
\end{theorem}
It is essential to assume that $p>41$. Without this assumption, the theorem is false (cf.\ Lemma \ref{l-P^2-B_1}(2)).
\begin{proof}
Suppose $H^0(X, \cO_X(K_X +  \Delta_{\red} + \rdown{p^{\ell}(K_X+\Delta)})) \neq 0$ 
for some $\ell \in \Z_{>0}$. 
Let us derive a contradiction. 
We have 
\[
K_X + \Delta_{\red}+ \llcorner p^{\ell}(K_X+\Delta) \lrcorner \sim N
\]
for some effective Weil divisor $N$. 
By running a $K_X$-MMP and replacing divisors by their push-forwards, 
we may assume that there is a $K_X$-Mori fibre space structure $\pi : X \to T$. 
Take the irreducible decomposition of $\Delta$ and define $D$ and $E$ as follows: 
\[
\Delta  = D+E, 
\qquad 
D := \sum_{i=1}^r \frac{d_i -1}{d_i} D_i,\qquad \text{and}\qquad
E := \sum_{j=1}^s \frac{e_j-1}{e_j} E_j, 
\]
where $D_i$ and $E_j$ are prime divisors, 
\[
2 \leq d_i \leq 42, \quad \text{and}\quad 
e_j \geq 43 \quad \text{for all}\quad i, j. 
\]
The following holds: 
\begin{align*}
 \Delta_{\red} + \rdown{p^{\ell} \Delta} 
&= \left( \sum_{i=1}^r D_i + \sum_{j=1}^s E_j\right) + 
\rdown{ p^{\ell}\left(\sum_{i=1}^r \frac{d_i -1}{d_i} D_i + \sum_{j=1}^s \frac{e_j-1}{e_j} E_j\right) }\\
&\leq  \left( \sum_{i=1}^r D_i + \sum_{j=1}^s  E_j\right) + 
 p^{\ell}\left(\sum_{i=1}^r \frac{d_i -1}{d_i} D_i + \sum_{j=1}^s \frac{e_j-1}{e_j} E_j\right)  - \sum_{i=1}^r \frac{1}{d_i} D_i\\
 &= (p^{\ell}+1)\left( D+ \frac{p^{\ell} E + E_{\red}}{p^{\ell}+1}\right). 
\end{align*}
Therefore, we obtain 
\[
0 \leq N \sim 
K_X +  \Delta_{\red} + \rdown{p^{\ell}(K_X+\Delta)}
\leq (p^{\ell}+1)\left( K_X+ D+ \frac{p^{\ell} E + E_{\red}}{p^{\ell}+1}\right).
\]
In particular, we obtain $E \neq 0$. 
Since  $-(K_X+D+E)$ is $\pi$-ample and 
\[
K_X+ D+ \frac{p^{\ell} E + E_{\red}}{p^{\ell} +1}
\]
is $\pi$-nef, 
we can find a $\Q$-divisor $E'$ such that 
\[
E \leq E' \leq \frac{p^{\ell} E + E_{\red}}{p^{\ell}+1}\qquad \text{and}\qquad 
K_X+D+E' \equiv_{\pi} 0. 
\]

If $\dim T=0$, then we get  $E'=0$ by Theorem \ref{t-ACC2}. 
However, this implies $E=0$, which is a contradiction. 
Assume $\dim T=1$.  
Pick a general fibre $W \simeq \P^1$ of $\pi : X \to T$. 
By $p\geq 5$, every $\pi$-horizontal prime divisor contained in $\Delta$ is generically \'etale over $T$. 
Hence it holds that  $\coeff (\Delta|_W) \subseteq \coeff(\Delta)$ and every coefficient of $E|_W$ is $\geq \frac{42}{43}$. 
We then get  $E'=0$ by Lemma \ref{l-curve-ACC}, which is a contradiction. 
\end{proof}

\subsection{Log del Pezzo pairs are quasi-F-split in characteristic $p>41$}\label{ss-LDP>41}


\begin{theorem} \label{t-BS-vanishing} 
Let $(X, \Delta)$ be a log del Pezzo pair. 
Let  $\varphi \colon V \to X$ be a birational morphism from a smooth projective surface $V$ and let $(V, \Delta_V)$ be a klt pair such that $\Delta_V$ is simple normal crossing, $-(K_V+\Delta_V)$ is nef, $\varphi_*\Delta_V = \Delta$, and $(V, \Supp \Delta_V)$ lifts to $W_2(k)$.

Let $\theta \colon W \to V$ be a birational morphism from a smooth projective surface $W$ such that the composition 
\[
\psi := \varphi \circ \theta \colon W \xrightarrow{\theta} V \xrightarrow{\varphi} X
\]
is a log resolution of $(X, \Delta)$. Write $K_W + \Delta_W = \theta^*(K_V+\Delta_V)$ and 
let $\bE_W$ be the union of all the $\psi$-exceptional prime divisors $E$ 
such that $\lfloor K_W + \Delta_W \rfloor \cdot E = 0$ and 
that the coefficient of $E$ in $\Delta_W$ is an integer. 
Assume that $p$ does not divide the determinant of the intersection matrix of $\bE_W$. 

Then 
\[
H^0(W, \Omega^1_W(\log  D) \otimes \MO_W( \lfloor K_W + \Delta_W \rfloor))=0
\] 
for the reduced divisor $D$ on $W$ satisfying $\Supp D = 
\mathrm{Exc}(\psi) \cup
\Supp \psi^{-1}_* \Delta$.
\end{theorem}

\begin{proof}
Recall that 
$\Omega^1_W(\log  D)( \lfloor K_W + \Delta_W \rfloor) = \Omega^1_W(\log  D) \otimes \MO_W( \lfloor K_W + \Delta_W \rfloor)$ under our notation. 
\begin{claim}\label{c-BS-vanishing}
Let $S$ be the sum of all the $\psi$-exceptional prime divisors whose coefficients in $\Delta_W$ are integers. Then 
\[
H^0(W, \Omega^1_W(\log  (D - S))( \lfloor K_W + \Delta_W \rfloor)) =
H^0(W, \Omega^1_W(\log  D)( \lfloor K_W + \Delta_W \rfloor)).
\]
\end{claim}
\vspace{0.3em}
The assertion  follows immediately from Claim \ref{c-BS-vanishing}. Indeed, since $\theta_{*}S$ contains all the $\varphi$-exceptional prime divisors whose coefficients in $\Delta_V$ are equal to zero, we have $\theta_{*}(D-S)\subseteq  \Supp \Delta_V$.
By $\theta_{*}\Delta_W=\Delta_V$, we have
\[
H^0(W, \Omega^1_W(\log (D - S))( \lfloor K_W + \Delta_W \rfloor)) \subseteq H^0(V, \Omega^1_V(\log  \Supp \Delta_V)( \lfloor K_V + \Delta_V \rfloor))=0,
\]
where the last equality follows from 
a nef-and-big Akizuki-Nakano vanishing \cite[Theorem 2.11]{Kaw3}.

\begin{proof}[Proof of Claim \ref{c-BS-vanishing}]
Let $E$ be a prime divisor contained in $\Supp S$.
Since $K_W+\Delta_W=\theta^{*}(K_V+\Delta_V)$ is anti-nef and 
the coefficient of $E$ in $\Delta_W$ is an integer, 
we have
\[
\rdown{K_W + \Delta_W} \cdot E\leq (K_W+\Delta_W)\cdot E\leq 0.
\]
First, we assume that $\rdown{K_W + \Delta_W} \cdot E<0$. Then, $H^0(E, \cO_E(\rdown{K_W+\Delta_W}))=0$, and so by tensoring the following short exact sequence
\[
0 \to \Omega^1_W(\log (D - E)) \to \Omega^1_W(\log  D) \to \mathcal{O}_E \to 0,
\]
with $\cO_W(\rdown{K_W+\Delta_W})$, 
we see that 
\[
H^0(W, \Omega^1_W(\log  (D - E))( \lfloor K_W + \Delta_W \rfloor)) = 
H^0(W, \Omega^1_W(\log  D)( \lfloor K_W + \Delta_W \rfloor)).
\]
Repeating this procedure, we can assume that every prime divisor $E$ in $S$ satisfies $\rdown{K_W + \Delta_W} \cdot E=0$, that is, that $S={\mathbb E_W}$.
Using essentially the same argument as in \cite[Subsection 8.C]{graf21}, we show that 
\[
H^0(W, \Omega^1_W(\log  (D - S))( \lfloor K_W + \Delta_W \rfloor)) = 
H^0(W, \Omega^1_W(\log  D)( \lfloor K_W + \Delta_W \rfloor)).
\]
By looking at the short exact sequence
\[
0 \to \Omega^1_W(\log  (D - S)) \to \Omega^1_W(\log  D) \to \bigoplus_{E \subseteq S} \cO_E \to 0, 
\]
tensored with $\cO_W(\lfloor K_W + \Delta_W \rfloor)$, it suffices to show that the homomorphism 
\begin{equation} \label{e1-BS-vanishing}
\overbrace{\bigoplus_{E \subseteq S} H^0(E, \mathcal{O}_E(\lfloor K_W + \Delta_W \rfloor))}^{=\ \bigoplus_{E \subseteq S} H^0(E, \mathcal{O}_E)} \to 
H^1(W, \Omega^1_W(\log (D -S))(\lfloor K_W + \Delta_W \rfloor))
\end{equation}
is injective. 

We note that $D-S$ is disjoint from $S$. 
Indeed, if there exists a prime divisor $E \subseteq S$ intersecting $D-S$, then 
$\lfloor K_W + \Delta_W \rfloor \cdot E <0$, because 
the prime divisors in $D-S$ have non-integral coefficients in $\Delta_W$. This is impossible by our construction of $D$.
Now, since $D-S$ is disjoint from $S$,  the map (\ref{e1-BS-vanishing}) factors through
\begin{equation} \label{e2-BS-vanishing}
\bigoplus_{E \subseteq S} H^0(E, \mathcal{O}_E) \to \bigoplus_{E \subseteq S} H^1(E, \Omega^1_E)
\end{equation}
via restriction $\Omega^1_W( \log (D - S))(\lfloor K_W + \Delta_W \rfloor) \to \bigoplus_{E \subseteq S} \Omega^1_E$.

Map (\ref{e2-BS-vanishing}) is given by the intersection matrix of $S$ (as explained in \cite[Subsection 8.C]{graf21}, this follows from \cite[Lemma 3.5, Fact 3.7, and Lemma 3.8]{gk14} which are essentially characteristic-free, just replace $\bC$ by the base field $k$), and so it is an isomorphism 
by the assumption that the determinant of the intersection matrix of $S={\mathbb{E}_W}$ is not divisible by $p$. In particular, 
map (\ref{e1-BS-vanishing}) is injective.
\end{proof}
\end{proof}

\begin{lemma}\label{l-log-lift-enough}
Assume $p>41$. 
Let $(X,\Delta)$ be a log del Pezzo pair with standard coefficients. 
Suppose that $(X,\Delta)$ is log liftable.
Then $(X, \Delta)$ is quasi-$F$-split.
\end{lemma}

\begin{proof}
Since $(X, \Delta)$ is log liftable, there exists a log resolution $f\colon Y\to X$ of $(X, \Delta)$ such that $(Y,f_{*}^{-1}\Delta+\Exc(f))$ lifts to $W(k)$. 
We define a $\Q$-divisor $\Delta_Y$ on $Y$ by 
\[
K_Y+\Delta_Y=f^{*}(K_X+\Delta).
\]
Since $X$ is $\Q$-factorial, 
we can find an effective $f$-exceptional $\Q$-divisor $F$ such that $-F$ is $f$-ample. 
For $0 < \epsilon \ll 1$, we set $B_Y := \Delta_Y + \epsilon F$. To prove that $(X, \Delta)$ is quasi-$F$-split, 
it suffices to show the following properties (A)--(C) 
(Theorem \ref{t-QFS-criterion}).  
\begin{enumerate}
\item[(A)] $\rdown{B_Y} \leq 0$, $f_*B_Y  =\Delta$, and $-(K_Y+B_Y)$ is ample. 
\item[(B)] $H^0(Y, \MO_Y(K_Y + (B_Y)_{\red} + p^{\ell}(K_Y+B_Y)) =0$ 
for any $\ell \in \Z_{>0}$. 
\item[(C)] $H^0(Y, \Omega_Y^1(\log\, (B_Y)_{\red}) \otimes \MO_Y(K_Y+B_Y))=0$. 
\end{enumerate}
The property (A) follows from the definition of $B_Y$. 
Since $(Y, (B_Y)_{\red})$ lifts to $W(k)$ and $-(K_Y+B_Y)$ is an  ample $\Q$-divisor satisfying $\Supp(\{K_Y+B_Y\})\subseteq \Supp (B_Y)_{\red}$, (C) follows from {\cite[Theorem 2.11]{Kaw3}}. 

It suffices to show (B). 
Since $f_{*}B_Y=\Delta$ and $f_{*}(B_Y)_{\red}= \Delta_{\red}$, we have 
\[
H^0(Y, \sO_Y(K_Y+(B_Y)_{\red}+p^{\ell}(K_Y+B_Y)))\subseteq 
H^0(X, \sO_X(K_X+\Delta_{\red}+p^{\ell}(K_X+\Delta))).
\]
Then the latter cohomology vanishes by Theorem \ref{t-42-vanishing}. 
\end{proof}

We are ready to prove the main theorem of this section.

\begin{theorem} \label{t-LDP-QFS}
Assume $p>41$. 
Let $(X,\Delta)$ be a log del Pezzo pair with standard coefficients. 
Then $(X,\Delta)$ is quasi-$F$-split.
\end{theorem}

\begin{proof}
Recall that a point $x \in X$ is called a \emph{special {point of}} $(X,\Delta)$ if 
Theorem \ref{t-hiromu-resol}(2) holds around $x$ 
(in particular, $X$ is smooth at $x$ and we have $\Delta = \frac{1}{2} C$ around $x$ for a prime divisor $C$). 
We consider two cases separately.\\

\noindent \textbf{Case 1.} 
For every special point $x$ of $(X,\Delta)$, 
$x$ is of type $a_x$ with $a_x < p$ (see Remark \ref{r-cusp-resol}). \\

{\cred 
By applying Theorem \ref{t-hiromu-resol} to each non-log smooth point of $(X, \Delta)$, 
we get  a resolution $\varphi \colon V \to X$  such that 
\begin{itemize}
\item $\Delta_V$ is effective and $(V, \Delta_V)$ is log smooth for the $\Q$-divisor $\Delta_V$ defined by $K_V+\Delta_V=\varphi^*(K_X+\Delta)$,  
\item over each special point $x$ of $(X, \Delta)$, $\varphi$ is obtained by 
a sequence of blowups  as in Theorem \ref{t-hiromu-resol}(2), and 
\item $\varphi$ is a log resolution over every point $x \in X$ except when $x$ is a special point {of} $(X, \Delta)$. 
\end{itemize}}
\noindent 
We construct  a birational morphism  $\theta \colon W \to V$ from a smooth projective surface ${\cred W}$ as follows. 
Over a non-special point of $(X, \Delta)$, $\theta$ is an isomorphism. 
Over a special point of $(X, \Delta)$, we blow up $X$ twice, so that the composition 
\[
\psi : W \xrightarrow{\theta} V \xrightarrow{\varphi} X
\]
is a log resolution of $(X, \Delta)$  (cf.\ Remark \ref{r-cusp-resol}). 

We define $\Q$-divisors $\Delta_V$ and $\Delta_W$ by 
\[
K_V+ \Delta_V=\varphi^*(K_X+\Delta)\qquad\text{and}\qquad 
K_W+\Delta_W = \theta^*(K_V+\Delta_V) (=\psi^*(K_X+\Delta)). 
\]
Pick an effective $\psi$-exceptional $\Q$-divisor $F$ on $W$ such that $-F$ is $\psi$-ample. 
Set $B_W := \Delta_W + \epsilon F$ for $0 < \epsilon \ll 1$. 
It is enough to show that the following conditions (A)--(C) hold (Theorem \ref{t-QFS-criterion}). 
\begin{enumerate}
\item[(A)] $\rdown{B_W} \leq 0$, $\psi_*B_W  =\Delta$, and $-(K_W+B_W)$ is ample. 
\item[(B)] $H^0(W, \MO_Y(K_W + (B_W)_{\red} + p^{\ell}(K_W+B_W))) =0$ 
for any $\ell \in \Z_{>0}$. 
\item[(C)] $H^0(W, \Omega_W^1(\log (B_W)_{\red}) \otimes \MO_W(K_W+B_W))=0$. 
\end{enumerate}
By construction, (A) holds. 
We have that 
{\cred 
\begin{eqnarray*}
&& H^0(W, \cO_W(K_W + (B_W)_{\red} + \rdown{p^{\ell}(K_W+B_W)}))\\
&\subseteq&  H^0(X, \cO_X(
\psi_*(K_W + (B_W)_{\red} + \rdown{p^{\ell}(K_W+B_W)})))\\
&=& H^0(X, \cO_X(K_X + \Delta_{\red} + \rdown{p^{\ell}(K_X+\Delta)}))=0
\end{eqnarray*}
}
for every $\ell \geq 1$ by Theorem \ref{t-42-vanishing}. 
Thus (B) holds. 

Let us show (C). 
By $\rdown{K_W+B_W} = \rdown{K_W+\Delta_W}$, we have 
 \[
H^0(W, \Omega^1_W(\log (B_W)_{\red})(\rdown{K_W+B_W}))\subseteq 
H^0(W, \Omega^1_W(\log D)(\rdown{K_W+\Delta_W})), 
\]
where $D$ is the reduced divisor with $\Supp D= \Supp \psi^{-1}_*\Delta + \Exc(\psi)$. 
Note that $(V, \Supp \Delta_V)$ lifts to $W_2(k)$ by Theorem \ref{thm:liftability-of-Hiromu-resolution}. 
By Remark \ref{r-cusp-resol}, we may apply Theorem \ref{t-BS-vanishing}, so that 
$H^0(W, \Omega^1_W(\log D)(\rdown{K_W+\Delta_W})) = 0$. 
Hence (C) holds. 
Thus $(X,\Delta)$ is quasi-$F$-split.\\

\noindent \textbf{Case 2.} There exists a special point $x$ of $(X,\Delta)$ of type $a_x$ with $a_x \geq p$. \\

Around $x$, we can write $\Delta  = \frac{1}{2}C$ for a prime divisor $C$ on $X$. 
We have $a_x \geq p \geq 43$. 
Let $\mu : X' \to X$ be the minimal resolution of $X$ and 
set $C_{X'} : =\mu_*^{-1}C$. 
The singularity $x$ of $C$ is resolved by $n_x$ blowups with $a_x = 2n_x +1$  (cf.\ Remark \ref{r-cusp-resol}). 
Until the singularity $x \in C$ is resolved, the multiplicity of the singular point lying over $x$ 
is equal to two. 
Since $C$ and $C_{X'}$ are isomorphic around $x$, 
we obtain $p_a(C_{X'}) \geq n_x$ \cite[Corollary V.3.7]{hartshorne77}. 
Then 
\[
p_a(C_{X'}) \geq n_x = \frac{a_x-1}{2} \geq \frac{43-1}{2} = 21 >18. 
\]
By Theorem \ref{t-bdd-or-liftable}, $(X, \Delta)$ is log liftable. 
Then Lemma \ref{l-log-lift-enough} implies that $(X, \Delta)$ is quasi-$F$-split. 
\end{proof}

%% file: section6.tex
\section{Log del Pezzo pairs which are not  quasi-$F$-split}\label{s-P^2-cex}

In this section, we prove that there exists a log del Pezzo pair $(X, \Delta)$ in characteristic $41$ such that $\Delta$ has standard coefficients and that $(X, \Delta)$ is not quasi-$F$-split (Theorem \ref{t-P^2-main}). Since our construction {can be applied to} some other characteristics, 
we start by {providing} a list of pairs $(X, \Delta)$ for which our argument works 
(Notation \ref{n-P^2-cex}). 
As explained in (\ref{sss-LDP-nonQFS}), {the} key part is to show that 
$h^1(X, B_m\Omega_X^1(\log \Delta_{\red})(p^{m}(K_X+\Delta))) =1$ (Proposition \ref{p-P^2-BZ-key}(5)). 
The following (Table \ref{Table:list of cohomologies}) is a list of the dimension of cohomologies that we will compute.

\begin{table}[h]
\caption{The following is a list of the dimension of cohomologies that we will compute. Here, we put $n, m \in \Z_{\geq 0}$, and $\Omega^{1}_{X,\log} := \Omega^{1}_{X}(\log \Delta_{\red})$. For the first column see: Lemma \ref{l-P^2-B_1} and Lemma \ref{l-P^2-H0Omega}. For the second column see Lemma \ref{l-P^2-H0Omega}, Lemma \ref{l-P^2-vanish}, and Lemma \ref{p-P^2-BZ-key}. For the last column see Lemma \ref{l-P^2-H0Omega}(3) and Lemma \ref{p-P^2-BZ-key}.}\label{Table:list of cohomologies}
     \centering
{\renewcommand{\arraystretch}{1.35}%
      \begin{tabular}{|l|c|c|c|}
      \hline
        & $h^0$ & $h^1$ & $h^{2}$ \\
       \hline
\multirow{2}{*}{$\Omega^{1}_{X}(-n)$} & \multirow{2}{*}{$0$} & 0 (if $n>0$) & \multirow{2}{*}{see (\ref{l-P^2-H0Omega}(3))} \\
&&1 (if $n=0$) &
       \\
       \hline
       $\Omega^{1}_{X} (p^{n} (K _{X}+ \Delta))$ & 0  &&
       \\
       \hline
$ \Omega^{1}_{X, \log}(p^{n}(K_{X}+ \Delta))$  & 0  & $0$ (if $n=1$)  &  0 (if $n=1$)  \\
\hline
\multirow{2}{*}{$Z_{m} \Omega^{1}_{X, \log}(p^{n}(K_{X}+ \Delta))$} & \multirow{2}{*}{0} & \hspace{0.4em} $0$ (if $n=m+1 \geq 2$)  & \\
&& \hspace{0.4em} 1 (if $n=m \geq 1$) \hspace{1.35em} &\\
\hline
\multirow{3}{*}{$B_{m} \Omega^{1}_{X,\log}(p^{n}(K_{X}+ \Delta))$}  &
\multirow{3}{*}{$0$}
 &
 \hspace{0.2em} $0$ (if $n=m+1 \geq 2$) \hspace{0.7em}
 &
 \\
 && \hspace{0.2em} $0$ (if $m=1$ and $n \geq 2$) &
 \\
 && \hspace{0.2em}  $1$ (if $n=m \geq 1$)  \hspace{2.35em} & \\
 \hline
 $B_{1}\Omega^{2}_{X,\log}(p(K_{X}+ \Delta))$ & 0  &&\\
 \hline
      \end{tabular}}
    \end{table}

\begin{notation}\label{n-P^2-cex}
{\cred 
Let $k$ be an algebraically closed field of characteristic $p>0$. 
Set $X := \mathbb P^2_k$ and 
let $L_0, L_1, L_2, L_3$ be distinct lines such that $L_0+L_1 + L_2 + L_3$ is simple normal crossing. 
{\cred Take positive integers $b_1, b_2, b_3$ satisfying 
\begin{itemize}
\item $\frac{1}{b_1}+\frac{1}{b_2}+\frac{1}{b_3} + \frac{1}{p+1} = 1$, and 
\item $p+1$ is divisible by $b_i$ for every $i \in \{ 1, 2, 3\}$.
\end{itemize}
}
\noindent Set 
\[
\Delta := \frac{p-1}{p} L_0 + \frac{b_1-1}{b_1}L_1+
\frac{b_2-1}{b_2}L_2+
\frac{b_3-1}{b_3}L_3 =:(p, b_1, b_2, b_3). 
\]
In what follows, for $n \in \Z$, $\MO_X(n)$ denotes the invertible sheaf on $X$ with $\deg \MO_X(n)=n$. 
Recall that $\Delta_{\red} = L_0 +L_1+L_2+L_3$.}
\end{notation}

{\cred 
The following is the list of $\Delta$ satisfying the conditions in Notation \ref{n-P^2-cex}.  
\begin{enumerate}
\renewcommand{\labelenumi}{(\roman{enumi})}
    \item $p=3$ and $\Delta = (3, 4, 4, 4)$.
    \item $p=5$ and $\Delta = (5, 3, 3, 6)$. 
    \item $p=5$ and $\Delta = (5, 2, 6, 6)$.
    \item $p=7$ and $\Delta = (7, 2, 4, 8)$.
    \item $p=11$ and $\Delta = (11, 2, 4, 6)$.
    \item $p=11$ and $\Delta = (11, 2, 3, 12)$
    \item $p=17$ and $\Delta =(17, 2, 3, 9)$.
    \item $p=19$ and $\Delta=(19, 2, 4, 5)$.
    \item $p=23$ and $\Delta = (23, 2, 3, 8)$.
    \item $p=41$ and $\Delta = (41, 2, 3, 7).$    
\end{enumerate}
}

\begin{proposition} \label{p-P^2-LDP}
{\cred 
We use the Notation \ref{n-P^2-cex}.
Then the following hold.
\begin{enumerate}
\item $-(K_X + \Delta)$ is ample.
\item $\cO_X(p (K_X + \Delta)) \simeq \cO_X(-1)$.
\end{enumerate}
}
\end{proposition}

\noindent
{\cred The following proof is due to Yuya Matsumoto.}

\begin{proof}
{\cred It holds that 
\begin{align*}
\deg \Delta &= 4 - \frac{1}{p}-\frac{1}{b_1}-\frac{1}{b_2}-\frac{1}{b_3} \\
&= 4 - \biggl( \frac{1}{b_1}+\frac{1}{b_2}+\frac{1}{b_3}+ \frac{1}{p+1} \biggr) - \biggl( \frac{1}{p} - \frac{1}{p+1} \biggr) \\
&= 3 - \frac{1}{p(p+1)} < 3.
\end{align*}
Thus (1) holds. 
Since $p+1$ is divisible by each $b_i$ (Notation \ref{n-P^2-cex}), 
we have 
$\lfloor -\frac{p}{b_i} \rfloor = -\frac{p+1}{b_i}$.
Then it holds that 
\begin{align*}
\deg \lfloor p \Delta \rfloor 
&= p (1 - \frac{1}{p}) + \sum_{i=1}^3 \lfloor p (1 - \frac{1}{b_i}) \rfloor  \\
&= 4 p - 1 + \sum_{i=1}^3 \lfloor -\frac{p}{b_i} \rfloor \\
&= 4 p - 1 - \sum_{i=1}^3 \frac{p+1}{b_i} \\
&\overset{(*)}{=} 4 p - 1 - (p+1) \biggl(1 - \frac{1}{p+1} \biggr) = 3 p - 1,
\end{align*}
where $(*)$ follows from $\frac{1}{p+1} + \sum_{i=1}^3 \frac{1}{b_i} = 1$  (Notation \ref{n-P^2-cex}). 
Thus  (2) holds.}
\end{proof}

\begin{lemma}\label{l-P^2-B_1}
We use Notation \ref{n-P^2-cex}. Then the following hold. 
\begin{enumerate}
    \item $(X, \{p^r \Delta \})$ is globally $F$-regular for all $r \in \Z_{>0}$.  
    \item $h^0(X, B_1\Omega_X^2(\log \Delta_{\red})(p(K_X+\Delta)))=1$. 
\end{enumerate}
\end{lemma}

\begin{proof}
Let us show (1). 
Fix $r \in \Z_{>0}$.  
By Notation \ref{n-P^2-cex}, there exists $a \in \Q$ such that 
$0 < a <1$ and we have $\{p^r \Delta \} \leq a(L_1 + L_2 + L_3)$. 
We may assume that $L_1, L_2, L_3$ are torus invariant divisors on 
a projective toric surface $X=\mathbb P^2_k$, 
and hence $(X, a(L_1+ L_2 + L_3))$ is globally $F$-regular. 
Therefore, also $(X, \{p^r \Delta \})$ is globally $F$-regular. 
Thus (1) holds. 

Let us show (2). 
By (\ref{exact1}), we have the following exact sequence: 
\begin{multline*}
0 \to B_1\Omega_X^2(\log \Delta_{\red})(p(K_X+\Delta)) \to F_*\Omega_X^2(\log \Delta_{\red})(p(K_X+\Delta)) \\
\overset{C}{\longrightarrow} 
\Omega_X^2(\log \Delta_{\red})(K_X+\Delta) \to 0. 
\end{multline*}
Hence, the equality $h^0(X, B_1\Omega_X^2(\log \Delta_{\red})(p(K_X+\Delta)))=1$ holds by  
\begin{align*}
H^0(X, \Omega_X^2(\log \Delta_{\red})(K_X+\Delta)) 
&=H^0(X, \Omega_X^2(\log \Delta_{\red})(K_X))\\
&= 
H^0(X, {\MO_X}(2K_X+ \Delta_{\red}))=0, \textrm{ and }\\
h^0(X, \Omega_X^2(\log \Delta_{\red})(p(K_X+\Delta)))
&=h^0(X, {\MO_X}(K_X+\Delta_{\red} + \llcorner p(K_X+\Delta)\lrcorner)) \\
&= h^0(X, \MO_X(-3+4-1))=1,   
\end{align*}
where $\MO_X(\rdown{p(K_X+\Delta)}) \simeq \MO_X(-1)$ follows from Proposition \ref{p-P^2-LDP}. 
Thus (2) holds. 
\end{proof}

\begin{lemma}\label{l-P^2-H0Omega}
We use Notation \ref{n-P^2-cex}. Then the following hold. 
\begin{enumerate}
\setlength{\itemsep}{6pt}
    \item $H^0(X ,\Omega_X^1 \otimes \MO_X(n))=0$ for all $n \leq 0$. 
    \item $\!\begin{aligned}[t]
    h^1(X ,\Omega_X^1 \otimes \MO_X(n))= 
    \begin{cases}
    0 \quad (n<0)\\
    1 \quad (n=0). 
    \end{cases}
    \end{aligned}$
    \item 
    $h^2(X,\Omega_X^1 \otimes \MO_X(n))= 3 h^0(X, \MO_X(-2-n)) - h^0(X, \MO_X(-3-n))$ for all $n \in \Z$. 
    \item $\!\begin{aligned}[t]
    &H^0(X, \Omega^1_X(p^n(K_X+\Delta))) = 0, \textrm{ and } \\ &H^0(X, \Omega^1_X(\log \Delta_{\red})(p^n(K_X+\Delta))) =0, \textrm{ for all } n \geq 0. 
    \end{aligned}
    $
    \item  
    $\!\begin{aligned}[t]
    &H^0(X, Z_m\Omega^1_X(\log \Delta_{\red})(p^n(K_X+\Delta))) = 0, \textrm{ and } \\
    &H^0(X, B_m\Omega^1_X(\log \Delta_{\red})(p^n(K_X+\Delta))) =0, \textrm{ for all } m\geq 0 \textrm{ and } n \geq 0. 
    \end{aligned}$
\end{enumerate}

\end{lemma}

\begin{proof}
The assertions (1)--(3) follow from the dual of the Euler sequence \cite[Ch. II, Theorem 8.13]{hartshorne77}: 
\[
0 \to \Omega_X^1 \to \MO_X(-1)^{\oplus 3} 
\to \MO_X \to 0. 
\]

Since (5) follows from (4), 
it is enough to show (4). 
By $\deg\,\llcorner p^n(K_X+\Delta)\lrcorner <0$ (Proposition \ref{p-P^2-LDP}(1)), 
it follows from (1) that $H^0(X, \Omega^1_X(p^n(K_X+\Delta))) =0$. 
We have an exact sequence: 
\[
0 \to \Omega_X^1 \to \Omega_X^1(\log \Delta_{\red}) \to \bigoplus_{i=1}^4 \MO_{L_i} \to 0. 
\]
Again by $\deg\,\llcorner p^n(K_X+\Delta)\lrcorner <0$, we have that 
\[
H^0(L_i, \MO_X(p^n(K_X+\Delta))|_{L_i}) =0
\]
for every $1 \leq i\leq 4$. 
Therefore, we get $H^0(X, \Omega^1_X(\log \Delta_{\red})(p^n(K_X+\Delta))) =0$. 
Thus (4) holds. 
\end{proof}

\begin{lemma}\label{l-P^2-vanish}
We use Notation \ref{n-P^2-cex}. 
Then the following hold. 
\begin{enumerate}
    \item $H^1(X, B_1\Omega_X^1(\log \Delta_{\red})(p^{s}(K_X+\Delta)))=0$ for all $s\geq 2$, and
    \item $H^1(X, B_m\Omega_X^1(\log \Delta_{\red})(p^{m+1}(K_X+\Delta)))=0$ for all $m \geq 1$. 
\end{enumerate}
\end{lemma}

\begin{proof}
Set $D :=K_X+\Delta$. Let us show (1). 
We have the following exact sequence (Lemma \ref{lem:Serre's map}): 
\[
0 \to \MO_X (p^{s-1}D ) \xrightarrow{F} F_* \MO_X(p^{s}D)  \to B_1\Omega_X^1(\log \Delta_{\red})(p^{s}D) \to 0. 
\]
By $H^1(X, F_* \MO_X(p^{s}D))=0$, it suffices to show that 
the $\MO_X$-module homomorphism 
\[
F: \MO_X (p^{s-1}D ) \to F_* \MO_X(p^{s}D) 
\]
splits. This is equivalent to the splitting of 
\[
F: \MO_X \to F_* \MO_X( p \{ p^{s-1} \Delta\}), 
\]
{given that $\rdown{p^sD}-p\rdown{p^{s-1}D} = \rdown{p(p^{s-1}D-\rdown{p^{s-1}D})} = \rdown{p\{p^{s-1}D\}} = \rdown{p\{p^{s-1}\Delta\}}$.}

This holds because $(X, \{ p^{s-1} \Delta\})$ is globally $F$-regular 
by Lemma \ref{l-P^2-B_1}(1) {(cf.\ \cite[Lemma 2.18]{KTTWYY1})}. 
Thus (1) holds. 

Let us show (2).  
We have the following exact sequence (\ref{exact2}):
\begin{multline*}
0 \to F_*^{m-1}(B_1\Omega_X^1(\log \Delta_{\red})(p^{m+1}D)) \to B_m\Omega_X^1(\log \Delta_{\red})(p^{m+1}D) \\
\overset{C}{\longrightarrow} B_{m-1}\Omega_X^1(\log \Delta_{\red})(p^{m}D) \to 0. 
\end{multline*}
By (1), it holds that $H^1(X, F_*^{m-1}(B_1{\Omega^1_X}(\log \Delta_{\red})(p^{m+1}D)))=0$ for any $m \geq 1$, 
which induces an injection: 
\[
H^1(X, B_m\Omega_X^1(\log \Delta_{\red})(p^{m+1}D)) \hookrightarrow H^1(X, B_{m-1}\Omega_X^1(\log \Delta_{\red})(p^{m}D)). 
\]
Therefore, we obtain a sequence of injections 
\begin{multline*}
H^1(X, B_m\Omega_X^1(\log \Delta_{\red})(p^{m+1}D)) \hookrightarrow H^1({X,} B_{m-1}\Omega_X^1(\log \Delta_{\red})(p^{m}D))
\hookrightarrow 
\\
\ \cdots\ \hookrightarrow
H^1(X, B_2\Omega_X^1(\log \Delta_{\red})(p^3D)) \hookrightarrow
H^1(X, B_1\Omega_X^1(\log \Delta_{\red})(p^2D))=0, 
\end{multline*}
where the last equality holds by (1). Thus (2) holds. 
\end{proof}

\begin{proposition}\label{p-P^2-BZ-key}
We use Notation \ref{n-P^2-cex}. 
Then the following hold. 
\begin{enumerate}
\setlength{\itemsep}{3pt}
    \item $H^j(X, \Omega_X^1(\log \Delta_{\red})(p(K_X+\Delta)))=0$,  for all $j \in \Z$. 
\item $\! \begin{aligned}[t]
H^j&(X, B_m\Omega_X^1(\log \Delta_{\red})(p^{m+1}(K_X+\Delta))) \\
&\simeq 
H^j(X, Z_m\Omega_X^1(\log \Delta_{\red})(p^{m+1}(K_X+\Delta))), \textrm{ for all } j \in \Z \textrm{ and } m \in \Z_{>0}.  
\end{aligned}$
\item 
$H^1(X, Z_m\Omega_X^1(\log \Delta_{\red})(p^{m+1}(K_X+\Delta)))=0$, 
for all $m \in \Z_{>0}$. 
\item 
$h^1(X, Z_m\Omega_X^1(\log \Delta_{\red})(p^{m}(K_X+\Delta))) =1$,   
for all $m \in \Z_{>0}$. 
\item 
$h^1(X, B_m\Omega_X^1(\log \Delta_{\red})(p^{m}(K_X+\Delta))) =1$   
for all $m \in \Z_{>0}$. 
\end{enumerate}
\end{proposition}


\begin{proof}
Set $D :=K_X+\Delta$. Let us show (1). 
Recall that $\MO_X(p(K_X+\Delta)) \simeq \MO_X(-1)$ 
(Proposition \ref{p-P^2-LDP}(2)). 
Therefore, we obtain 
$H^j(X, \Omega_X^1 \otimes \MO_X(p(K_X+\Delta))) =0$ 
for any $j \in \Z$ (Lemma \ref{l-P^2-H0Omega}(1)--(3)). 
By  the exact sequence 
\[
0 \to \Omega_X^1 \to \Omega_X^1(\log \Delta_{\red}) 
\to \bigoplus_{i=1}^4 \MO_{L_i} \to 0, 
\]
we obtain 
\begin{align*}
H^j(X, \Omega_X^1(\log \Delta_{\red})(p(K_X+\Delta)))
&\simeq 
\bigoplus_{i=1}^4 
H^j(L_i, \MO_X(p(K_X+\Delta))|_{L_i})\\
&\simeq 
\bigoplus_{i=1}^4 
H^j(\mathbb P^1, \MO_{\mathbb P^1}(-1))=0
\end{align*}
for any $j \in \Z$. Thus (1) holds. 

The assertion (2) holds by (1) and (\ref{exact1}). 
Moreover, it follows from Lemma \ref{l-P^2-vanish}(2) that (2) implies (3). 
Let us show (4). 
By (\ref{exact4}), we have the following exact sequence: 
\begin{multline*}
H^0( X,  Z_{m-1}\Omega_X^1(\log \Delta_{\red})(p^mD)) 
\to 
H^0(X,  B_1\Omega_X^2(\log \Delta_{\red})(pD)) \\
\longrightarrow 
H^1(X,  Z_m\Omega_X^1(\log \Delta_{\red})(p^mD)) \to 
H^1(X,  Z_{m-1}\Omega_X^1(\log \Delta_{\red})(p^mD)). 
\end{multline*}
It holds that 
\begin{align*}
H^0(X,  Z_{m-1}\Omega_X^1(\log \Delta_{\red})(p^mD))&=0 \textrm{ (Lemma \ref{l-P^2-H0Omega}(5)), and } \\ 
h^0(X,  B_1\Omega_X^2(\log \Delta_{\red})(pD)) &=1 \textrm{ (Lemma \ref{l-P^2-B_1}(2)).}
\end{align*}
By {(1) and} (3), we have $H^1(X, 
 Z_{m-1}\Omega_X^1(\log \Delta_{\red})(p^mD))=0$, 
and therefore we also get $h^1(X,  Z_m\Omega_X^1(\log \Delta_{\red})(p^mD))=1$. 
Thus (4) holds.

Let us show (5). 
We have the following commutative diagram in which each horizontal sequence is exact (\ref{exact2}): 
{\small\[
\begin{tikzcd}[column sep=10pt, nodes={inner sep=3pt}]
0 \arrow{r} & F_*^{m-1}B_1\Omega_X^1(\log \Delta_{\red})(p^mD) \arrow{d}{=} \arrow{r} &
B_m\Omega_X^1(\log \Delta_{\red})(p^mD) \arrow{d}{\alpha_m} \arrow{r}{C} & B_{m-1}\Omega_X^1(\log \Delta_{\red})(p^{m-1}D) \arrow{d}{\alpha_{m-1}} \arrow{r} & 0\\
0 \arrow{r} & F_*^{m-1}B_1\Omega_X^1(\log \Delta_{\red})(p^mD) \arrow{r} &
Z_m\Omega_X^1(\log \Delta_{\red})(p^mD) \arrow{r}{C} & Z_{m-1}\Omega_X^1(\log \Delta_{\red})(p^{m-1}D) \arrow{r} & 0,
\end{tikzcd}
\]}
\!\!where the vertical arrows are the natural inclusions. 
We then obtain the following commutative diagram in which each horizontal sequence is exact:  
\[
\begin{tikzcd}[column sep=10pt, nodes={inner sep=4pt}]
H^0(B_{m-1}\widetilde{\Omega}_X^1) \arrow{r} \arrow{d}{H^0(\alpha_{m-1})} & H^1(F_*^{m-1}B_1\overline{\Omega}_X^1) \arrow{d}{=} \arrow{r} & H^1(B_m\widetilde{\Omega}_X^1) \arrow{r} \arrow{d}{H^1(\alpha_{m})} & H^1(B_{m-1}\widetilde{\Omega}_X^1) \arrow{r} \arrow{d}{H^1(\alpha_{m-1})} & H^{{\cred 2}}(F_*^{m-1}B_1\overline{\Omega}_X^1) \arrow{d}{=} \\
H^0(Z_{m-1}\widetilde{\Omega}_X^1) \arrow{r} & H^1(F_*^{m-1}B_1\overline{\Omega}_X^1) \arrow{r} & H^1({Z_m\widetilde{\Omega}_X^1}) \arrow{r} & H^1(Z_{m-1}\widetilde{\Omega}_X^1) \arrow{r} & H^{{\cred 2}}(F_*^{m-1}B_1\overline{\Omega}_X^1)
\end{tikzcd}
\]
where 
\[
B_1\overline{\Omega}^1_X := B_1\Omega^1_X(\log \Delta_{\red})(p^mD),  
\]
\[
B_{\ell}\widetilde{\Omega}^1_X := B_{\ell}\Omega^1_X(\log \Delta_{\red})(p^{{\ell}}D), 
\quad\text{and} \quad 
Z_{\ell}\widetilde{\Omega}^1_X := Z_{\ell}\Omega^1_X(\log \Delta_{\red})(p^{{\ell}}D). 
\]
Since $H^0(X,  B_{m-1}\widetilde{\Omega}_X^1) = H^0(X,  Z_{m-1}\widetilde{\Omega}_X^1) =0$ (Lemma \ref{l-P^2-H0Omega}(5)), $H^0(\alpha_{m-1})$ is an isomorphism. 
By the 5-lemma and induction on $m$, it suffices to show that 
\[
H^1(\alpha_1) : 
H^1(X, B_1\Omega_X^1(\log \Delta_{\red})(pD)) \to 
H^1(X, Z_1\Omega_X^1(\log \Delta_{\red})(pD))
\]
is an isomorphism. 
To this end, it is enough to prove the following. 
\begin{enumerate}
\renewcommand{\labelenumi}{(\roman{enumi})}
    \item $h^1(X, Z_1\Omega_X^1(\log \Delta_{\red})(pD))=1$. 
    \item $h^1(X, B_1\Omega_X^1(\log \Delta_{\red})(pD))=1$. 
    \item $H^1(\alpha_1)$ is injective. 
\end{enumerate}
By (4), (i) holds. 

Let us show (ii). 
We have the following exact sequence (Lemma \ref{lem:Serre's map}):
\[
0 \to \MO_X(D) \to F_*\MO_X(pD) \to 
B_1\Omega_X^1(\log \Delta_{\red})(pD) \to 0, 
\]
which induces the following exact sequence: 
\[
H^1(X,  \MO_X(pD)) \to 
H^1(X,  B_1\Omega_X^1(\log \Delta_{\red})(pD))
\to H^2(X,  \MO_X(D)) \to H^2(X,  \MO_X(pD)). 
\]
{We have $h^2(X, \MO_X(D))=h^2(X, \MO_X(K_X+\Delta)) = h^2(X, \MO_X(K_X)) = h^0(X, \MO_X)=1$. 
This, together with $h^1(X, \MO_X(pD))=0$ and $h^2(X, \MO_X(pD))=0$ (Proposition \ref{p-P^2-LDP}(2)), implies that $h^1(X,  B_1\Omega_X^1(\log \Delta_{\red})(pD))=1$. 
}
Thus (ii) holds. 

Let us show (iii). 
We have the following exact sequence (\ref{exact1}): 
\[
0 \to 
B_1\Omega_X^1(\log \Delta_{\red})(pD))
\xrightarrow{\alpha_1} 
Z_1\Omega_X^1(\log \Delta_{\red})(pD)) 
\xrightarrow{C} 
\Omega_X^1(\log \Delta_{\red})(D) \to 0. 
\]
Hence (iii) follows from $H^0(X,  \Omega_X^1(\log \Delta_{\red})(D))=0$ (Lemma \ref{l-P^2-H0Omega}(4)). 
Thus (5) holds. 
\end{proof}

\begin{lemma}\label{l-Bn-Cartier}
Let $X$ be a smooth variety, 
$E$ a reduced simple normal crossing divisor, 
and 
$D$ a $\Q$-divisor satisfying $\Supp\,\{D\} \subseteq E$. 
Let $\zeta : B_{n+1}\Omega^1_X(\log E) \to B_{n}\Omega^1_X(\log E)$ be the homomorphism 
that completes the following commutative diagram in which each horizontal sequence is exact (cf.\ Lemma \ref{lem:Serre's map}): 
\[
\begin{tikzcd}
0 \arrow{r} & W_{m+1}\MO_X(D) \arrow{r}{F} \arrow{d}{R} & F_*W_{m+1}\MO_X(pD) \arrow{r}{s} \arrow{d}{R} & B_{m+1}\Omega_X^1(\log E)(p^{m+1}D) \arrow{d}{\zeta} \arrow{r} & 0 \\
0 \arrow{r} & W_{m}\MO_X(D) \arrow{r}{F} & F_*W_{m}\MO_X(pD) \arrow{r}{s} & B_{m}\Omega_X^1(\log E)(p^mD) \arrow{r} & 0
\end{tikzcd}
\]
Then the {equality} $\zeta =C$ holds for the Cartier operator 
$C: B_{m+1}\Omega^1_X(\log E)(p^{m+1}D) \to B_{m}\Omega^1_X(\log E)(p^mD)$.
\end{lemma}

\begin{proof}
For open subsets $V \subset U \subset X$, the restriction map  
\[
\Gamma(U, B_{m}\Omega^1_X(\log E)(p^mD)) \to 
\Gamma(V, B_{m}\Omega^1_X(\log E)(p^mD))
\]
is injective \cite[Lemma 5.10]{KTTWYY1}. 
Therefore, it is enough to show the equality $\zeta =C$ after replacing $X$ by an open subset $X \setminus (\Supp\,E \cup \Supp\,D)$. 
Hence we may assume that $D=E=0$. 
In this case, the required equality $\zeta =C$ 
follows from \cite[Ch.\ I, (3.11.6) in Proposition 3.11]{illusie_de_rham_witt}. 
Note that, for $n \in \{m, m+1\}$, the above horizontal arrow $s: F_*W_n\MO_X \to B_n\Omega_X^1$ 
coincides with the homomorphism $F^{n-1}d: W_n\Omega^0_X \to B_n\Omega_X^1$ appearing in 
\cite{illusie_de_rham_witt} as pointed out in 
\cite[Ch.\ I, Remarques 3.12(a)]{illusie_de_rham_witt}, although \cite{illusie_de_rham_witt} omits $F_*$ as a minor difference.
\end{proof}

\begin{theorem}\label{t-P^2-main}
We use Notation \ref{n-P^2-cex}. 
Then $(X=\mathbb P^2, \Delta)$ is a log del Pezzo pair with standard coefficients which is not quasi-$F$-split. 
\end{theorem}

\begin{proof}
By Notation \ref{n-P^2-cex} and Proposition \ref{p-P^2-LDP}, 
$(X=\mathbb P^2, \Delta)$ is a log del Pezzo pair with standard coefficients. 

It suffices to show that $(X, \Delta)$ is not quasi-$F$-split. 
Set $D :=K_X+\Delta$. 
It follows from Lemma \ref{l-Bn-Cartier} that 
we have the following commutative diagram in which each horizontal sequence is exact: 
\[
\begin{tikzcd}
  0 \arrow{r} & \MO_X(D) \arrow{r}{\Phi} \arrow{d}{=} & Q_{X, D, m+1} \arrow{r} \arrow{d} & B_{m+1}\Omega_X^1(\log \Delta_{\red})(p^{m+1}D) \arrow{r} \arrow{d}{C} & 0 \\
  0 \arrow{r} & \MO_X(D) \arrow{r}{\Phi} & Q_{X, D, m} \arrow{r} & B_{m}\Omega_X^1(\log \Delta_{\red})(p^{m}D) \arrow{r} & 0.
\end{tikzcd}
\]
We then obtain a commutative diagram 
\[
\begin{tikzcd}
H^1(X,   B_{m+1}\Omega_X^1(\log \Delta_{\red})(p^{m+1}D)) \arrow{r}{\delta_{m+1}} \arrow{d}{C} & H^2(X,  \MO_X(D)) \arrow{d}{=} \\
H^1(X,  B_{m}\Omega_X^1(\log \Delta_{\red})(p^mD)) \arrow{r}{\delta_m} & H^2(X,  \MO_X(D)).
\end{tikzcd}
\]
By (\ref{eq:coh-def-of-qfsplit}), to show that $(X,\Delta)$ is not quasi-$F$-split, 
 it suffices to prove that $\delta_m$ is nonzero for every $m>0$. Hence it is enough to prove (1) and (2) below. 
\begin{enumerate}
    \item $\delta_1 : H^1(X,  B_1\Omega_X^1(\log \Delta_{\red})(pD)) 
    \to H^2(X,  \MO_X(D))$ is nonzero. 
    \item $C: H^1(X, B_{m+1}\Omega_X^1(\log \Delta_{\red})(p^{m+1}D)) 
    \to H^1(X,  B_{m}\Omega_X^1(\log \Delta_{\red})(p^mD))$ is an isomorphism for every $m \geq 1$. 
\end{enumerate}

Let us show (1). 
By $Q_{X, D, 1} = F_*\MO_X(pD)$, we have the following exact sequence: 
\[
0=H^1(X,  Q_{X, D, 1}) \to H^1(X,  B_1\Omega_X^1(\log \Delta_{\red})(pD)) 
    \xrightarrow{\delta_1} H^2(X,  \MO_X(D)). 
\]
Hence $\delta_1$ is injective. 
By $H^1(X,  B_1\Omega_X^1(\log \Delta_{\red})(pD)) \neq 0$ 
(Proposition \ref{p-P^2-BZ-key}(5)), 
$\delta_1$ is nonzero. 
Thus (1) holds. 

Let us show (2). 
We have an exact sequence (\ref{exact2}): 
\begin{multline*}
0 \to F_*^mB_1\Omega_X^1(\log \Delta_{\red})(p^{m+1}D) \to 
B_{m+1}\Omega_X^1(\log \Delta_{\red})(p^{m+1}D) \\
\overset{C}{\longrightarrow} 
B_{m}\Omega_X^1(\log \Delta_{\red})(p^{m}D) \to 0. 
\end{multline*}
By $H^1(X, B_1\Omega_X^1(\log \Delta_{\red})(p^{m+1}D))=0$ (Lemma \ref{l-P^2-vanish}(1)),  we get an injection: 
\[
C: H^1(X,  B_{m+1}\Omega_X^1(\log \Delta_{\red})(p^{m+1}D)) 
\hookrightarrow H^1(X,  B_{m}\Omega_X^1(\log \Delta_{\red})(p^{m}D)).
\]
It follows from $h^1(X, B_{m+1}\Omega_X^1(\log \Delta_{\red})(p^{m+1}D)) 
=h^1(X, B_{m}\Omega_X^1(\log \Delta_{\red})(p^{m}D)) =1$ (Proposition \ref{p-P^2-BZ-key}(5)) that this map $C$ is an isomorphism. 
Thus (2) holds. 
\end{proof}


\begin{remark}\label{r-small-counter}
\begin{enumerate}
\item 
Take $p \in \{2, 3, 5\}$ and let $k$ be an algebraically closed field of characteristic $p$. 
Then there exists a log del Pezzo pair $(X, 0)$ over $k$ which is not quasi-$F$-split. 
Indeed, $X$ is quasi-$F$-split if and only if $X$ is log liftable by \cite[Theorem 6.3]{KTTWYY1}, and there exists a klt del Pezzo surface over $k$ which is not log liftable by \cite[Example 6.1]{Kaw3}. 
For the construction of the examples, we refer to \cite[Theorem 4.2(6)]{CT19-2} ($p=2$), \cite[Theorem 1.1]{Ber} ($p=3$), and \cite[Proposition 5.2]{ABL20} ($p=5$). 
\item 
If $p \in \{2, 3, 5, 7, 11, 17, 19, 23, 41\}$ and $k$ is an algebraically closed field of characteristic $p$, 
then there exists, by (1) and Theorem \ref{t-P^2-main}, a log del Pezzo pair $(X, \Delta)$ with standard coefficients which is not quasi-$F$-split. 
On the other hand, the authors do not know whether such examples exist in the other characteristics $p\leq 41$, that is, $p \in \{  13, 29, 31, 37\}$. 
\end{enumerate}
\end{remark}

%% file: section7.tex
\section{Cone correspondence for quasi-F-splitting}\label{s-Qcone-QFS}

Let $(X, \Delta)$ be a projective log pair such that $-(K_X+\Delta)$ is ample and 
$\Delta$ has standard coefficients with $\rdown{\Delta}=0$. 
The purpose of this section is to prove that $(X, \Delta)$ is quasi-$F$-split if and only if its cone $R_{\m}$ is quasi-$F$-split (Theorem \ref{cor:corresponding}), 
where 
\[
R := \bigoplus_{d \geq 0} H^0(X, \MO_X(-d(K_X+\Delta))) 
\quad\text{and}
\quad 
\m := \bigoplus_{d > 0} H^0(X, \MO_X(-d(K_X+\Delta))) \subseteq R. 
\]
The strategy is to generalise an analogue of this result for $F$-splittings which was established by \cite{Watanabe91}. 
Most of this section is devoted to introducing 
functorial graded structures on several key modules that feature in the theory of quasi-$F$-splittings. 

\subsection{Notation on graded modules}

\begin{enumerate}
\item 
    Let $Z$ be a submonoid of $\Q$, that is, 
    $Z$ is a subset of $\Q$ such that $0 \in Z$ and $Z$ is closed under the addition (e.g. $\Z_{\geq 0}, \Z$). 
    For a $Z$-graded ring $R = \bigoplus_{d \in Z} R_d$, 
    we set $R_{d'} :=0$ for any $d' \in \Q \setminus Z$ and  
    we consider $R$ as a $\Q$-graded ring via $R = \bigoplus_{d \in \Q} R_d$. 
    Similarly, $Z$-graded $R$-modules can be considered as $\Q$-graded $R$-modules. 
    \item 
    Given $\Q$-graded rings $R = \bigoplus_{d \in \Q} R_d$ and $S = \bigoplus_{d \in \Q} S_d$, 
    we say that $\varphi : R \to S$ is a {\em graded ring homomorphism}  
    if $\varphi$ is a ring homomorphism such that $\varphi(R_d) \subseteq S_d$ 
    for every $d \in \Q$. 
    \item 
    Given a $\Q$-graded ring $R = \bigoplus_{d \in \Q} R_d$ and 
    $\Q$-graded $R$-modules $M =  \bigoplus_{d \in \Q} M_d$ and $N =  \bigoplus_{d \in \Q} N_d$, 
    we say that $\psi : M \to N$ is a {\em graded $R$-module  homomorphism} 
    if $\psi$ is an $R$-module homomorphism such that 
    $\psi(M_d) \subseteq N_d$ for any $d \in \Q$. 
   Note that we assume that $\psi$ preserves degrees, although in the literature it is sometimes allowed for the
   degrees to be shifted. 
\end{enumerate}

\subsection{Witt rings of graded rings}

\begin{proposition}\label{p-Witt-of-graded}
Let $R = \bigoplus_{d \in \Q} R_d$ be a $\Q$-graded $\F_p$-algebra. 
Then 
\begin{enumerate}
    \item 
    $W_n(R)$ has the $\Q$-graded ring structure 
\[
W_n(R) = \bigoplus_{e \in \Q} W_n(R)_e
\]
such that, for every $e \in \Q$:
\begin{align*}
W_n(R)_e = \{ (a_0, &a_1, \cdots, a_{n-1}) \in W_n(R) \,|\\
&\, a_0 \in R_e, a_1 \in R_{pe}, a_2\in R_{p^2e}, \cdots , a_{n-1} \in R_{p^{n-1}e}\}.
\end{align*}
In particular, $W_n(R)_0 = W_n(R_0)$ holds. 
\end{enumerate}
Note that, by (1), 
$W_n(R)_e$ is a $W_n(R_0)$-module for every $e \in \Q$. 
For all $r \in \Z_{\geq 0}$ and $e \in \Q$, we set 
\[
(F^r_*W_n(R))_e := F_*^r(W_n(R)_{p^re}), 
\]
which is a $W_n(R_0)$-module. Then the following hold. 
\begin{enumerate}
\setcounter{enumi}{1}
\item 
For every $r \in \Z_{\geq 0}$, 
$F_*^rW_n(R)$ has the $\Q$-graded $W_n(R)$-module structure 
\[
F_*^rW_n(R)  = \bigoplus_{e \in \Q} (F^r_*W_n(R))_e.
\]
\item 
The Frobenius ring homomorphism of $W_n(R)$ 
\[
F : W_n(R) = \bigoplus_{e \in \Q} W_n(R)_e \to F_*W_n(R) = \bigoplus_{e \in \Q} (F_*W_n(R))_{e}
\]
is a $\Q$-graded ring homomorphism. 
\item 
The Verschiebung homomorphism of $W_n(R)$ 
\[
V : F_*W_n(R) = \bigoplus_{e \in \Q} (F_*W_n(R))_{e} \to 
W_{n+1}(R) = \bigoplus_{e \in \Q} W_{n+1}(R)_{e}
\]
is a $\Q$-graded $W_n(R)$-module homomorphism. 
\item 
For any $n, m \in \Z_{> 0}$, 
the ring homomorphism 
\[
R : W_{n+m}(R) = \bigoplus_{e \in \Q} W_{n+m}(R)_e 
\to W_n(R) = \bigoplus_{e \in \Q} W_n(R)_e 
\]
is a $\Q$-graded ring homomorphism. 
\end{enumerate}
\end{proposition}

\begin{proof}
We omit a proof, as the same argument as in Proposition \ref{p-Witt-of-graded-can-div} works by setting $q=0$.
\qedhere
\end{proof}

\begin{remark}\label{r-Witt-of-graded}
In what follows, we are mainly interested in 
$\Z_{\geq 0}$-graded $\F_p$-algebras $R = \bigoplus_{d \in \Z_{\geq 0}} R_d$. 
Then $R$ can be considered as a $\Q$-graded $\F_p$-algebra $R= \bigoplus_{d \in \Q} R_d$, where we set $R_{d'} :=0$ for all $d' \in \Q \setminus \Z_{\geq 0}$. 
In this case, $W_n(R)$ can be considered as a $p^{-(n-1)}\Z_{\geq 0}$-graded ring and 
$F_*W_n(R)$ is a $p^{-n}\Z_{\geq 0}$-graded $W_n(R)$-module. 
\end{remark}

\begin{nothing}[Graded structure on $R(qK_R)$]\label{graded-can-div}
Let $k$ be a field. 
Let $R = \bigoplus_{d \in \Z} R_d$ be a $\Z$-graded 
normal integral domain such that $R_d=0$ for every $d < 0$, $R_0=k$, and $R$ is a finitely generated $k$-algebra. 
Fix a canonical divisor $K_R$ on $R$. 
We equip $R(K_R)$ with the $\Z$-graded $R$-module structure as in \cite[Definition (2.1.2)]{Got-Wat78}. 
We then define the $\Z$-graded $R$-module structure on {$R(qK_R)$} by taking the reflexive hull of $R(K_R)^{\otimes q}$ for every $q \in \Z_{\geq 0}$ 
(recall that given $\Z$-graded $R$-modules $M = \bigoplus_{d \in \Z} M_d$ and $N = \bigoplus_{d \in \Z} N_d$, 
the tensor product is a $\Z$-graded $R$-module given by 
$(M \otimes_R N)_d =\bigoplus_{d = d_1 +d_2}\mathrm{Im}( M_{d_1} \otimes_{R_0} N_{d_2} \to M \otimes_R N)$).  
\end{nothing}

\begin{proposition}\label{graded-can-div2}
We use the notation as in (\ref{graded-can-div}). 
Then the natural map
\[
\beta \colon R(q_1K_R) \otimes_R R(q_2K_R) \to R((q_1+q_2)K_R), 
\]
induced by the product $\mathrm{Frac}\,R \otimes_R \mathrm{Frac}\,R \to \mathrm{Frac}\,R, f \otimes g \mapsto fg$, 
is a $\Z$-graded $R$-module homomorphism for all $q_1, q_2 \in \Z_{>0}$.
\end{proposition}

\begin{proof}
We have the following commutative diagram of $R$-module homomorphisms: 
\[
\begin{tikzcd}
R(K_R)^{\otimes q_1} \otimes_R R(K_R)^{\otimes q_2} 
\arrow[r, "\alpha", "\text{graded}"'] \arrow[d, "\zeta", "\text{graded}"'] 
\arrow[rr, "\gamma := \beta \circ \alpha", bend left=16, "\text{graded}"'] 
&R(q_1K_R) \otimes_R R(q_2K_R) \arrow[r, "\beta"] \arrow[d, "\xi", "\text{graded}"']
& R((q_1 + q_2)K_R)\arrow[d, "\eta", "\text{graded, bij.}"']
\\
(R(K_R)^{\otimes q_1} \otimes_R R(K_R)^{\otimes q_2})^{**} 
\arrow[r, "\alpha^{**}", "\text{graded, bij.}"'] \arrow[rr, "\gamma^{**}", bend right=16, "\text{graded, bij.}"']
& (R(q_1K_R) \otimes_R R(q_2K_R))^{**}  \arrow[r, "\beta^{**}", "\text{bij.}"'] 
& R((q_1 + q_2)K_R))^{**},
\end{tikzcd}
\]
where the lower horizontal arrows and all the vertical arrows are obtained by taking double duals. 
It is clear that $\eta, \alpha^{**}, \beta^{**}, \gamma^{**}$ are bijective. 
All the vertical arrows are graded homomorphisms. 
Since $\alpha$ and $\gamma$ are graded homomorphisms by definition, 
so are $\alpha^{**}$ and $\gamma^{**}$. 
Then $\beta^{**}$ is a graded homomorphism. 
The composition $\beta^{**} \circ \xi = \eta \circ \beta$ is a graded homomorphism, 
and hence $\beta = \eta^{-1} \circ ( \eta \circ \beta)$ is a graded homomorphism. 
\end{proof}

\begin{proposition}\label{p-Witt-of-graded-can-div}
Let $k$ be an $F$-finite field of characteristic $p>0$.
Let $R = \bigoplus_{d \in \Q} R_d$ be a $\Q$-graded normal integral domain such that  $R_d=0$ for $d \in \Q \setminus \Z_{\geq 0}$, $R_0=k$, and $R$ is a finitely generated $k$-algebra. 
Fix $q \in \Z_{>0}$ and a canonical divisor $K_R$ on $R$. 
We define the graded structure on $R(qK_R)$ for $q \in \Z_{>0}$ 
as in (\ref{graded-can-div}). 
Then the following hold. 
\begin{enumerate}
    \item 
    $W_n(R)(qK_R)$ has the $\Q$-graded $W_n(R)$-module structure 
\[
W_n(R)(qK_R) = \bigoplus_{e \in \Q} W_n(R)(qK_R)_e
\]
such that, for every $e \in \Q$: 
\begin{align*}
\quad W_n(R)(qK_R)_e = \{ (a_0, a_1, &\cdots, a_{n-1}) \in W_n(R)(qK_R) \,|\,
\\
&a_i \in R(p^iqK_R)_{p^ie}\text{ for every }1 \leq i \leq n-1 \}.
\end{align*}
 
\end{enumerate}
For all $r \in \Z_{\geq 0}$ and $e \in \Q$, we set 
\[
(F^r_*W_n(R)(qK_R))_e := F_*^r(W_n(R)(qK_R)_{p^re}). 
\]
Then the following hold. 
\begin{enumerate}
\setcounter{enumi}{1}
\item 
For every $r \in \Z_{\geq 0}$, 
$F_*^rW_n(R)(qK_R)$ has the $\Q$-graded $W_n(R)$-module structure 
\[
F_*^rW_n(R)(qK_R)  = \bigoplus_{e \in \Q} (F^r_*W_n(R)(qK_R))_e.
\]
\item 
The induced map 
\[
F : W_n(R)(qK_R)  \to F_*W_n(R)(pqK_R) 
\]
by the Frobenius homomorphism 
is a $\Q$-graded $W_n(R)$-module homomorphism. 
\item 
The induced map 
\[
V : F_*W_n(R)(pqK_R)  \to 
W_{n+1}(R)(qK_R) 
\]
by the Verschiebung homomorphism 
is a $\Q$-graded $W_n(R)$-module homomorphism. 
\item 
For any $n, m \in \Z_{> 0}$, 
the restriction map
\[
R: W_{n+m}(R)(qK_R) 
\to W_n(R)(qK_R) 
\]
is a $\Q$-graded $W_n(R)$-module homomorphism. 
\end{enumerate}

\end{proposition}

\begin{proof}
Let us show (1). 
It suffices to show the following assertions (i)--(iii). 
\begin{enumerate}
\renewcommand{\labelenumi}{(\roman{enumi})}
    \item $W_n(R)(qK_R)_e$ is a subgroup of $W_n(R)(qK_R)$. 
    \item $W_n(R)(qK_R)_{e} \cdot W_n(R)_{e'} \subseteq W_n(R)(qK_R)_{e+e'}$ holds for all $e, e' \in \Q$. 
    \item $W_n(R)(qK_R) = \bigoplus_{e \in \Q} W_n(R)(qK_R)_e$. 
\end{enumerate}
Let us show (i).
We take $\varphi=(\varphi_0,\ldots,\varphi_{n-1}) \in W_n(R)(qK_R)_e$ and $\psi=(\psi_1,\ldots,\psi_{n-1}) \in W_n(R)(qK_R)_e$.
By \cite[Lemma 2.3]{tanaka22}, we have
\[
\varphi+\psi=(S_0(\varphi_0,\psi_0),S_1(\varphi_0,\psi_0,\varphi_1,\psi_1),\ldots)
\]
for some polynomials
\[
S_m(x_0,y_0,\ldots,x_m,y_m) \in \Z[x_0,y_0,\ldots,x_m,y_m].
\]
By \cite[Lemma 2.4]{tanaka22}, if we pick a monomial $x_0^{a_0}y_0^{b_0} \cdots x_m^{a_m}y_m^{b_m}$ appearing in the monomial decomposition of $S_m$, then we have
\[
\sum_{i=0}^m p^i(a_i+b_i)=p^m.
\]
By $\varphi_i,\psi_i \in R(p^iqK_R)_{p^ie}$ and Proposition \ref{graded-can-div2},
it holds that 
\[
\varphi_0^{a_0}\psi_0^{b_0} \cdots \varphi_m^{a_m}\psi_m^{b_m} \in R(p^mK_R)_{p^me}.
\]
Therefore, we have $\varphi+\psi \in R(qK_R)_{e}$.

Let us show (ii). 
Fix $e, e' \in \Q$. 
Take two elements $V^m\underline{b} \in W_n(R)(qK_R)_e$ and $V^{m'}\underline{b'} \in W_n(R)_{e'}$, 
where $0 \leq m \leq n-1, 0 \leq m' \leq n-1, b \in R(p^mqK_R)_{p^me}, b' \in R_{p^{m'}e'}$. 
{\cred 
It is enough to prove 
$(V^m\underline{b}) \cdot (V^{m'}\underline{b'})   \in W_n(R)(qK_R)_{e+e'}$, 
because $\beta \in  W_n(R)(qK_R)_e$ (resp.\ $\beta'  \in  W_n(R)_{e'}$) 
can be written by $\beta = \sum_m V^m \underline{b_m}$ (resp.\ $\beta'  = \sum_{m'}V^m\underline{b'_{m'}}$), 
where $b_m \in R(p^mqK_R)_{p^me}$ and $b'_{m'} \in R_{p^{m'}e'}$.} 
We 
obtain 
\[
(V^m\underline{b}) \cdot (V^{m'}\underline{b'}) 
= V^{m+m'}( ( F^{m'}\underline{b})\cdot ( F^{m}\underline{b'})) 
= V^{m+m'}(\underline{ b^{p^{m'}} \cdot b'^{p^m} }).
\]
We have $b^{p^{m'}} \in R(p^{m+m'}qK_R)_{p^{m+m'}e}$ and $b'^{p^m} \in R_{p^{m+m'}e'}$, which implies 
$b^{p^{m'}} \cdot b'^{p^m} \in R(p^{m+m'}qK_R)_{p^{m+m'}(e+e')}$ (Proposition \ref{graded-can-div2}). 
Therefore, we get 
$V^{m+m'}(\underline{ b^{p^{m'}} \cdot b'^{p^m} })  \in W_n(R)(qK_R)_{e+e'}$. 
Thus (ii) holds. 

Let us show (iii). 
It is easy to check that the equality $W_n(R)(qK_R) = \sum_{e \in \Q} W_n(R)(qK_R)_e$ holds. 
It suffices to show that this sum $\sum_{e \in \Q} W_n(R)(qK_R)_e$ is a direct sum. 
Assume that 
\[
\zeta_1 + \cdots + \zeta_r =0
\]
holds for $\zeta_1 \in W_n(R)(qK_R)_{e_1}, ..., \zeta_r \in W_n(R)(qK_R)_{e_r}$ 
with $e_1 < \cdots <e_r$. 
It is enough to prove that $\zeta_1 = \cdots = \zeta_r=0$. 
We have 
\[
\zeta_i =(\zeta_{i, 0}, \zeta_{i, 1}, ..., \zeta_{i,n-1}) \in W_n(R)(qK_R)
\]
for some $\zeta_{i, j} \in R(p^jqK_R)$. 
By 
\[
(0, 0, ...) = 0 = \sum_{i=1}^r \zeta_i = 
\left( \sum_{i=1}^r \zeta_{i, 0}, ...\right),  
\]
we get $\sum_{i=1}^r \zeta_{i, 0}=0$, which implies 
$\zeta_{1, 0} = \cdots =\zeta_{r, 0}=0$. 
It follows from 
\[
(0, 0, ...) = 0 = \sum_{i=1}^r \zeta_i = 
\left(0, \sum_{i=1}^r \zeta_{i, 1}, ...\right)  
\]
that $\sum_{i=1}^r, \zeta_{i, 1}=0$, which implies 
$\zeta_{1, 1} = \cdots = \zeta_{r, 1}=0$. 
Repeating this procedure, we obtain $\zeta_1 = \cdots = \zeta_r =0$. 
Thus (iii) holds. 
This completes the proof of (1).

We now show the following $(*)$. 
\begin{enumerate}
    \item[$(*)$] For $r \in \Z_{\geq 0}$ and the Frobenius action
    $F^r : W_n(R)(qK_R) \to W_n(R)(p^rqK_R)$, it holds that $F^r ( W_n(R)(qK_R)_e) \subseteq  W_n(R)(p^rqK_R)_{p^re}$. 
\end{enumerate}
To this end, we may assume that $r=1$. 
Fix $e \in \Q$ and $\alpha \in W_n(R)(qK_R)_e$. 
We have 
\[
\alpha=  (a_0, a_1, ..., a_{n-1}) \qquad \text{for} \quad \text{some} \quad a_i \in R(p^iqK_R)_{p^i e}. 
\]
It holds that 
\[
F(\alpha) = F(a_0, a_1, ..., a_{n-1}) = (a_0^p, a_1^p, ..., a_{n-1}^p). 
\]
Then $a_i \in R(p^iqK_R)_{p^i e}$ implies $a_i^p \in R(p^{i+1}qK_R)_{p^{i+1} e}$ by Proposition \ref{graded-can-div2}, that is,  $F(\alpha) \in W_n(R)(pqK_R)_{pe}$. 
Thus $(*)$ holds.

Let us show (2). 
It follows from (1) that we have the direct sum decomposition 
$F_*^rW_n(R)(qK_R) = \bigoplus_{e \in \Q} (F_*^rW_n(R)(qK_R))_e$  
as an additive group. 
Hence it suffices to show that $W_n(R)_e \cdot (F_*^rW_n(R)(qK_R))_{e'} 
\subseteq (F_*^rW_n(R)(qK_R))_{e+e'}$. 
Take $\zeta \in W_n(R)_e$ and $F^r_* \xi \in (F_*^rW_n(R)(qK_R))_{e'}$, 
where $F^r_* \xi$ denotes the same element as $\xi \in W_n(R)(qp^rK_R)_{p^re'}$ 
via the set-theoretic equality $(F_*^rW_n(R)(qK_R))_{e'} = W_n(R)(qp^rK_R)_{p^re'}$. 
We have $F^r(\zeta) \in W(R)_{p^re}$. 
Hence we obtain 
\[
\zeta \cdot (F^r_* \xi) = F^r_*( F^r(\zeta) \cdot \xi) 
\in F_*^r( W(R)_{p^re} \cdot  W_n(R)(qp^rK_R)_{p^re'}) 
\]
\[
\subseteq F_*^r(W_n(R)(qp^rK_R)_{p^r(e+e')}) =(F^r_* W_n(R)(qK_R))_{e+e'}. 
\]
Thus (2) holds.

The assertions (3)--(5) follow from $(*)$. 
\end{proof}

\subsection{Cones and quasi-F-splittings}
We start by making a general comment about graded $W_n(R)$-modules.

\begin{nothing}[Graded structures on cohomologies]\label{n-graded-str}
Let $R = \bigoplus_{d \in \Z_{\geq 0}}R_d$ be an $F$-finite noetherian 
$\Z_{\geq 0}$-graded $\F_p$-algebra such that $R_0$ is a field. 
Set $\m := \bigoplus_{d \in \Z_{> 0}}R_d$  and 
$U :=\Spec\,R \setminus \{\m\}$.
Fix $n \in \Z_{>0}$ and $\nu \in \Z_{\geq 0}$. 
Let $M$ be a 
$p^{-\nu}\Z$-graded $W_n(R)$-module. 
Fix $i \in \Z_{\geq 0}$. 
In what follows, we introduce 
\begin{enumerate}
\item a $p^{-\nu}\Z$-graded $W_n(R)$-module structure on $H^i(U, \widetilde{M}|_U)$, and 
\item a $p^{-\nu}\Z$-graded $W_n(R)$-module structure on $H^{i+1}_{\m}(M)$. 
\end{enumerate}
Moreover, we show that 
\begin{enumerate}
\item[(3)] the connecting homomorphism  $H^i(U, \widetilde{M}|_U) \to H^{i+1}_{\m}(M)$ 
is a graded $W_n(R)$-module homomorphism. 
\end{enumerate}

For every nonzero homogeneous element $ f \in W_n(R)_{>0}$, we have 
\[
\Gamma(D(f), \widetilde{M}|_{U}) = \Gamma(D(f), \widetilde{M}) = M_f, 
\]
which is a $p^{-\nu}\Z$-graded $W_n(R)_f$-module. 
Fix  nonzero homogeneous elements $f_1,\ldots,f_r \in W_n(R)$ such that $(f_1, ..., f_r) = W_n(R)_{>0}$, 
that is, $D(f_1)\cup \cdots \cup D(f_r)=U$. 
Then we have the \v{C}ech complex of $p^{-\nu}\Z$-graded ${\cred W_n(R)}$-modules 
\[
C(f_1,\ldots,f_r):=\left( \bigoplus_{1 \leq i \leq r} M_{f_i} \to \bigoplus_{1 \leq i < j \leq r} M_{f_if_j} \to \cdots \to M_{f_1\cdots f_r}\right).
\]
Since the $i$-th cohomology of this complex is $H^i(U,\widetilde{M}|_U)$ 
\cite[Theorem III.4.5]{hartshorne77},  
the complex $C(f_1,\ldots,f_r)$ gives a $p^{-\nu}\Z$-graded $W_n(R)$-module 
structure on $H^i(U,\widetilde{M}_U)$. 
Note that this does not depend on the choice of $f_1,\ldots,f_r \in W_n(R)$. 
Indeed, if we pick another nonzero homogeneous element $f_{r+1} \in W_n(R)_{>0}$, 
then  $C(f_1,\ldots,f_r)\to C(f_1,\ldots,f_{r+1})$ is a homomorphism of complexes of $p^{-\nu}\Z$-graded $W_n(R)$-modules. 

Similarly, we equip $H^{i+1}_{\m}(M)$ with a $p^{-\nu}\Z$-graded $W_n(R)$-module structure 
by using the fact that $H^{i+1}_{\m}(M)$ is the $i$-th cohomology of $M \to C(f_1, ..., f_r)$. 
Then  the connecting homomorphism $H^i(U,\widetilde{M}|_U) \to H^{i+1}_\m(M)$ 
is a graded $W_n(R)$-module homomorphism, because it is induced by the following commutative diagram: 
\[
\begin{tikzcd}
    0 \arrow{r}\arrow{d} & C(f_1,\ldots,f_r) \arrow{d}\\
    M \arrow{r} & C(f_1,\ldots,f_r).
\end{tikzcd}
\]
\end{nothing}

\begin{notation}\label{n-cone-QFS}
Let $k$ be an $F$-finite field 
of characteristic $p>0$. 
Let $X$ be a projective normal variety over $k$
with $\dim X \geq 1$  and $H^0(X, \MO_X)=k$. 
Let $D$ be an ample $\Q$-Cartier $\Q$-divisor. 
Set 
\[
R := R(X, D) := \bigoplus_{d \in \Z_{\geq 0}} H^0(X, \MO_X(dD))t^d \subseteq {K(X)}[t],  
\]
which is a $\Z_{\geq 0}$-graded subring of the standard $\Z_{\geq 0}$-graded polynomial ring ${K(X)}[t]$. 
Note that $R$ is a finitely generated $\Z_{\geq 0}$-graded  $k$-algebra 
and $\Spec R$ is an affine normal variety (Theorem \ref{t-Qcone-birat}). 
Let $D = \sum_{i=1}^r \frac{\ell_i}{d_i}D_i$ be the irreducible decomposition, 
where $\ell_i$ and $d_i$ are coprime integers satisfying $d_i>0$ for each $1 \leq i \leq r$. 
Set $D':=\sum_{i=1}^r \frac{d_i-1}{d_i}D_i$.  

For the graded maximal ideal $\m := \bigoplus_{d >0} H^0(X, \MO_X(dD))t^d \subseteq R$, 
we set $U := \Spec R \setminus \{\m\}$. 
{\cred 
For the quasi-coherent finitely generated $\MO_X$-algebras
\[
\mathcal A_{X, D} := \bigoplus_{d \geq 0} \MO_X(dD)t^d \qquad 
{\rm and} \qquad 
\mathcal A^{\circ}_{X, D} := \bigoplus_{d \in \Z} \MO_X(dD)t^d, 
\]
we set 
\[
W_{X, D} := \Spec_X \mathcal A_{X, D} \qquad 
{\rm and} \qquad  W^{\circ}_{X, D} := \Spec_X \mathcal A^{\circ}_{X, D}.
\]
For the natural  morphism $\mu \colon W_{X, D} \to \Spec R$ (cf.\ Definition \ref{d-Qcone2}),} 
we have $\mu|_{W^{\circ}_{X, D}} : W^{\circ}_{X, D} \xrightarrow{\simeq} U$ (Theorem \ref{t-Qcone-birat}), 
and hence we have the induced affine morphism: 
\[
\rho : U \xrightarrow{(\mu|_{W^{\circ}_{X, D}})^{-1}, \simeq} W^{\circ}_{X, D} \xrightarrow{\pi^{\circ}} X, 
\]
{\cred where $\pi^{\circ} \colon  W^{\circ}_{X, D} \to X$ denotes the natural projection.} 
For $n \in \Z_{>0}$, let 
\[
W_n(R) = \bigoplus_{e \in p^{-(n-1)}\Z_{\geq 0}} W_n(R)_e
\]
be  the $p^{-(n-1)}\Z_{\geq 0}$-graded ring structure induced by 
Proposition \ref{p-Witt-of-graded}. 
Set $W_n(R)_{>0} := \bigoplus_{e \in p^{-(n-1)}\Z_{>0}} W_n(R)_{{e}}$, 
which is a graded primary ideal such that $\sqrt{W_n(R)_{>0}}$ is a  maximal ideal. 
\end{notation}

\begin{proposition}\label{prop: affine morphism}
We use Notation \ref{n-cone-QFS}. 
Then the equality 
\[
\rho_*\cO_U =\bigoplus_{d \in \Z} \cO_X(d{D})t^d
\]
induces the following isomorphism of graded $R$-modules
\[
H^i(U,\cO_U) \simeq \bigoplus_{d \in \Z} H^i(X,\cO_X(d{D}))t^d, 
\]
where the graded structure of $H^i(U,\cO_{U})$ is defined in Notation \ref{n-graded-str}.
\end{proposition}

\begin{proof}
Note that $\rho_*\cO_U =\bigoplus_{d \in \Z} \cO_X(d{D})t^d$ 
{\cred holds} 
inside the constant sheaf ${K(X)[t,t^{-1}]}$.
By Theorem \ref{t-Qcone-birat}, we have
\[
W_{X,D}^{\circ} \simeq W_{X,D} \backslash \Gamma_{X,D} \simeq V_{X,D} \backslash \{v_{X,D}\} =U.
\]
By the definition of $W^{\circ}_{X,D}$, we have
\[
\rho_*\cO_U =\bigoplus_{d \in \Z} \cO_X(d{D})t^d
\]
holds inside the constant sheaf ${K(X)[t,t^{-1}]}$.
In what follows, for $d \in \Z_{>0}$ and a nonzero element $f \in H^0(X, \MO_X(dD))$,  set 
\[
\overline{f} := ft^{d} \in H^0(X, \MO_X(dD))t^{d} \subseteq 
\bigoplus_{d=0}^{\infty} H^0(X, \MO_X(dD))t^d =R(X, D). 
\]
Since we have $\Proj\, R(X,D) \simeq X$ by \cite[Proposition 3.3 (a)]{Demazure}, we can define an open subset $D_{+}(\overline{f})$ as in \cite[Proposition II.2.5]{hartshorne77}.
Let $f_1,\ldots , f_r \in R(X,D)$ be homogeneous elements of positive degree such that $\{D(f_i) \mid 1 \leq i \leq r\}$ is an open covering of $U$, then $\{D_+(\overline{f_i}) \mid 1 \leq i \leq r\}$ is an open covering of $X$.
By $\rho_*\cO_U=\oplus \cO_X(dD)t^d$ and $\rho^{-1}(D_+(\overline{f_i}))=D(f_i)$, 
we have the isomorphism of complexes of $\Z$-graded $R$-modules: 
\[
\begin{tikzcd}
 \bigoplus_i R_{f_i} \arrow[r] \arrow[d, "\simeq"] & \bigoplus_{i<j} R_{f_if_j} \arrow[r] \arrow[d, "\simeq"] & \cdots \\
\bigoplus_{d \in \Z} \bigoplus_i H^0(D_+(\overline{f_i}), \cO_X(dD)) \arrow[r] & \bigoplus_{d \in \Z} \bigoplus_{i<j}  H^0(D_+(\overline{f_if_j}), \cO_X(dD)) \arrow[r] & \cdots.   
\end{tikzcd}
\]
Thus, taking cohomologies, we have the following isomorphism of graded $R$-modules: 
\[
H^i(U,\cO_U) \simeq \bigoplus_{d \in \Z} H^i(X, \cO_X(dD))t^d.
\]
This completes the proof of Proposition \ref{prop: affine morphism}. 
\end{proof}

\begin{nothing}[Canonical divisors]
We use Notation \ref{n-cone-QFS}. 
Fix a canonical divisor $K_X$, that is, a Weil divisor $K_X$ such that $\MO_X(K_X) \simeq \omega_X$. 
By \cite[3.2]{Watanabe91} (for the convenience of the reader, we also provide a proof below), 
there exists a canonical divisor 
$K_R (=K_{\Spec R})$ on $\Spec R$ such that 
the following equality holds as a $\Z$-graded $R$-module for every $q \in \Z_{>0}$: 
\[
R(qK_R) 
= \bigoplus_{d \in \Z} H^0(X, \cO_X(q(K_X+D')+dD))t^d 
\subseteq K(X)[t, t^{-1}],  
\]
where the graded structure on $R(qK_R)$ is defined as in (\ref{graded-can-div}). 
We note that $H^i(U, \MO_U(qK_U))$ is a $\Z$-graded 
$R$-module for any $i \geq 0$ (\ref{n-graded-str}). 
\end{nothing}

\begin{proposition}
We use Notation \ref{n-cone-QFS}.
Fix a canonical divisor $K_X$, that is, a Weil divisor $K_X$ such that $\MO_X(K_X) \simeq \omega_X$. 
Then there exists a canonical divisor 
$K_R (=K_{\Spec R})$ on $\Spec R$ such that 
the following equality holds as a $\Z$-graded $R$-module for every $q \in \Z_{>0}$: 
\begin{equation}\label{eq: canonical module}
    R(qK_R) 
= \bigoplus_{d \in \Z} H^0(X, \cO_X(q(K_X+D')+dD))t^d 
\subseteq K(X)[t, t^{-1}],
\end{equation}
where the $\Z$-graded $R$-module structure on $R(qK_R)$ is defined as in (\ref{graded-can-div}).
\end{proposition}
\begin{proof}
By \cite[Theorem 2.8]{wat81}, we have an isomorphism 
\[
\omega_R \simeq \bigoplus_{d \in \Z} H^0(X, \underbrace{\cO_X(K_X+D'+dD)}_{=\cO_X({\cred \lceil} K_X+dD {\cred \rceil})})t^d.
\]
Therefore, there exists a canonical divisor $K_R$ on $\Spec R$ such that
\[
R(K_R) 
= \bigoplus_{d \in \Z} H^0(X, \cO_X(K_X+D'+dD))t^d 
\subseteq K(X)[t, t^{-1}].
\]
Thus, we obtain the case of $q=1$.
Next, we prove the equality (\ref{eq: canonical module}),
The right hand side of (\ref{eq: canonical module}) is denoted by $M_q$, which is a finite $\Z$-graded $R$-module.
\begin{claim}\label{claim:S_2}
We have $H^i_\m(M_q)=0$ for $i=0,1$.
\end{claim}
\begin{proof}[Proof of Claim \ref{claim:S_2}]
By a similar argument to the proof of Proposition \ref{prop: affine morphism},
\[
(\widetilde{M}_q)|_{U} \simeq \bigoplus_{d \in \Z} \cO_X(q(K_X+D')+dD)t^d
\]
as $\cO_U$-modules.
We consider the exact sequence
\[
0 \to H^0_\m(M_q) \to M_q \to H^0(U,(\widetilde{M}_q)|_{U}) \to H^1_{\m}(M_q) \to 0.
\]
Since
\[
H^0(U,(\widetilde{M}_q)|_{U}) \simeq \bigoplus_{d \in \Z} H^0(X,\cO_X(q(K_X+D')+dD))t^d=M_q,
\]
we get $H^i_\m(M_q)=0$ for $i=0,1$.
\end{proof}
By Claim \ref{claim:S_2}, it is enough to show that
\[
\widetilde{(R(qK_R)}|_{U} = \bigoplus_{d \in \Z} \cO_X(q(K_X+D')+dD)t^d.
\]
This follows from Proposition \ref{prop:S_2} and the case of $q=1$.
\end{proof}

\begin{proposition}\label{prop: canonical divisor}\textup{(cf.\ \cite{Watanabe91})}
We use Notation \ref{n-cone-QFS}. 
Fix $q \in \Z_{>0}$. 
Then the following hold. 
\begin{enumerate}
\item  
For $m \in \Z_{>0}$ and a nonzero element $f \in H^0(X, \MO_X(mD))$, 
we have the following isomorphism of $\Z$-graded $R$-modules: 
\[
R(qK_R)_{ft^m} \xrightarrow{\simeq} \bigoplus_{d \in \Z}  H^0(D_+(\overline{f}), \cO_X(q(K_X+D')+dD))t^d.
\]
\item  
The equality 
\[
\rho_*\cO_U(qK_U)=\bigoplus_{d \in \Z}\cO_X(q(K_X+D')+dD) t^d \subseteq K(X)[t, t^{-1}] 
\]
holds. 
Furthermore, it induces the isomorphism of $\Z$-graded $R$-modules
\[
H^i(U,\cO_U(qK_U)) \simeq \bigoplus_{d \in \Z} H^i(X,\cO_X(q(K_X+D')+dD))t^d.
\]
\end{enumerate}

\end{proposition}

\begin{proof}
We have the homomorphism of graded $R$-modules
\[
R(qK_R)_{ft^m} \to \bigoplus_{d \in \Z}  H^0(D_+(\overline{f}), \cO_X(q(K_X+D')+dD))t^d.
\]
By a similar argument to the proof of Proposition \ref{prop: affine morphism}, this is an isomorphism.
Then the assertions (1) and (2) follow from a similar argument to that of  Proposition \ref{prop: affine morphism}.
\end{proof}

\begin{remark}[Graded structure on $Q_{R, K_R, n}$]\label{rem:coh of Q}
We use Notation \ref{n-cone-QFS}. 
\begin{enumerate}
    \item 
    We equip $W_n(R)(qK_R)$ with the $p^{-(n-1)}\Z$-graded $W_n(R)$-module structure as in Proposition \ref{p-Witt-of-graded-can-div}. 
    \item 
    We define a $p^{-(n+r-1)}\Z$-module structure on $F_*^rW_{n}(R)(qK_R)$ by
    \[
    (F^r_*W_n(R)(qK_R))_e:=W_n(R)(qK_R)_{p^re}.
    \]
    Then Proposition \ref{p-Witt-of-graded-can-div} implies that the following maps are graded $W_n(R)$-module homomorphisms: 
    \begin{align*}
    &V \colon F_*W_{n-1}(R)(pqK_R) \to W_n(R)(qK_R),  \text{ and}\\
    &F \colon W_n(R)(qK_R) \to F_*W_n(R)(pqK_R). 
    \end{align*}
\item  
We introduce the $p^{-n}\Z$-graded $W_n(R)$-module structure on $Q_{R,K_R,n}$ via the 
following exact sequence (that is, the $e$-th graded piece is the quotient of the $e$-th graded pieces of the terms on the left): 
\[
0 \to F_*W_{n-1}(R)(pK_R) \xrightarrow{FV} F_*W_n(R)(pK_R) \to Q_{R,K_R,n} \to 0.
\]
\item
    It follows from (1) and (\ref{n-graded-str}) that $H^i(U,W_n\cO_U(qK_U))$ and 
    $H^{i+1}_\m(W_n(R)(qK_R))$ have  $p^{-(n-1)}\Z$-graded $W_n(R)$-module structures 
    such that the connecting map
    \[
    H^i(U,W_n\cO_U(qK_U)) \to H^{i+1}_\m(W_n(R)(qK_R))
    \]
    is a  graded $W_n(R)$-module homomorphism.
    Similarly, by (3) and (\ref{n-graded-str}), the cohomologies $H^i(U, Q_{U, K_U, n})$ and 
    $H^{i+1}_\m(Q_{R, K_R, n})$ have  $p^{-n}\Z$-graded $W_n(R)$-module structures 
    such that the connecting map
    \[
    H^i(U, Q_{U, K_U, n}) \to H^{i+1}_\m(Q_{R, K_R, n})
    \]
    is a  graded $W_n(R)$-module homomorphism.
\end{enumerate}
\end{remark}

\begin{proposition}\label{prop:coh of Q}
We use Notation \ref{n-cone-QFS}. 
Fix $d, n \in \Z_{>0}$. 
Then 
the following commutative diagram 
consists of graded $W_n(R)$-module homomorphisms of $p^{-n}\Z$-graded $W_n(R)$-modules
\[
\begin{tikzcd}
H^d(U,\omega_U) \arrow[r, "\Phi_{U, K_U, n}"] \arrow[d, "\simeq"] & H^d(Q_{U,K_U,n}) \arrow[d, "\simeq"] \\
H^{d+1}_\m(\omega_R) \arrow[r, "\Phi_{R, K_R, n}"] & H^{d+1}_\m(Q_{R,K_R,n}), 
\end{tikzcd}
\]
where all the maps are the natural ones and each vertical map is an isomorphism. 
\end{proposition}

\begin{proof}
The commutativity of the diagram follows from the functoriality of the connecting map. 
It is clear that the vertical maps are bijective. 
It suffices to show that each map is a graded $W_n(R)$-module homomorphism. 
As for the vertical arrows, we may apply Remark \ref{rem:coh of Q}(4). 
Recall that the horizontal maps are given by the right vertical maps in the following diagrams in which each horizontal sequences are exact 
\[
\begin{tikzcd}
0 \arrow[r] & F_*W_{n-1}\cO_U(pK_U) \arrow[r, "V"] \arrow[d, equal] & W_n\cO_U(K_U) \arrow[r] \ar[d, "F"] & \omega_U \arrow[r] 
\arrow[d, "\Phi_{U, K_U, n}"] & 0 \\
0 \arrow[r] & F_*W_{n-1}\cO_U(pK_U) \arrow[r, "FV"] & F_*W_n\cO_U(pK_U) \arrow[r] & Q_{U,K_U,n} \arrow[r] & 0
\end{tikzcd}
\]
and
\[
\begin{tikzcd}
0 \arrow[r] & F_*W_{n-1}(R)(pK_R) \arrow[r, "V"] \arrow[d, equal] & W_n(R)(K_R) \arrow[r] \arrow[d, "F"] & \omega_R \arrow[r] \arrow[d, "\Phi_{R, K_R, n}"] & 0 \\
0 \arrow[r] & F_*W_{n-1}(R)(pK_R) \arrow[r, "FV"] & F_*W_n\cO_R(pK_R) \arrow[r] & Q_{R,K_R,n} \arrow[r] & 0.     
\end{tikzcd}
\]
Hence we are done by Remark \ref{rem:coh of Q}(2). 
\end{proof}

\begin{proposition}\label{prop:commutativity}
We use Notation \ref{n-cone-QFS}. 
Set $d := \dim X$. 
Then we have the commutative diagram
\[
\begin{tikzcd}[column sep=1.5in]
H^d(U,\omega_U)_0 \arrow[r, "{H^d(U, \Phi_{U, K_U, n})_0}"]  & H^d(U, Q_{U,K_U,n})_0  \\
H^d(X,\cO_X(K_X+D')) \arrow[r, "{H^d(X, \Phi_{X, K_X+D', n})}"] \arrow[u, "\simeq"] & H^d(X,Q_{X,K_X+D',n}).
\arrow[u, "\simeq"]
\end{tikzcd}
\]
\end{proposition}

\begin{proof}
We take nonzero homogeneous elements 
$\overline{f}_1 := f_1t^{d_1},\ldots , \overline{f}_r := f_rt^{d_r}$ of $R$ 
such that $D(\overline{f}_1)\cup \cdots \cup D(\overline{f}_r)=U$. 
By Proposition \ref{prop: canonical divisor}(1), we obtain 
\begin{align*}
R(qK_R)_{\overline{f}_i}&=
    \bigoplus_{m \in \Z} H^0(D_+(\overline{f}_i), \cO_X(q(K_X+D')+mD))t^m\qquad\text{and}\\
\left(R(qK_R)_{\overline{f}_i}\right)_0&= 
H^0(D_+(\overline{f}_i), \cO_X(q(K_X+D'))) 
\end{align*}
for any $1 \leq i \leq r$. 
Then the latter equality implies 
\begin{equation}\label{e1:commutativity}
\left(W_n(R)(qK_R)_{[\overline{f}_i]}\right)_0
= H^0(D_+(\overline{f}_i), W_n\cO_X(q(K_X+D')))
\end{equation}
as $W_n(R)$-submodules of $W_n(K(R))$ by Proposition \ref{prop: canonical divisor}. 
Indeed, the left hand side is
\[
\left(W_n(R)(qK_R))_{[\overline{f}_i]}\right)_0= 
\prod_{m=0}^{n-1} (R(p^mqK_R)_{\overline{f}_i})_0 
\]
and the right hand side is
\[
H^0(D_+(\overline{f}_i), W_n\cO_X(q(K_X+D'))) 
=\prod_{m=0}^{n-1} H^0(D_+(\overline{f}_i), \cO_X(p^mq(K_X+D'))). 
\]
Therefore, (\ref{e1:commutativity}) induces the isomorphism of graded $W_n(R)$-modules
\[
\theta: H^d(X,W_n\cO_X(q(K_X+D')))\xrightarrow{\simeq} H^d(U,W_n\cO_U(qK_U) )_0,
\]
which commutes with $F$ and $V$. 
Therefore, we obtain the following commutative diagram of graded $W_n(R)$-modules 

{\small 
\[
\begin{tikzcd}[row sep=1.5em, column sep =0.8em]
 H^d(F_*W_{n-1}\cO_U(pK_U))_0
 & &  & 
 H^d(W_n\cO_U(K_U))_0
 & & \\
 & 
 H^d(F_*W_{n-1}\cO_U(pK_U))_0
 &  &  & 
 H^d(F_*W_n\cO_U(pK_U))_0
 &  & \\
 H^d(F_*W_{n-1}\cO_X(pE))
 & & & H^d(W_n\cO_X(E))
 & & \\
 & H^d(F_*W_{n-1}\cO_X(pE))
 & & & H^d(F_*W_n\cO_X(pE))
 &  
%
\arrow[from=1-1,to=2-2, equal]
\arrow[from=1-1,to=1-4, "V" pos=.6]
\arrow[from=1-4,to=2-5, "F"]
%
\arrow[from=3-1,to=4-2, equal]
\arrow[from=3-1,to=3-4, "V" pos=.6]
\arrow[from=4-2,to=4-5, "FV" pos=.4]
\arrow[from=3-4,to=4-5, "F"]
%
\arrow[from=3-1,to=1-1, "\theta" pos=.3, "\simeq"' pos=.3]
\arrow[from=3-4,to=1-4, "\theta" pos=.3, "\simeq"' pos=.3]
\arrow[from=4-5,to=2-5, "\theta" pos=.3, "\simeq"' pos=.3]
%
\arrow[from=4-2,to=2-2, crossing over, "\theta" pos=.3, "\simeq"' pos=.3]%
\arrow[from=2-2,to=2-5, crossing over, "FV" pos=.4]
\end{tikzcd}
\]}
\noindent where $E:=K_X+D'$.
Taking the cokernels of horizontal maps, we obtain the required diagram.
\end{proof}

\begin{theorem}\label{thm:height equality}
We use Notation \ref{n-cone-QFS}.  
Then, for every $n \in \Z_{>0}$,   
 $R_{\m}$ is $n$-quasi-$F$-split if and only if $(X, D')$ is $n$-quasi-$F$-split. 
In particular, $R_{\m}$ is quasi-$F$-split if and only if $(X, D')$ is quasi-$F$-split. 
\end{theorem}

\begin{proof}
Set $d := \dim X$. By Proposition \ref{prop:coh of Q} and Proposition \ref{prop:commutativity}, we obtain the following commutative diagram
\[
\begin{tikzcd}[column sep=3cm]
H^d(X,\cO_X(K_X+D')) \arrow[r, "H^d(\Phi_{X, K_X+D', n})"] \arrow[d, "\simeq"] & H^d(Q_{X,K_X+D',n}) \arrow[d, "\simeq"] \\
H^{d+1}_\m(\omega_R)_0 \arrow[r, "H^{d+1}_{\m}(\Phi_{R, K_R, n})_0"] \arrow[d, hookrightarrow] & H^{d+1}_\m(Q_{R,K_R,n})_0 \arrow[d, hookrightarrow] \\
H^{d+1}_\m(\omega_R) \arrow[r, "H^{d+1}_{\m}(\Phi_{R, K_R, n})"] &  H^{d+1}_\m(Q_{R,K_R,n}).    
\end{tikzcd}
\]

By \cite[Lemma 3.13]{KTTWYY1} (cf.\ Remark \ref{r-QFS}), 
the upper map 
$H^d(\Phi_{X, K_X+D', n})$  
is injective if and only if 
$(X, D')$ is $n$-quasi-$F$-split. 
Again by \cite[Lemma 3.13]{KTTWYY1}, 
the bottom map $H^{d+1}_{\m}(\Phi_{R, K_R, n})$ 
is injective if and only if $R_{\m}$ is $n$-quasi-$F$-split. 


By diagram chase, 
if the bottom map $H^{d+1}_{\m}(\Phi_{R, K_R, n})$ is injective, then the upper map 
$H^d(\Phi_{X, K_X+D', n})$ is injective. 
Conversely, assume that the upper map 
$H^d(\Phi_{X, K_X+D', n})$ is injective. 
Since the degree $0$ part 
$H^{d+1}_\m(\omega_R)_0$ 
of $H^{d+1}_\m(\omega_R)$ is the socle of $H^{d+1}_\m(\omega_R)$  (cf.\ \cite[(3.2)]{Watanabe91}), the bottom map 
$H^{d+1}_{\m}(\Phi_{R, K_R, n})$ 
is also injective.
\end{proof}

\begin{corollary}\label{cor:corresponding}
Let $k$ be an $F$-finite field of characteristic $p>0$. 
Let $X$ be a projective normal variety over $k$ with $H^0(X, \MO_X)=k$. 
Let $\Delta$ be an effective $\Q$-divisor on $X$ with standard coefficients 
such that $\rdown{\Delta}=0$.  
Assume that $-(K_X+\Delta)$ is ample.
Set 
\[
R :=\bigoplus_{d \geq 0} H^0(X,\cO_X(-d(K_X+\Delta)))t^d \subseteq K(R)[t] 
\]
and 
$\m := \bigoplus_{d > 0} H^0(X,\cO_X(-d(K_X+\Delta)))t^d \subseteq R$. 
Then 
$R_{\m}$ is quasi-$F$-split if and only if $(X, \Delta)$ is quasi-$F$-split. 
\end{corollary}

\begin{proof}
Set $d:=\dim X$ and $D := -(K_X+\Delta)$. 
We use Notation \ref{n-cone-QFS}. 
For the irreducible decomposition $\Delta = \sum_{i=1}^r \frac{d_i -1}{d_i}D_i$, 
we obtain 
\[
\{ D \} =  \{ -(K_X+\Delta)\} = \sum_{i=1}^r \frac{1}{d_i} D_i, 
\]
which implies  $D' = \sum_{i=1}^r \frac{d_i -1}{d_i}D_i =\Delta$ 
(for the definition of $D'$, see Notation \ref{n-cone-QFS}). 
It follows from Theorem \ref{thm:height equality} 
that $R_{\m}$ is quasi-$F$-split if and only if $(X, \Delta)$ is quasi-$F$-split. 
\end{proof}

%% file: section8.tex
\section{Klt threefolds and quasi-F-splitting}\label{s-klt3}


We start by generalising Theorem \ref{t-LDP-QFS} to the case when the base field is perfect. 

\begin{theorem} \label{t-LDP-QFS2}
Let $k$ be a perfect field of characteristic $p>41$. 
Let $(X,\Delta)$ be a log del Pezzo pair with standard coefficients. 
Then $(X,\Delta)$ is quasi-$F$-split.
\end{theorem}

\begin{proof}
If $k$ is algebraically closed, then the assertion follows from Theorem \ref{t-LDP-QFS}. 
The general case is reduced to this case by taking the base change to the algebraic closure \cite[Corollary 3.20]{KTTWYY1}. 
\end{proof}

We are ready to prove main theorems of this paper: 
Theorem \ref{t-3dim-klt-QFS} and Theorem \ref{t-3dim-klt-nonQFS}. 

\begin{theorem} \label{t-3dim-klt-QFS}
Let $k$ be a perfect field of characteristic $p >41$. 
Let $(X, \Delta)$ be a three-dimensional $\Q$-factorial affine klt pair 
over $k$, where $\Delta$ has standard coefficients. 
Then $(X, \Delta)$ is quasi-$F$-split. 
In particular, $X$ lifts to $W_2(k)$. 
\end{theorem}

\begin{proof}
The same argument as in \cite[Theorem 6.19]{KTTWYY1} works 
after replacining  \cite[Theorem 6.18]{KTTWYY1} by Theorem \ref{t-LDP-QFS}. 
\end{proof}

\begin{theorem} \label{t-3dim-klt-nonQFS}
Assume $p \in \{11, 17, 19, 23, 41\}$. 
Let $k$ be an algebraically closed field of characteristic $p$. 
Then there exists a three-dimensional affine $\Q$-factorial klt variety $V$ which is not quasi-$F$-split. 
\end{theorem}

\begin{proof}
Let $(X:=\bP^2_k, \Delta)$ be as in Notation \ref{n-P^2-cex}. 
Set $D:=-(K_X+\Delta)$ and $V := V_{X, D}$. 
By Theorem \ref{t-klt}, 
$V$ is a three-dimensional affine  $\Q$-factorial klt variety. 
It follows from Theorem \ref{t-P^2-main} that $(X, \Delta)$ is not quasi-$F$-split. 
Then Corollary \ref{cor:corresponding} implies that $V$ is not quasi-$F$-split. 
\end{proof}



%% file: section9.tex
\section{Application to extension theorem for differential forms}

In this section, we apply 
Theorem \ref{t-3dim-klt-QFS} 
to prove the logarithmic extension theorem for one-forms on three-dimensional terminal singularities (Corollary \ref{c-3dim-1form}).

\begin{definition}[Logarithmic extension theorem]\label{def:log ext thm}
Let $X$ be a normal variety over a perfect field and let $D$ be a reduced divisor on $X$.
We say that $(X,D)$ satisfies \textit{the logarithmic extension theorem for $i$-forms} if, for every proper birational morphism $f\colon Y\to X$ from a normal variety $Y$, the natural 
restriction injection 
\[
f_{*}\Omega^{[i]}_Y(\log\,(E+f_{*}^{-1}D))\hookrightarrow \Omega^{[i]}_X(\log\,D)
\]
is an isomorphism, where $E$ is the largest reduced $f$-exceptional 
divisor. 
\end{definition}

\begin{remark}
Let $i \in \Z_{\geq 0}$. 
\begin{enumerate}
    \item If $(X,D)$ is a log canonical pair over $\mathbb{C}$, then $(X,D)$ satisfies  the logarithmic extension theorem for $i$-forms  (\cite[Theorem 1.5]{GKKP}).
    \item If $(X,D)$ is a two-dimensional log canonical pair over a perfect field of characteristic $p>5$, then $X$ satisfies the logarithmic extension theorem for $i$-forms (\cite[Theorem 1.2]{graf21}).
\end{enumerate}
\end{remark}

\begin{definition}\label{def:reflexive Carter operators}
Let $X$ be a normal variety over a perfect field 
of characteristic $p>0$ and let $D$ be a reduced divisor on $X$.
Let $j\colon U\to X$ be the inclusion of the log smooth locus $U$ of $(X,D)$. 
For $D_U \coloneqq D|_U$, set 
\begin{align*}
    &B_n\Omega_X^{[i]}(\log\,D)\coloneqq j_{*}B_n\Omega^{{\cred i}}_U(\log\,D_U)\quad \text{and}\\
    &Z_n\Omega_X^{[i]}(\log\,D)\coloneqq j_{*}Z_{n}\Omega_U^{i}(\log\,D_U).
\end{align*}
We define \textit{the $i$-th reflexive Cartier operator} by
\begin{align*}C^{[i]}_{X,D}:=j_{*}C^{i}_{U,D}\colon Z_1\Omega_X^{[i]}(\log\,D_U)\to \Omega^{[i]}_X(\log\,D_U).
\end{align*}
\end{definition}

We have the following criterion for the logarithmic extension theorem in positive characteristic.

\begin{theorem}[\textup{\cite[Theorem A]{Kaw4}}]\label{thm:criterion for log ext thm}
Let $X$ be a normal variety over a perfect field of characteristic $p>0$ and let $D$ be a reduced divisor on $X$. 
Fix $i \in \Z_{\geq 0}$. 
If the $i$-th reflexive Cartier operator 
\[
C^{[i]}_{X,D}\colon Z_1\Omega_X^{[i]}(\log\,D)\to \Omega_X^{[i]}(\log\,D)
\]
is surjective, then $(X,D)$ satisfies the logarithmic extension theorem for $i$-forms.
\end{theorem}

\begin{theorem}\label{thm:log ext thm for qFs}
Let $X$ be a normal variety 
over a perfect field of characteristic $p>0$. 
Let $Z$ be the singular locus of $X$ and $j\colon U:=X\setminus Z\hookrightarrow X$ the inclusion.
Suppose that $X$ is quasi-$F$-split, $X$ satisfies $S_3$, and $\codim_X(Z)\geq 3$.
Then the first reflexive Cartier operator
\[
C^{[1]}_{X}\colon Z_{1}\Omega_X^{[1]}\to \Omega^{[1]}_X
\]
is surjective.
In particular, $X$ satisfies the logarithmic extension theorem for one-forms.
\end{theorem}

\begin{proof}
By Theorem \ref{thm:criterion for log ext thm}, it suffices to prove the surjectivity of $C_X^{[1]}$. 
To this end, we may assume that $X$ is affine.
Since $F_{*}\sO_X$ satisfies $S_3$ and $\codim_X(Z)\geq 3$, it follows from \cite[Proposition 1.2.10 (a) and (e)]{BH93} that $H_Z^{2}(F_{*}\sO_X)=0$.

We define $\sO_X$-modules $\mathcal{B}_m$ by $\mathrm{Coker}(F\colon W_m\sO_X\to F_{*}W_m\sO_X)$.
Then $B_m\Omega_X^{[1]}$ and $\mathcal{B}_m$ coincide with each other on the smooth locus of $X$ by Lemma \ref{lem:Serre's map}.
Since $X$ is quasi-$F$-split, we can take $n\in\Z_{>0}$ such that 
\[
0 \to \sO_X\to Q_{X,n}\to \mathcal{B}_n\to 0 
\] 
splits.
Then we have the following commutative diagram
\[
\begin{tikzcd}
H^2_{Z}(Q_{X,n})\arrow[r, twoheadrightarrow]\arrow[d] & H^2_{Z}(\mathcal{B}_n)\arrow[d]\\
  \mathllap{0\,=\,\, } H^2_{Z}(F_{*}\sO_X)\arrow[r]& H^2_{Z}(\mathcal{B}_1),
\end{tikzcd}
\]
where the right vertical arrow induced by the $n$-times iterated Cartier operator $(C_U^{1})^{n}$ by Lemma \ref{l-Bn-Cartier} and Remark \ref{r-QFS}.
The top horizontal arrow is surjective since $Q_{X,n}\to \mathcal{B}_n$ is a splitting surjection.
Thus, the above diagram shows that the right vertical map is zero.
Since $X$ is affine, we have 
\[
H^0(X, R^1j_{*}B_n\Omega_U^{1})\cong H^1(U, B_n\Omega_U^{1})\cong H^2_{Z}(\mathcal{B}_n),
\]
and thus, $R^1j_{*}(C_U^{1})^{n-1}\colon R^1j_{*}B_n\Omega_U^{1} \to R^1j_{*}B_1\Omega_U^{1}$ is zero.
By the following commutative diagram
\[
\begin{tikzcd}[column sep=3cm]
Z_{n}\Omega_X^{[1]}\arrow[r, "(C_{X}^{[1]})^{n}"]\arrow[d, "(C_{X}^{[1]})^{n-1}"'] & \Omega_X^{[1]}\ar[r]\arrow[d, equal] & R^1j_{*}B_{n}\Omega_U^{1}\arrow[d, "R^1j_{*}(C_U^{1})^{n-1}=\,0"]\\
 Z_{1}\Omega_X^{[1]}\arrow[r, "C_X^{[1]}"] & \Omega_X^{[1]}\arrow[r] & R^1j_{*}B_{1}\Omega_U^{1},
\end{tikzcd}
\]
where horizontal exact sequences are obtained by (\ref{exact1}).
Now, by the above diagram, we obtain the surjectivity of $C_{X}^{[1]}$, as desired.
\end{proof}

\begin{corollary}\label{c-3dim-1form}
Let $X$ be a three-dimensional variety over a perfect field of characteristic $p>41$. 
Assume that one of the following holds. 
\begin{enumerate}
\item $X$ is terminal. 
\item 
The singular locus of $X$ is zero-dimensional and $X$ is $\Q$-factorial and klt. 
\end{enumerate}
Then, $X$ satisfies the logarithmic extension theorem for one-forms.
\end{corollary}

\begin{proof}
We first treat the case when (2) holds. 
In this case, $X$ is Cohen--Macaulay by \cite[Corollary 1.3]{ABL20}. 
Then the assertion follows from 
Theorem \ref{t-3dim-klt-QFS} 
and Theorem \ref{thm:log ext thm for qFs}. 
This completes the proof for the case when (2) holds. 

Assume (1). 
By taking a $\Q$-factorialisation of $X$ (see, e.g., \cite[Theorem 2.14]{GNT16}), we may assume that $X$ is $\Q$-factorial. 
Since $X$ is terminal, it follows from \cite[Corollary 2.13]{kollar13} that $X$ has only isolated singularities. 
Then (2) holds. 
\end{proof}

%% file: section10.tex
\section{Appendix: Cone construction for $\Q$-divisors}\label{s-Qcone}

Throughout this section, 
we work over an arbitrary  field $k$. 
All the results and proofs in this section work in any characteristic.
Let $X$ be a projective normal variety $X$ and let $D$ be an ample $\Q$-Cartier $\Q$-divisor. 
We shall summarise some foundational results on the {\emph{orbifold cone}}
\[
V_{X, D} := \Spec \left(\bigoplus_{d \geq 0} H^0(X, \MO_X(dD)){t^d}
\right). 
\]
In order to compare $X$ and $V_{X, D}$, we also introduce: 
\begin{alignat*}{3}
&\mathbb A^1\textrm{-fibration:} \quad &&\pi : W_{X, D} := \Spec_X \Big( \bigoplus_{d \geq 0} \MO_X(dD)t^d\Big) \to X, \textrm{ and }\\
&\mathbb G_m\textrm{-fibration:} \quad &&\pi^{\circ} : W^{\circ}_{X, D} := \Spec_X \Big( \bigoplus_{d \in \bZ} \MO_X(dD)t^d\Big) \to X.
\end{alignat*}
Note that if $D$ is a very ample Cartier divisor, then all the results in this section are well known. 
In this special case, $\pi : W_{X, D} \to X$ is an $\mathbb A^1$-bundle, {$\pi^{\circ} : W^{\circ}_{X, D} \to X$ is a $\mathbb G_m$-bundle}, and $V_{X, D}$ is the {usual} cone.

\subsection{Foundations}

\begin{definition}\label{d-Qcone1}
Let $X$ be a normal variety 
and let $D$ be a $\Q$-Cartier $\Q$-divisor. 
For the quasi-coherent finitely generated $\MO_X$-algebras
\[
\mathcal A_{X, D} := \bigoplus_{d \geq 0} \MO_X(dD)t^d \qquad 
{\rm and} \qquad 
\mathcal A^{\circ}_{X, D} := \bigoplus_{d \in \Z} \MO_X(dD)t^d, 
\]
we set 
\[
W_{X, D} := \Spec_X \mathcal A_{X, D} \qquad 
{\rm and} \qquad  W^{\circ}_{X, D} := \Spec_X \mathcal A^{\circ}_{X, D}.
\]
Note that 
$\mathcal A_{X, D}$ and $\mathcal A^{\circ}_{X, D}$ 
are sheaves of graded subrings of the constant sheaves 
$\bigoplus_{d \geq 0}K(X)t^d$ and $\bigoplus_{d \in \Z}K(X)t^d$, respectively. 
We have the induced morphisms: 
\[
\pi^{\circ} : W^{\circ}_{X, D} \xrightarrow{j} W_{X, D} \xrightarrow{\pi} X. 
\]
Let $\Gamma_{X, D} \subseteq W_{X, D}$ be the {\em $0$-section}, i.e., 
$\Gamma_{X, D}$ is the section of $\pi : W{_{X, D}}\to X$ corresponding to the ideal sheaf 
\[
(\mathcal A_{X, D})_+ := \bigoplus_{d > 0}\MO_X(dD)t^d
\]
of $\mathcal A_{X, D} = \bigoplus_{d \geq 0}\MO_X(dD)t^d$. 
We drop the subscript $(-)_{X, D}$ when no confusion arises, for example, $W := W_{X, D}$. 
\end{definition}



\begin{lemma}\label{l-VW-normal}
Let $X$ be a normal variety and let $D$ be a $\Q$-Cartier $\Q$-divisor. 
Then the following hold. 
\begin{enumerate}
    \item $W_{X, D}$ is a normal variety. 
    \item 
    The induced morphism $j: W^{\circ}_{X, D} \to W_{X, D}$ is an open immersion. 
    \item The set-theoretic equality $W_{X, D} \setminus \Gamma_{X, D} = j(W^{\circ}_{X, D})$ holds. 
\end{enumerate}
\end{lemma}

\begin{proof}
 The assertion (1) follows from \cite[3.1 in page 48]{Dem88}. 
The assertions (2) and (3) hold by \cite[Lemme 2.2(i) in page 40]{Dem88}. 
\end{proof}

In what follows, we consider $W^{\circ}_{X, D}$ as an open subscheme of $W_{X, D}$.

\begin{definition}\label{d-Qcone2}
Let $X$ be a projective normal variety and let $D$ be an ample $\Q$-Cartier $\Q$-divisor. 
For the finitely generated $k$-algebra 
\[
R(X, D) := \bigoplus_{d \geq 0} H^0(X, \MO_X(dD)t^d, 
\]
we set 
\[
V_{X, D} := \Spec R(X, D). 
\]
Let $v_{X, D}$  be the {\em vertex} of $V_{X, D}$, that is, 
$v_{X, D}$ is the closed point of $V_{X,D}$ corresponding to 
the maximal ideal $\bigoplus_{d > 0} H^0(X, \MO_X(dD))t^d$ 
of $R(X, D) = \bigoplus_{d \geq 0} H^0(X, \MO_X(dD))t^d$. 
By the ring isomorphism 
\[
\Gamma(V_{X, D}, \MO_{V_{X, D}}) =  
\bigoplus_{d \geq 0} H^0(X, \MO_X(dD)t^d) \xrightarrow{\simeq}
H^0(X, \bigoplus_{d \geq 0}\MO_X(dD)t^d) = \Gamma(W_{X, D}, \MO_{W_{X, D}}), 
\]
we obtain a morphism $\mu: W_{X, D} \to V_{X, D}$
with $\mu_*\MO_{W_{X, D}}=\MO_{V_{X, D}}$ (cf.\ \cite[Ch. II, Exercise 2.4]{hartshorne77}); {here the latter condition can be checked on global sections as $V_{X,D}$ is affine and $\mu_*\MO_{W_{X, D}}$ is quasi-coherent}. 
To summarise, we have the following {diagram}:  
\[
\begin{tikzcd}
	W^{\circ}_{X, D} \arrow[rd, "\pi^{\circ}"']& W_{X, D} & V_{X, D} \\
	& X. 
	\arrow["j", hook, from=1-1, to=1-2]
	\arrow["\mu", from=1-2, to=1-3]
	\arrow["\pi", from=1-2, to=2-2]
\end{tikzcd}
\]
We drop the subscript $(-)_{X, D}$ when no confusion arises, for example, $V:= V_{X, D}$. 
\end{definition}

\begin{theorem}\label{t-Qcone-birat}
Let $X$ be a projective normal variety and 
let $D$ be an ample $\Q$-Cartier $\Q$-divisor. 
Then the following hold. 
\begin{enumerate}
    \item $V_{X, D}$ is an affine normal variety. 
    \item  the induced morphism $\mu : W_{X, D} \to V_{X, D}$ 
is a projective birational morphism such that the following set-theoretic equalities are valid:
\begin{equation}\label{e1-Qcone-birat}
\Ex(\mu) = \Gamma_{X, D} 
\qquad {\rm and} \qquad 
\mu(\Ex(\mu)) = \mu(\Gamma_{X, D})=v_{X, D}.
\end{equation}
\end{enumerate}
\end{theorem}

\begin{proof}
By $(\mu_{X, D})_*\MO_{W_{X, D}} = \MO_{V_{X, D}}$, 
it follows from Lemma \ref{l-VW-normal} that $V$ is an affine normal variety. 
Thus (1) holds. 
The assertion (2) follows from \cite[3.4 in page 48]{Dem88}. 
\end{proof}

\subsection{Functoriality}

Let $\alpha : Y \to X$ be a dominant morphism of normal varieties. 
Let $D$ be a $\Q$-Cartier $\Q$-divisor and set $D_Y := \alpha^*D$. 
Then we have the following commutative diagram (for the definitions of $W_{X, D}$ and $W_{Y, D_Y}$, see Definition \ref{d-Qcone1}): 
\[
\begin{tikzcd}
Y \arrow[d, "\alpha"'] & W_{Y, D_Y} \arrow[l, "\pi_Y"']\arrow[d, "\beta"]\\
X & W_{X, D}. \arrow[l, "\pi"]
\end{tikzcd}
\]
Furthermore, if $\alpha : Y \to X$ is a finite  surjective morphism of projective normal varieties and $D$ is an ample $\Q$-Cartier $\Q$-divisor, then we get the following commutative diagram  (for the definitions of $V_{X, D}$ and $V_{Y, D_Y}$, see Definition \ref{d-Qcone2}): 
\[
\begin{tikzcd}
    Y \arrow[d, "\alpha"'] & W_{Y, D_Y} \arrow[l, "\pi_Y"']\arrow[d, "\beta"] \arrow[r, "\mu_Y"] & V_{Y, D_Y} \arrow[d, "\gamma"]\\
X & W_{X, D} \arrow[l, "\pi"] \arrow[r, "\mu"'] & V_{X, D}.
\end{tikzcd}
\]

\begin{lemma}\label{l-funct-Gamma}
Let $\alpha : Y \to X$ be a dominant morphism of normal varieties. 
Let $D$ be a $\Q$-Cartier $\Q$-divisor on $X$ and set $D_Y := \alpha^*D$. 
Then the set-theoretic equality 
\[
\beta^{-1}(\Gamma_{X, D}) = \Gamma_{Y, D_Y}
\]
holds for the induced morphism $\beta : W_{Y, D_Y} \to W_{X, D}$. 
\end{lemma}

\begin{proof}
Since the problem is local on $X$ and $Y$, 
we may assume that $X=\Spec R$, $Y = \Spec R_Y$, and $\MO_X(d_0D)|_{\Spec R} \simeq \MO_{\Spec R}$ for some $d_0 \in \Z_{>0}$. 
There exists $f \in K(X) \setminus \{0\}$ such that 
\[
R(d_0dD)= f^dR \qquad {\rm and} \qquad 
R_Y(d_0dD)= f^dR_Y
\]
for every $d \in \Z$. 
We have 
the induced graded
ring homomorphism: 
\[
A := \bigoplus_{d \geq 0} R(dD)t^d \to \bigoplus_{d \geq 0} R_Y(dD_Y)t^d=:A_Y. 
\]
Since this ring homomorphism is injective, 
we consider $A$ as a subring of $A_Y$. 
It suffices to show that 
\[
(A_Y)_+ = \sqrt{A_+ \cdot A_Y},
\]
where $A_+ :=  \bigoplus_{d >0} R(dD)t^d$ and 
$(A_Y)_+ :=  \bigoplus_{d >0} R_Y(dD_Y)t^d$. 
By $(A_Y)_+ \supseteq A_+$ and  $\sqrt{(A_Y)_+} =(A_Y)_+$, 
it holds that 
\[
(A_Y)_+ =\sqrt{(A_Y)_+}  \supseteq 
\sqrt{A_+ \cdot A_Y}. 
\]

It is enough to prove the opposite inclusion: 
$(A_Y)_+  \subseteq \sqrt{A_+ \cdot A_Y}$.  
Fix $d \in \Z_{>0}$. 
Take a homogeneous element $\psi t^d \in 
R_Y(dD_Y)t^d$ with $\psi \in R_Y(dD_Y)$. 
We obtain 
\[
\psi^{d_0} \in R_Y(d_0d{D_Y}) = f^d R_Y. 
\]
Hence we can write $\psi^{d_0} =f^d s$ for some 
$s \in R_Y$. 
Then it holds that 
\[
(\psi t^d)^{d_0} = 
\psi^{d_0} t^{d_0d} 
= f^ds \cdot t^{d_0d} = (ft^{d_0})^d \cdot s 
\in A_+ \cdot A_Y, 
\]
which implies $\psi t^{d} \in \sqrt{A_+ \cdot A_Y}$, as required. 
\end{proof}

\begin{lemma}\label{l-Qcone-etale}
Let $\alpha : Y \to X$ be a smooth morphism of normal varieties. 
Let $D$ be a $\Q$-Cartier $\Q$-divisor on $X$ and set $D_Y := \alpha^*D$. 
Then the following diagram
\[
\begin{tikzcd}
Y \arrow[d, "\alpha"'] & W_{Y, D_Y} \arrow[l, "\pi_Y"']\arrow[d, "\beta"]\\
X & W_{X, D}. \arrow[l, "\pi"]
\end{tikzcd}
\]
is cartesian, where $\beta, \pi, \pi_Y$ are the induced morphisms. 
\end{lemma}

\begin{proof}
Since both $W_{Y, D_Y}$  and $W_{X, D} \times_X Y$ are affine over $Y$, 
the problem is local on $Y$. 
In particular, we may assume that $X= \Spec R$ and $Y= \Spec R_Y$. 
It suffices to show that 
\[
\left( \bigoplus_{d \geq 0} R(dD) \right) \otimes_R R_Y \to  \bigoplus_{d \geq 0} R_Y(dD_Y)
\]
is an isomorphism. 

Fix $d \in \Z_{\geq 0}$. 
Set $E := \llcorner d D \lrcorner$ and $E_Y := \llcorner dD_Y \lrcorner$. 
We have 
\[
E_Y = \llcorner dD_Y \lrcorner = \llcorner \alpha^*(dD) \lrcorner = 
\alpha^*(\llcorner dD \lrcorner) = \alpha^* E, 
\]
where the third equality holds, because $\alpha$ is smooth. 
Then it is enough to prove that 
\[
\alpha^*\MO_X(E) \to \MO_Y(\alpha^*E) 
\]
is an isomorphism. 
Outside the singular locus, these coincide. 
Hence it suffices to show that both hand sides are reflexive. 
It is well known that the right hand side $\MO_Y(\alpha^*E)$ is reflexive. 
For the smooth locus $X'$ of $X$ and its inverse image $Y' := \alpha^{-1}(X')$, 
we have the cartesian diagram: 
\[
\begin{tikzcd}
    Y' \arrow[r, "i"]\arrow[d, "\alpha'"'] & Y \arrow[d, "\alpha"]\\
    X' \arrow[r, "j"] & X.
\end{tikzcd}
\]
By the isomorphism 
\[
j_*j^*\MO_X(E) \xrightarrow{\simeq} \MO_X(E), 
\]
we obtain 
\[
\alpha^*
\MO_X(E)
\xleftarrow{\simeq}
\alpha^*j_*j^*\MO_X(E)
\simeq i_*\alpha'^*j^*\MO_X(E) 
\simeq i_*i^* \alpha^*\MO_X(E), 
\]
where the isomorphism 
$\alpha^*j_*j^*\MO_X(E)
\simeq i_*\alpha'^*j^*\MO_X(E)$ follows from the flat base change theorem. 
Therefore, $\alpha^*\MO_X(E)$ is reflexive, as required. 
\end{proof}



\subsection{Miscellanies} 

\begin{proposition}\label{prop:S_2}
    Let $X$ be a normal variety with $\dim X \geq 1$, 
    {let} $D$ {be} an ample $\Q$-Cartier $\Q$-divisor, 
    and 
    {let} $L$ {be} a $\Q$-divisor.
    Let $\mathcal{F}_L$ be the quasi-coherent $\MO_{W^{\circ}_{X,D}}$-submodule of the constant sheaf $K(X)[t,t^{-1}]$ corresponding to the quasi-coherent $\pi^{\circ}_*\MO_{W^{\circ}_{X,D}}$-module
    \[
    \bigoplus_{d \in \Z} \cO_X(L+dD)t^d
    \]
    (see \cite[Exercise II.5.17(e)]{hartshorne77}),
    that is, $\mathcal{F}_L$ satisfies 
    \[
    \pi^{\circ}_*\mathcal{F}_L \simeq \bigoplus_{d \in \Z} \cO_X(L+dD)t^d.
    \]
    {Then the following hold.}
    \begin{enumerate}
        \item If $X$ is regular or $D$ is Cartier, then $\pi^{\circ} \colon W_{X,D}^{\circ} \to X$ is flat.
        \item $\mathcal{F}_{L}$ is coherent and satisfies the condition $S_2$.
        \item For every codimension one point $w \in W^{\circ}_{X,D}$, the codimension of $\pi^{\circ}(w)$ is one or zero.
        \item Let $D=\sum^r_{i=1} \frac{l_i}{d_i} D_i$ be the irreducible decomposition, where $l_i$ and $d_i$ are coprime integers satisfying $d_i >0$ for each $1 \leq i \leq r$. Set $ D':=\sum^r_{i=1}\frac{d_i-1}{d_i}D_i$. 
        Then, for every $q \geq 1$  {\cred and every Weil divisor $K$}, we have
        \[
        (\mathcal{F}_{{\cred K+}D'}^{\otimes q})^{**}=\mathcal{F}_{q({\cred K+}D')}
        \]
        via the natural inclusions into $K(X)[t,t^{-1}]$, where $(-)^{**}$ denotes the reflexive hull.
    \end{enumerate}
\end{proposition}
\begin{proof}
    By definition,
    \[
    \pi^{\circ}_*\cO_{W^{\circ}_{X,D}} =\bigoplus_{d \in \Z} \cO_X(dD)t^d.
    \]
    If $X$ is regular or $D$ is Cartier,  then $\pi^{\circ}_*\cO_{W^{\circ}_{X,D}}$ is a locally free $\cO_X$-module.
    Since $\pi$ is affine, $\pi$ is flat.
    Next, we prove the assertion (2).
    Let $N$ be a positive integer such that $ND$ is Cartier.
    Then the natural inclusion
    \[
    \bigoplus_{d \in \Z} \cO_X(dND)s^d \to \bigoplus_{d \in \Z} \cO_X(dD) t^d
    \qquad\qquad s \mapsto t^N
    \]
    induces a finite surjective morphism $f \colon W^{\circ}_{X,D} \to W^{\circ}_{X,ND}$.
    We have
    \[
    (\pi^{\circ}_N)_*f_*\mathcal{F}_{L} \simeq \bigoplus_{0 \leq i \leq N-1} \cO_X(L+iD) \otimes_{\cO_X} \cO_{W^{\circ}_{X,ND}} \simeq \bigoplus_{0 \leq i \leq N-1} \cO_{X}(L+iD) \otimes_{\cO_X} \mathcal{A}^{\circ}_{X,ND},
    \]
    where $\pi^{\circ}_N \colon W^{\circ}_{X,ND} \to X$ is a canonical morphism induced in Definition \ref{d-Qcone1}.
    In particular,
    \[    (\pi^{\circ}_N)^*\bigoplus_{0 \leq i \leq N-1} \cO_{X}(L+iD) \simeq f_*\mathcal{F}_L.
    \]
    Since $f_*\mathcal{F}_L$ is coherent, so is $\mathcal{F}_L$. 
    By the assertion (1), $\pi^{\circ}_N$ is flat.
    Furthermore, since $\cO_X(L+iD)$ satisfies the condition $S_2$, so does $\mathcal{F}_{L}$.
    Next, we prove the assertion (3).
    Since the above $f$ is finite and surjective, $f(w)$ is also a codimension one point.
    Therefore, the assertion follows from the flatness of $\pi^{\circ}_N$.
    Indeed, since $\pi^{\circ}_N$ is flat, every fiber of $\pi^{\circ}_N$ has dimension one.
    By \cite[III.Proposition 9.5]{hartshorne77}, we have
    \[
    \dim \cO_{X,\pi^{\circ}_N(f(w))} \leq \dim \cO_{W^{\circ}_{X,ND},f(w)}=1.
    \]
    Finally, we prove the assertion (4).
    First, we describe the natural inclusions.
    We consider the composition of maps
    \[
    \mathcal{F}^{\otimes q}_{{\cred K+}D'} \hookrightarrow K(X)[t,t^{-1}]^{\otimes q} \simeq K(X)[t,t^{-1}],
    \]
    where the first inclusion induced by the natural injection and the second isomorphism induced by taking products.
    Therefore, it induces the natural injection
    \[
    (\mathcal{F}_{{\cred K+}D'}^{\otimes q})^{**} \hookrightarrow K(X)[t,t^{-1}].
    \]
    By the assertion (3) and the $S_2$-condition of $\mathcal{F}_{{\cred q(K+D')}}$, we may assume that $X$ is a spectrum of a discrete valuation ring.
    In particular, {\cred $K=cE$ for some integer $c$, 
    and} $D=\frac{a}{b} E$ is the irreducible decomposition, where $a$ and $b$ are coprime integers satisfying 
    $b > 0$ and $D'=\frac{b-1}{b} E$.
    For the assertion, it is enough to show that for every $d \in \Z$, the map
    \[
     \bigoplus_{d_1+\cdots+d_q=d}( \cO_X({\cred K+}D'+d_1D) \otimes \cdots \otimes \cO_X({\cred K+}D'+d_qD)) \to \cO_X(q{\cred (K+D')}+dD)
    \]
    induced by taking the products is surjective.
    Therefore, it is enough to show that for every $d \in \Z$, there exists $d_1,\ldots,d_q \in \Z$ such that $d_1+\cdots +d_q=d$ and
    \[
    \sum_{1 \leq i \leq q}\lfloor {\cred c+} \frac{b-1}{b}+d_i\frac{a}{b} \rfloor =\lfloor q {\cred (c+\frac{b-1}{b})}+d\frac{a}{b} \rfloor.
    \]
    Since $a$ and $b$ are coprime integers, we can find $x, y \in \Z$ such that
    \[
    (b-1)+ax \equiv q(b-1)+da \mod b,\ \text{and}\ (b-1)+ay \equiv 0 \mod b.
    \]
    Then $x+(q-1)y \equiv d \mod b$.
    Set $d_1:=d-(q-1)y$ and $d_i:=y$ for $2 \leq i \leq q$. Then  $d_1+\cdots +d_q=d$, and
    \[
    (b-1)+d_1a \equiv (b-1)+ax \equiv q(b-1)+da \mod b.
    \]
    Therefore, if we set
    \begin{align*}
        (b-1)+ad_1&= bQ_1+r, \\
        (b-1)+ad_i&= bQ_i\ &\text{for $i \geq 2$}, \\
        q(b-1)+ad &= bQ+r,
    \end{align*}
    then $Q_1+\cdots +Q_q=Q$.
    Furthermore, we have
    \[
    \sum_{1 \leq i \leq q}\lfloor {\cred c+}\frac{b-1}{b}+d_i\frac{a}{b} \rfloor
        = {\cred qc+} Q_1+\cdots +Q_q =Q=\lfloor q {\cred (c+\frac{b-1}{b})}+d\frac{a}{b} \rfloor.
    \]
\end{proof}

%% file: section11.tex
\section{Appendix: Cone singularities for log smooth pairs}\label{s-Qcone-klt}

Throughout this section, 
we work over an algebraically closed field $k$ of arbitrary characteristic. However, the main results of this section should extend to perfect fields by base change.

Let $X$ be a smooth projective variety and let $\Delta$ be a simple normal crossing $\Q$-divisor. 
Assume that $\rdown{\Delta}=0$, $\Delta$ has standard coefficients, and $D :=-(K_X+\Delta)$ is ample.  
The purpose of this section is to prove that the {orbifold} cone $V := V_{X, D}$ is $\Q$-factorial and klt (Theorem \ref{t-klt}). 
For the $0$-section $\Gamma := \Gamma_{X, D}$ of $\pi : W:=W_{X, D} \to X$, 
we shall prove (I) and (II) below. 
\begin{enumerate}
\item[(I)] $W$ is $\Q$-factorial, $(W, \Gamma)$ is plt, and {so} $(W, 0)$ is klt. 
\item[(II)] $(K_W + \Gamma)|_{\Gamma} = K_{\Gamma}+\Delta_{\Gamma}$, where $\Delta_{\Gamma} := (\pi|_{\Gamma})^*\Delta$ via 
$\pi|_{\Gamma} : \Gamma \hookrightarrow W \to X$. 
\end{enumerate}
The main tool is the toric geometry. 
Both (I) and (II) are \'etale local problems on $X$. 
Then the proofs will be carried out by introducing a toric structure on $(W, \Gamma)$ for the case when $X=\mathbb A^n$ and $\Supp \Delta$ is a union of coordinate hyperplanes. 
Subsection \ref{ss-toric-str} and Subsection \ref{ss-Diff} are devoted to proving (I) and (II), respectively. 
In Subsection \ref{ss-notation-toric}, we recall some foundations on toric geometry. 



\subsection{Notation and results on toric varieties}\label{ss-notation-toric}

When we work with toric varieties, we use the following notation, 
which is extracted from \cite{CLS11}. 
In this paper, we only need normal affine toric varieties $U_{\sigma}$ associated to a rational {strongly} 
 convex polyhedral cone $\sigma$. 

\begin{notation}\label{n-toric}
\hphantom{a}\\
\vspace{-1em}
\begin{enumerate}
    \item Let $N$ be a finitely generated free $\Z$-module. Let $M$ be its dual, 
    that is, $M := \Hom_{\Z}(N, \Z)$. 
    \item 
    Set $N = M = \Z^r$. 
    We identify $M$ and $\Hom_{\Z}(N, \Z)$ via the standard inner product: 
    \[
    N \times M \to \Z, \qquad (x, y) \mapsto \langle x, y\rangle, 
    \]
    that is, 
    \[
    \langle (x_1, ..., x_r), (y_1, ..., y_r)\rangle := x_1y_1 + \cdots +x_ry_r. 
    \]
    \item Set $N_{\R} := N \otimes_{\Z} \R$ and $M_{\R} := M \otimes_{\Z} \R$. 
    \item For a rational {strongly} convex polyhedral cone $\sigma \subseteq N_{\R}$, 
    let $U_{\sigma}$ be the associated affine normal toric variety, that is, 
    $U_{\sigma} = \Spec k[M \cap \sigma^{\vee}]$. 
\end{enumerate}
\end{notation}

\begin{definition}\label{d-toric}
\hphantom{a}\\
\vspace{-1em}
\begin{enumerate}
\item We say that $T$ is a {\em torus} if $T =\mathbb G_m^n$, 
which is an algebraic group over $k$. 
\item 
We say that $(T, X)$ is a {\em toric variety} if 
$X$ is a variety and $T$ is an open subscheme of $X$ such that 
\begin{itemize}
    \item $T$ is a torus and 
    \item the action of $T$ on itself extends to $X$. 
\end{itemize}
By abuse of notation, also $X$ is called a toric variety. 
In this case, $(T, X)$ is called the {\em toric structure} on $X$. 
\item 
A prime divisor $D$ is {\em torus-invariant} if $D$ corresponds to a ray of $\sigma$. 
\end{enumerate}
\end{definition}

Given a rational strongly convex polyhedral sone $\sigma \subseteq N_{\R}$ with $d:=\dim \sigma = \dim N_{\R}$, 
the following hold for the associated affine normal toric variety $U_{\sigma} = \Spec\,k[M \cap \sigma^{\vee}]$. 
\begin{enumerate}
    \item $U_{\sigma}$ is smooth if and only if we can write $\sigma = \R_{\geq 0} u_1+ \cdots + \R_{\geq 0} u_d$ for some $\Z$-linear basis $u_1, ..., u_d \in N \cap \sigma$ of $N$ 
    \cite[Theorem 1.3.12]{CLS11}. 
    \item $U_{\sigma}$ is $\Q$-factorial if and only if $\sigma$ is simplicial, that is, 
    we can write $\sigma = \R_{\geq 0} u_1+ \cdots + \R_{\geq 0} u_d$ for some 
    $u_1, ..., u_d \in N \cap \sigma$ \cite[Proposition 4.2.7]{CLS11}. 
    In this case,  $u_1, ..., u_d$ is an $\R$-linear basis of $N_{\R}$. 
\end{enumerate}

\begin{lemma}\label{l-anQfac-toric}
Let $W$ be an affine normal $\Q$-factorial toric variety. 
Then the following hold. 
\begin{enumerate}
    \item There exists a finite surjective morphism $\mathbb A^r \times \mathbb G_m^s \to W$ for some $r, s \in \Z_{\geq 0}$. 
    \item 
    $W$ is analytically $\Q$-factorial, that is, 
    $\Spec \widehat{\MO}_{W, w}$ is $\Q$-factorial for every (possibly non-closed) point $w \in W$, where $\widehat{\MO}_{W, w}$ denotes the completion of the local ring $\MO_{W, w}$.  
\end{enumerate}
\end{lemma}


\begin{proof}
Let us show (1). 
We have $W =U_{\sigma}$ for some simplicial cone $\sigma \subseteq N_{\R}$. 
Set $d := \dim W$. 
Fix a minimal collection of generators $v_1, ..., v_d \in N \cap \sigma$ of $\sigma$ 
and set $N'$ to be the sublattice of $N$ generated by $v_1, ..., v_d$. 
Then the inclusion $N' \hookrightarrow N$ induces 
a finite surjective morphism $\mathbb A^r \times \mathbb G_m^s \to W$ for some $r, s \in \Z_{\geq 0}$ (cf.\ \cite[page 127]{CLS11}). 
Thus (1) holds. 

Let us show (2). 
By (1), we have a finite surjective morphism 
$V:= \mathbb A^r \times \mathbb G_m^s \to W$. 
Let $v_1, ..., v_r \in V$ be the points lying over $w$. 
We have a finite injective ring homomorphism: 
\[
\widehat{\MO}_{W, w} \hookrightarrow \MO_{V}(V) \otimes_{\MO_W(W)} \widehat{\MO}_{W, w} 
\simeq \varprojlim_n \MO_{V}(V)/\m_w^n \MO_{V}(V) \simeq 
\widehat{\MO}_{V, v_1} \times \cdots \times \widehat{\MO}_{V, v_r}, 
\]
where $\m_w$ denotes the maximal ideal of $\MO_{W, w}$, 
the first isomorphism holds by  \cite[Proposition 10.13]{AM69}, 
and the second one follows from \cite[Theorem 17.7]{Nag62}. 
In particular, the composite ring homomorphism 
\[
\widehat{\MO}_{W, w} \to \widehat{\MO}_{V, v_1} \times \cdots \times \widehat{\MO}_{V, v_r} \xrightarrow{{\rm pr}_1} \widehat{\MO}_{V, v_1}
\]
is a finite injective ring homomorphism to a regular local ring $\widehat{\MO}_{V, v_1}$. 
Then $\Spec \widehat{\MO}_{W, w}$ is $\Q$-factorial by the same argument as in \cite[Propostion 5.16]{KM98}. 
Thus (2) holds. 
\end{proof}



\subsection{Toric structure}\label{ss-toric-str}

\begin{notation}\label{n-W-toric}
Set $X := \mathbb A^n = \Spec k[x_1, ..., x_n]$ and 
$H_1 :=V(x_1), ..., H_n := V(x_n)$, which are the coordinate hyperplanes. 
We equip $X$ with the toric structure $(T_X, X)$ 
given by $T_X := X \setminus \bigcup_{i=1}^n H_i = \Spec k[x_1^{\pm 1}, ..., x_n^{\pm 1}]$. 

For $a_1, ..., a_n \in \Q \cap [0, 1)$ and $D := a_1 H_1 + \cdots + a_n H_n$, 
recall that we have the induced morphism (cf.\ Definition \ref{d-Qcone1}): 
\[
\pi : W := W_{X, D} = \Spec A \to X
\quad {\rm where}
\quad A := \bigoplus_{d \in \Z_{\geq 0}} R(dD)t^d \subseteq k[x_1^{\pm 1}, ..., x_n^{\pm 1}, t]. 
\]
\vspace{-1.2em}

\noindent Set 
\begin{itemize}
    \item $E_W := (\pi^*H_1)_{\red} + \cdots + (\pi^*H_n)_{\red}$, 
    \item $T_W := \Spec k[x_1^{\pm 1}, ..., x_n^{\pm 1}, t^{\pm 1}]$, and 
    \item $\Gamma := \Gamma_{X, D}$, which is the $0$-section of $\pi : W \to X$. 
\end{itemize}
We have the commutative diagrams:  
\[
\begin{tikzcd}
A= \bigoplus_{d \geq 0} R(dD)t^d \arrow{r} & k[x_1^{\pm 1}, ..., x_n^{\pm 1}, t^{\pm 1}] \\
R= k[x_1, ..., x_n]\arrow{r}\arrow[u] & k[x_1^{\pm 1}, ..., x_n^{\pm 1}]\arrow[u]
\end{tikzcd}
\hspace{20mm} 
\begin{tikzcd}
W \arrow[d, "\pi"'] & T_W \arrow[l, "j_W"'] \arrow[d] \\
X & T_X \arrow[l, "j_X"],   
\end{tikzcd}
\]
where the right diagram is obtained from the left one by applying $\Spec (-)$. 
\end{notation}

\begin{remark}\label{r-An-inductive}
We use Notation \ref{n-W-toric}. 
Fix a closed point $x \in X$. 
The following arguments (1) and (2) enable us to reduce many problems 
to the case when $x$ is the origin.  
\begin{enumerate}
    \item 
    Assume $x \not\in H_1, ..., x \not\in H_n$. 
    Then $W$ and $W_{X, 0} = X \times \mathbb A^1$ are isomorphic over some open neighbourhood of $x \in X$. 
    More precisely, for the projection $\pi': W_{X, 0} \to X$, 
    there exists an open neighbourhood $X'$ of $x \in X = \mathbb A^n$ 
    such that $
    \pi^{-1}(X')$ and $\pi'^{-1}(X') = X' \times \mathbb A^1$ 
    are isomorphic over $X'$. 
    In particular, $\pi^{-1}(x) \simeq \mathbb A^1$. 
    
    \item 
    Assume $x \in H_1, ..., x \in H_r, x \not\in H_{r+1}, ..., x \not\in H_n$ 
    for some $1 \leq r <n$. 
    Then, $W$ and $W' :=W_{X, a_1H_1+ \cdots +a_r H_r}$ are isomorphic over 
    some  open neighbourhood of $x \in X$. 
    Furthermore, 
    we have the following cartesian diagram (Lemma \ref{l-Qcone-etale}): 
    \[
    \begin{tikzcd}
        W_{X, a_1H_1+ \cdots +a_r H_r} \arrow[d, "\pi"'] \arrow[r] & W_{X', a_1H'_1+ \cdots + a_rH'_r} \arrow[d, "\pi '"]\\ 
    X=\mathbb{A}^n \arrow[r, "\rho"] & X':=\mathbb{A}^r
    \end{tikzcd}
    \]
    for the projection $\rho :\mathbb A^n \to \mathbb A^r, (c_1, ..., c_n)  \mapsto (c_1, ..., c_r)$ and 
    the coordinate hyperplanes $H'_1, ..., H'_r$ of $X' = \mathbb A^r$. 
\end{enumerate}
\end{remark}

%

\begin{proposition}\label{p-toricity1}
We use Notation \ref{n-W-toric}. 
Then the following hold. 
\begin{enumerate}
\item The induced morphism $j_W : T_W \to W$ is an open immersion. 
\item $W$ is an affine normal toric variety with the toric structure $(T_W, W)$. 
\item $\pi : (T_W, W) \to (T_X, X)$ is a toric morphism. 
\item $W$ is $\Q$-factorial. 
\end{enumerate} 
\end{proposition}

\begin{proof}
The assertion (1) follows from  
\[
A[x_1^{-1}, ..., x_n^{-1}, t^{-1}] = 
\left( \bigoplus_{d \in \Z_{\geq 0}} R(dD)t^d\right) [x_1^{-1}, ..., x_n^{-1}, t^{-1}]
= k[x_1^{\pm 1}, ..., x_n^{\pm 1}, t^{\pm 1}]. 
\]
Let us show (2). 
For any $d \in \Z_{\geq 0}$, it holds that 
\[
R(dD)t^d = R(\llcorner dD \lrcorner) t^d =  R x_1^{-\alpha_{d, 1}} \cdots x_n^{-\alpha_{d, n}} t^d, 
\]
where $\alpha_{d, 1}, ..., \alpha_{d, n} \in \Z_{\geq 0}$ are defined as follows: 
\begin{equation}\label{e1-toricity1}
\alpha_{d, 1} := \llcorner da_1 \lrcorner,\quad \ldots,\quad \alpha_{d, n} := \llcorner da_n \lrcorner. 
\end{equation}
Therefore, we obtain 
\[
A = \bigoplus_{d \geq 0} R(dD) t^d 
= \sum_{d=0}^{\infty} Rx_1^{-\alpha_{d, 1}} \cdots x_n^{-\alpha_{d, n}} t^d 
= k[\Lambda], 
\]
where $k[\Lambda]$ denotes the $k$-subalgebra of 
$k[x_1^{\pm 1}, ..., x_n^{\pm 1}, t]$  generated by 
\[
\Lambda := \{x_1, ..., x_n\} \cup \{ x_1^{-\alpha_{d, 1}} \cdots x_n^{-\alpha_{d, n}} t^d \,|\, d \in \Z_{\geq 0} \}.  
\]
Since $k[\Lambda]$ is a finitely generated $k$-algebra, 
there exists a finite subset $\Lambda' \subseteq \Lambda$ such that 
$\{x_1, ..., x_n \} \subseteq \Lambda'$ and $k[\Lambda] = k[\Lambda']$. 
Hence it follows from \cite[Proposition 1.1.14]{CLS11} that $W$ is an affine toric variety. 
Thus (2) holds. 
The assertion (3) follows from the fact that the induced morphism $T_W \to T_X$ is a homomorphism of algebraic groups. 

Let us show (4). 
Under the identification $M := \Z^{n+1} = \Hom_{\Z}(\Hom (T_W, \mathbb G_m), \Z)$, 
let $e_1, ..., e_{n+1} \in M$ be the standard basis, where $e_1, ..., e_n, e_{n+1}$ correspond to 
the variables $x_1, ..., x_n, t$.

For each $d \in \Z_{\geq 0}$, set 
\[
v_d := -\alpha_{d, 1}e_1 - \cdots - \alpha_{d, n}e_n + d e_{n+1} 
= -\sum_{i=1}^n \alpha_{d, i}e_i + d e_{n+1} . 
\]
Fix $d_0 \in \Z_{>0}$ such that $d_0 D$ is a Cartier divisor, that is,  
$d_0 a_1, ..., d_0 a_n \in \Z$. 
Set 
\[
\tau := \R_{\geq 0} e_1 + \cdots + \R_{\geq 0} e_n + \sum_{d \geq 0} \R_{\geq 0}v_d,\quad \text{ and }
\quad
\tau' := \R_{\geq 0} e_1 + \cdots + \R_{\geq 0} e_n + \R_{\geq 0} v_{d_0}. 
\]
For the time being, let us finish the proof under assumption $\tau = \tau'$. 
By $k[M \cap \tau'] \subseteq k[\Lambda] \subseteq k[M \cap \tau] = k[M \cap \tau']$, 
we obtain $A = k[\Lambda] = k[M \cap \tau']$. 
Since $\tau'$ is simplicial,  so is $\sigma:=(\tau ')^{\vee}=\tau^{\vee}$. 
Hence $W=U_{\sigma}$ is $\Q$-factorial.

It is enough to show that $\tau = \tau'$. 
The inclusion $\tau \supseteq \tau'$ is clear. 
Fix $d \in \Z_{\geq 0}$. 
It suffices to prove that $v_d \in \tau'$, which follows from 
\begin{eqnarray*}
v_d &=& -\sum_{i=1}^n \alpha_{d, i}e_i + d e_{n+1} \\
&\overset{({\rm i})}{=}&-\sum_{i=1}^n \llcorner da_i \lrcorner e_i + d e_{n+1}\\
&=&\sum_{i=1}^n (da_i - \rdown{da_i})e_i -\sum_{i=1}^n  da_i e_i + d e_{n+1}\\
&=&\sum_{i=1}^n (da_i - \rdown{da_i})e_i 
+ \frac{d}{d_0} \left(- \sum_{i=1}^n  d_0a_i e_i + d_0 e_{n+1} \right)\\
&\overset{({\rm ii})}{=}&\sum_{i=1}^n (da_i - \rdown{da_i})e_i 
+ \frac{d}{d_0} \left(- \sum_{i=1}^n  \alpha_{d_0, i} e_i + d_0 e_{n+1} \right)\\
&=&\sum_{i=1}^n (da_i - \rdown{da_i})e_i 
+\frac{d}{d_0}v_{d_0}\\
&\in & \R_{\geq 0} e_1 + \cdots + \R_{\geq 0} e_n + \R_{\geq 0} v_{d_0} = \tau', 
\end{eqnarray*}
where (i) and (ii) follow from (\ref{e1-toricity1}).  
Thus (4) holds. 
\end{proof}

\begin{remark}\label{r-toric-str}
We use Notation \ref{n-W-toric} 
and the same notation as in the proof of Proposition \ref{p-toricity1}. 
In what follows, we summarise the properties of the toric structures on $X$ and $W$. 
\begin{enumerate}
    \item Set $N := M:= \Z^{n+1}$ and $N' := M' := \Z^n$. 
    Let $e_1, ..., e_{n+1}$ (resp.\ $e'_1, ..., e'_n$) be the standard basis of $N=M=\Z^{n+1}$ (resp.\ $N'=M'=\Z^n$). 
    
    Take the projection: 
    \[
    \varphi : N \to N', \qquad e_1 \mapsto e'_1, \quad \ldots,\quad e_n \mapsto e'_n, \quad e_{n+1} \mapsto 0.
    \]
    We identify $M$ and $\Hom_{\Z}(N, \Z)$ 
    via the standard inner product (cf.\ Notation \ref{n-toric}). 
    Similarly, we identify $M'$ and $\Hom_{\Z}(N', \Z)$. 
    \item 
    For $v := v_{d_0}$, we set 
    \[
    \tau := \R_{\geq 0}e_1 + \cdots + \R_{\geq 0} e_n + \R_{\geq 0} v \subseteq M_{\R} 
    \qquad {\rm and}\qquad 
    \sigma := \tau^{\vee} \subseteq N_{\R}.
    \]
    It follows from the proof of Proposition \ref{p-toricity1}(4) that  
    \[
    W = U_{\sigma} = \Spec k[\tau \cap M]. 
    \]
    \item 
    We have 
    \[
    \sigma = \R_{\geq 0} u_1 + \cdots +\R_{\geq 0} u_n + \R_{\geq 0}e_{n+1}, 
    \]
    where each $u_i \in N$ is the primitive element such that 
    \begin{itemize}
        \item $\langle u_i, e_i \rangle >0$, 
        \item $\langle u_i, e_j \rangle =0$ for any $j \in \{1, ..., n\} \setminus \{i\}$, and 
        \item $\langle u_i, v \rangle =0$. 
    \end{itemize}
    For each $i$, it holds that 
    \[
    u_i = d_i e_i + \ell_i e_{n+1}
    \]
    for the coprime integers $d_i$ and $\ell_i$ such that 
    $d_i > \ell_i \geq 0$ and $a_i = \ell_i/d_i$.  In particular, we get 
    \[
    \varphi(u_1) = d_1 e'_1,  \quad \ldots,\quad \varphi(u_n) = d_n e'_n, \quad \varphi(e_{n+1}) =0. 
    \]
    \item 
    For $\sigma' := \sum_{i=1}^n \R_{\geq 0} e'_i \subseteq N'_{\R}$, 
    we have $\varphi_{\R}(\sigma) \subseteq \sigma'$, which induces the toric morphism 
    \[
    \pi : W = U_{\sigma} \to U_{\sigma'} = \mathbb A^n =X.
    \]
\end{enumerate}
\end{remark}

\begin{proposition}\label{p-toricity2}
We use Notation \ref{n-W-toric}. 
Then the following hold. 
\begin{enumerate}
\item For any point $x \in X$, the fibre $\pi^{-1}(x)$ is one-dimensional and 
geometrically irreducible.  
\item $\pi : W \to X$ is flat. 
\item $\Gamma$ is a torus-invariant prime divisor on $W$. 
\item For every $1 \leq i \leq n$, $(\pi^*H_i)_{\red}$ is a torus-invariant prime divisor on $W$. 
\end{enumerate} 
\end{proposition}

\begin{proof}
Let us show (1). 
We may assume that $x \in X$ is a closed point. 
If $x \not\in H_1, ..., x\not\in H_n$, then we have $\pi^{-1}(x) \simeq \mathbb A^1$ (Remark \ref{r-An-inductive}(1)). 
Hence, the problem is reduced to the case when 
$x \in H_1, ..., x \in H_r, x \not\in H_{r+1}, ..., x\not\in H_n$ 
for some $1 \leq r \leq n$. 
By Remark \ref{r-An-inductive}(2), we may assume that $r=n$, that is, 
$x$ is the origin. 
In what follows, we use the notation introduced in Remark \ref{r-toric-str}. 
It follows from 
\cite[Theorem 3.2.6 and Lemma 3.3.21]{CLS11} that 
\[
\pi^{-1}(0) = \pi^{-1}( O_X({\sigma'})) = O_W(\sigma) \amalg O_W(\widetilde{\sigma}) \subseteq \overline{O_W(\widetilde{\sigma})}
\]
for the orbits  $O_X(-)$ and $O_W(-)$, 
where $\sigma$ and $\widetilde{\sigma}$ are the faces such that 
\begin{itemize}
    \item $\dim \sigma =n+1, \dim \widetilde{\sigma} =n$,
    \item $\varphi(\sigma ) =\varphi(\widetilde{\sigma}) = \sigma'$, and 
    \item $\widetilde{\sigma} \prec \sigma$. 
\end{itemize}
Since $\pi^{-1}(0)$ is a closed subset of $W$, we have $\pi^{-1}(0) = \overline{O_W(\widetilde{\sigma})}$, which is one-dimensional and irreducible. 
Thus (1) holds. 
The assertion (2) follows from (1) and the fact that 
$X$ is smooth and $W$ is Cohen--Macaulay.

Let us show (3). 
By $\Gamma \simeq X$, $\Gamma$ is a prime divisor on $W$. 
Hence it suffices to show that $\Gamma$ is torus-invariant. 
We have ring homomorphisms: 
\[
A = \bigoplus_{d \geq 0} R(dD) t^d \hookrightarrow k[x_1^{\pm 1}, ..., x_n^{\pm 1}, t] 
\hookrightarrow k[x_1^{\pm 1}, ..., x_n^{\pm 1}, t^{\pm 1}], 
\]
corresponding to $T_W$-equivariant open immersions: 
\[
W = \Spec A \hookleftarrow \mathbb G_m^n \times \mathbb A^1\hookleftarrow \mathbb G_m^{n+1} =T_W. 
\]
Since $\Gamma|_{\mathbb G_m^n \times \mathbb A^1}$ is $T_W$-invariant, 
also $\Gamma$ itself is $T_W$-invariant. 
Thus (3) holds. 

Let us show (4). 
Fix $1 \leq i \leq n$. 
Recall that $H_i \subseteq X$ is a torus-invariant prime divisor on $X$. 
Therefore, its set-theoretic inverse image $\pi^{-1}(H_i)$ is stable under the $T_W$-action. 
Therefore, it suffices to show that 
the effective Cartier divisor $\pi^*H_i$ is irreducible, 
which follows from (1) and (2). 
Thus (4) holds. 
\end{proof}

\begin{proposition}\label{p-toricity3}
We use Notation \ref{n-W-toric}. 
Then the following hold. 
\begin{enumerate}
\item $(W, \Gamma +E_W)$ is lc. 
    \item $(W, \Gamma + cE_W)$ is plt for any $0 \leq c <1$. 
    In particular, $W$ is klt.  
\end{enumerate} 
\end{proposition}

\begin{proof}
Let us show (1). 
Since $(\pi^*H_1)_{\red}, \cdots, (\pi^*H_n)_{\red}, \Gamma$ are torus-invariant prime divisors, 
$\Gamma + E_W = \Gamma + (\pi^*H_1)_{\red} + \cdots + (\pi^*H_n)_{\red}$ 
is a sum of torus-invariant prime divisors (Proposition \ref{p-toricity2}(3)(4)). 
Therefore, $(W, \Gamma +E_W)$ is lc by \cite[Proposition 11.4.24]{CLS11}. 
Thus (1) holds.

Let us show (2). 
By construction, 
$\Gamma|_{W \setminus E_W}$ is a smooth prime divisor on $W \setminus E_W$. 
Hence $(W \setminus E_W, (\Gamma + cE_W)|_{W \setminus E_W})$ is plt. 
This, together with (1), implies that $(W, \Gamma + cE_W)$ is plt. 
Thus (2) holds. 
\qedhere

\end{proof}

\begin{theorem}\label{t-toricity}
Let $X$ be a smooth variety and let $D$ be a $\Q$-divisor 
such that $\{D\}$ is simple normal crossing. 
For $W := W_{X, D}$ and $\Gamma := \Gamma_{X, D}$, 
we set $E := \{ D\}_{\red}$ and $E_W := (\pi^*E)_{\red}$, where $\pi : W \to X$ denotes the induced morphism. 
Then the following hold. 
\begin{enumerate}
    \item $W$ is $\Q$-factorial. 
    \item $(W, \Gamma +E_W)$ is lc. 
    \item $(W, \Gamma +cE_W)$ is plt for any $0 \leq c <1$. 
\end{enumerate} 
\end{theorem}

\begin{proof}
Set $n := \dim X$. 
The problem is local on $X$. 
Fix a closed point $x \in X$, around which we shall work. 
Let $ D= a_1 D_1 + \cdots + a_m D_m$ be the  irreducible decomposition. 
We may assume that the following hold. 
\begin{itemize}
    \item $D=\{D\}$, that is, $0 \leq a_1 < 1, ..., 0 \leq a_m <1$. 
    \item $x \in D_1, ...,  x \in D_m$. In particular, $m \leq n$. 
\end{itemize}
Fix an \'etale morphism $\alpha : X \to Z :=\mathbb A^n$ such that 
$D_1 = \alpha^*H_1, ..., D_m =\alpha^*H_m$ hold and $\alpha(x)$ is the origin, 
where $H_1, ..., H_n$ denote the coordinate hyperplanes. 
For $D_Z := a_1 H_1+ \cdots + a_m H_m$, we have the cartesian diagram 
(Lemma \ref{l-Qcone-etale}): 
\[
\begin{tikzcd}
    X \arrow[d, "\alpha"'] & W_{X, D} = W \arrow[l, "\pi"']\arrow[d, "\beta"]\\
Z=\mathbb A^n & W_{Z, D_Z}.\arrow[l, "\pi_Z"] 
\end{tikzcd}
\]
Since $W_{Z, D_Z}$ is a normal $\Q$-factorial toric variety (Proposition \ref{p-toricity1}), 
$W_{Z, D_Z}$ is analytically $\Q$-factorial (Lemma \ref{l-anQfac-toric}). 
Since $\beta : W \to W_{Z, D_Z}$ is \'etale, also $\Spec \MO_{W, w}$ is 
$\Q$-factorial for any closed point $w \in W$. 
Therefore, $W$ is $\Q$-factorial, that is, (1) holds. 
By $\alpha^*\Gamma_{Z, D_Z} = \Gamma$ (Lemma \ref{l-funct-Gamma}, Lemma \ref{l-Qcone-etale}), 
(2) and (3) follow from Proposition \ref{p-toricity3}. 
\end{proof}

\subsection{Different}\label{ss-Diff}

\begin{definition}\label{d-Qfac-index}
\hphantom{a}\\
\vspace{-1em}
\begin{enumerate}
    \item Let $W$ be a klt surface such that $W$ has a unique singular point $Q$. 
Let $\psi : W' \to W$ be the minimal resolution of $W$. 
For the $\psi$-exceptional prime divisors $E_1, ..., E_r$ and 
the intersection matrix $(-E_i \cdot E_j)$, we set  
\[
d_{W, Q} := \det (-E_i \cdot E_j), 
\]
Note that 
$d_Q$ does not depends on the choice of  the order of $E_1, ..., E_n$
\item 
Let $W$ be a klt surface. 
For a singular point $Q$ of $W$ and 
an open neighbourhood $W'$ of $Q \in W$ such that 
$Q$ is a unique singular point of $W',$ 
we set 
\[
d_Q := d_{W, Q} := d_{W', Q}, 
\]
where $d_{W', Q}$ is the integer defined in (1). 
\end{enumerate}
\end{definition}

\begin{lemma}\label{l-plt-chain}
Let $(W, \Gamma)$ be a two-dimensional plt pair, where $\Gamma$ is a prime divisor. 
Assume that $W$ has a  unique singular point $Q$ and 
that $Q \in \Gamma$. 
Then the following hold. 
\begin{enumerate}
\item 
For the different $\Diff_{\Gamma}$ defined by 
$(K_W+\Gamma)|_{\Gamma} = K_{\Gamma} + \Diff_{\Gamma}$, 
it holds that 
\[
\Diff_{\Gamma} = \frac{d_Q -1}{d_Q}. 
\]
    \item 
    $d_Q$ is equal to the $\Q$-factorial index of $W$, that is, 
    $d_Q$ is the minimum positive integer satisfying the following property $(*)$. 
    \begin{enumerate}
        \item[$(*)$] If $D$ is a Weil divisor on $W$, then $d_Q D$ is Cartier. 
    \end{enumerate}
\end{enumerate}
\end{lemma}



\begin{proof}
The assertion (1) follows from \cite[Theorem 3.36]{kollar13}. 
Let us show (2). 
By \cite[Lemma 2.2]{CTW15b}, 
the $\Q$-factorial index $d'_Q$ is a positive integer with $d_Q \geq d'_Q$. 
Let $\psi : W' \to W$ be the minimal resolution of $W$. 
Set $\Gamma' := \psi^{-1}_*\Gamma$. 
For the irreducible decomposition $\Ex(\psi) = E_1 \cup \cdots \cup E_r$, the extended dual graph 
is a chain \cite[Theorem 3.36]{kollar13}: 
\[
\Gamma' - E_1 - \cdots - E_n. 
\]
There exist $c_1, ..., c_n \in \Q_{\geq 0}$ such that 
\[
K_{W'} + \Gamma' + c_1E_1+ \cdots + c_n E_n = \psi^*(K_W+ \Gamma). 
\]
By (1), we have $c_1 = \frac{d_Q-1}{d_Q}$. 
Since $d'_Q(K_W + \Gamma)$ is Cartier, we obtain 
\[
d'_Q \cdot \frac{d_Q-1}{d_Q}  = d'_Q c_1 \in \Z,
\]
which implies $d'_Q \in d_Q \Z$. 
Combining with $d_Q \geq d'_Q$, we get $d_Q = d'_Q$. 
Thus (2) holds. 
\end{proof}

\begin{proposition}\label{p-Diff}
Set $X :=\mathbb A^1$ and let $P$ be a closed point. 
For coprime integers $d > \ell >0$, we set $D:= \frac{\ell}{d} P$. 
Set $(W, \Gamma):=(W_{X, D}, \Gamma_{X, D})$. 
Let $Q \in \Gamma$ be the closed point lying over $P$. 
Then the following hold. 
\begin{enumerate}
\item $(W, \Gamma)$ is a two-dimensional affine plt pair. 
\item $Q$ is a unique singular point of $W$. 
    \item $d_Q = d$. 
    \item For the different $\Diff_{\Gamma}$ given by 
$(K_W+\Gamma)|_{\Gamma} = K_{\Gamma} + \Diff_{\Gamma}$, 
it holds that 
\[
\Diff_{\Gamma} = \frac{d-1}{d}  Q. 
\]
\end{enumerate}
\end{proposition}

\begin{proof}
We may assume that $P$ is the origin of $\mathbb A^1$. 
By setting $(n, H_1, a_1) :=(1, P, \frac{\ell}{d})$, 
we may use Notation \ref{n-W-toric}. 
Hence $(W, \Gamma)$ is a plt pair (Proposition \ref{p-toricity3}(2)). 
Thus (1) holds. 
By Remark \ref{r-toric-str}, 
we have $W=U_{\sigma}$ for 
\[
\sigma := \R_{\geq 0} u + \R_{\geq 0} e_2, \qquad u :=de_1 + \ell e_2. 
\]
By $d \geq 2$, we have $\Z u +\Z e_2 = \Z(de_1) + \Z e_2 \neq \Z^2$, 
and hence $W$ is a singular affine normal toric surface with a unique singular point $Q'$. 
Since $Q'$ is the torus-invariant point, 
$Q'$ is set-theoretically equal to the intersection of the torus invariant prime divisors $\Gamma$ and $(\pi^*P)_{\red}$ (Proposition \ref{p-toricity2}(3)(4)). 
Therefore, we obtain $Q = Q'$. 
Thus (2) holds.

It follows from (1), (2), and Lemma \ref{l-plt-chain}(1) that (3) implies (4). 
It suffices to show (3). 
We use the notation in Remark \ref{r-toric-str}. 
In particular, we equip $\Z^2 =N = M$ with the standard inner product: $\langle e_i, e_j \rangle =\delta_{ij}$. 
Note that $(W, \Gamma)$ is a plt pair and $\Gamma$ passes through $Q$. 
Therefore, it suffices to show that 
the $\Q$-factorial index is equal to $d$ (Lemma \ref{l-plt-chain}(2)).

Take $a_1, a_2 \in \Z$ and set 
\[
D := a_1 D_1 + a_2 D_2,
\]
where $D_1$ and $D_2$ are the torus invariant prime divisors on $W$ corresponding to $u$ and $e_2$, respectively.
We define $\beta_1, \beta_2 \in \Q$ such that 
the equality 
\[
(a_1, a_2) = (\langle m, u \rangle, \langle m, e_2 \rangle) 
\quad \text{ holds}\quad\text{for }\quad
m := \beta_1 e_1 + \beta_2 e_2 \in M. 
\]
Then we have that 
\[
\begin{pmatrix}
a_1\\
a_2
\end{pmatrix}
=
\begin{pmatrix}
\langle m, u \rangle \\
\langle m, e_2 \rangle
\end{pmatrix}
=
\begin{pmatrix}
\langle \beta_1 e_1 + \beta_2 e_2, de_1 + \ell e_2\rangle \\
\langle \beta_1 e_1 + \beta_2 e_2, e_2\rangle
\end{pmatrix}
=
\begin{pmatrix}
\beta_1 d + \beta_2\ell \\
\beta_2
\end{pmatrix}
=
\begin{pmatrix}
d & \ell\\
0 & 1
\end{pmatrix}
\begin{pmatrix}
\beta_1\\
\beta_2
\end{pmatrix}, 
\]
which implies 
\[
\begin{pmatrix}
\beta_1\\
\beta_2
\end{pmatrix}
=
\begin{pmatrix}
d & \ell\\
0 & 1
\end{pmatrix}^{-1}
\begin{pmatrix}
a_1\\
a_2
\end{pmatrix} = 
\frac{1}{d}
\begin{pmatrix}
1 & -\ell\\
0 & d
\end{pmatrix}
\begin{pmatrix}
a_1\\
a_2
\end{pmatrix} = 
\begin{pmatrix}
\frac{1}{d} & -\frac{\ell}{d}\\
0 & 1
\end{pmatrix}
\begin{pmatrix}
a_1\\
a_2
\end{pmatrix}
=
\begin{pmatrix}
\frac{a_1 - \ell a_2}{d}\\
a_2
\end{pmatrix}.
\]
By \cite[Theorem 4.2.8]{CLS11}, $D= a_1 D_1 + a_2 D_2$ is Cartier if and only if $\frac{a_1-\ell a_2}{d} \in \Z$ and $a_2 \in \Z$.
Since every Weil divisor on $W$ is linearly equivalent to 
a torus invariant divisor \cite[Theorem 4.1.3]{CLS11},
the $\Q$-factorial index of $W$ is equal to $d$. 
Thus (3) holds. 
\end{proof}

\begin{theorem}\label{t-Diff}
Let $X$ be a smooth variety and let $D$ be a $\Q$-divisor such that 
$\{ D\} = a_1 D_1+ \cdots + a_m D_m$ is simple normal crossing, 
where $a_1, ..., a_m \in \Q_{>0}$. 
For each $1 \leq i \leq m$, 
let $\ell_i$ and $d_i$ be coprime integers such that 
$a_i = \ell_i /d_i$ and $0 < \ell_i <d_i$. 
Set $(W, \Gamma) :=(W_{X, D}, \Gamma_{X, D})$. 
Then, for the different $\Diff_{\Gamma}$ defined by 
\[
(K_W+\Gamma)|_{\Gamma} = K_{\Gamma} + \Diff_{\Gamma},  
\]
it holds that 
\[
\Diff_{\Gamma} = \sum_{i=1}^m \frac{d_i -1}{d_i}D_{\Gamma, i}, 
\]
where $D_{\Gamma, i}$ denotes the pullback of $D_i$ under 
the composite isomorphism $\Gamma \hookrightarrow W \xrightarrow{\pi} X$. 
\end{theorem}

\begin{proof}
In order to compute the coefficient of $D_{\Gamma, 1}$ in $\Diff_{\Gamma}$, 
we may replace $X$ by an open subset of $X$ which intersects $D_1$. 
Hence, the problem is reduced to the case when $D = \{ D\} = a_1D_1$. 
Taking a suitable \'etale coordinate, 
we may assume that 
$X=\mathbb A^n$ and $D_1$ is a coordinate hyperplane (Lemma \ref{l-Qcone-etale}). 
 Since the problem is stable under smooth base change, 
we may further assume that $X= \mathbb A^1$ and $D_1$ is a closed point (Lemma \ref{l-Qcone-etale}). 
Then the assertion follows from Proposition \ref{p-Diff}. 
\end{proof}

\begin{theorem}\label{t-klt}
Let $X$ be a smooth projective variety and 
let $\Delta$ be a simple normal crossing $\Q$-divisor such that 
$\llcorner \Delta \lrcorner=0$ and $\Delta$ has standard coefficients. 
Assume that 
\begin{enumerate}
\renewcommand{\labelenumi}{(\roman{enumi})}
    \item $D:=-(K_X+\Delta)$ is ample, and 
    \item $\rho(X)=1$. 
\end{enumerate}
For $W := W_{X, D}$ and $V:=V_{X, D}$, let
$\mu :W \to V$ be the induced morphism. 
Then the following hold. 
\begin{enumerate}
    \item $\rho(W/V)=1$. 
    \item $-\Gamma$ is $\mu$-ample. 
    \item $V$ is $\Q$-factorial. 
    \item $V$ is klt. 
\end{enumerate}
\end{theorem}

\begin{proof}
Note that the following property (iii) holds 
by (i), (ii), and \cite[Corollary 6.5]{tanakainseparable}. 
\begin{enumerate}
\item[(iii)] Any numerically trivial Cartier divisor on $X$ is torsion. 
\end{enumerate}
The assertion (1) follows from (ii) and the fact that 
$\mu : W \to V$ is a projective birational morphism 
with $\Ex(\mu) = \Gamma$ (Theorem \ref{t-Qcone-birat}). 
Fix a curve $C$ satisfying $C \subseteq \Gamma$. 

Let us show (2). 
Pick an effective Cartier divisor $H_V$ on $V$ passing through the vertex $v \in V$. 
It holds that 
\[
\mu^* H_V = \mu^{-1}_* H_V + \alpha \Gamma 
\]
for some $\alpha \in \Z_{>0}$. 
Since $\mu^{-1}_*H_{V}$ intersects $\Gamma$, 
it follows from $\rho(\Gamma)=\rho(X)=1$ that $\mu^{-1}_*H \cdot C>0$. 
Therefore, it holds that  
\[
0 = \mu^* H_V \cdot C = \mu^{-1}_* H_V\cdot C+ \alpha \Gamma \cdot C 
> \alpha \Gamma \cdot C. 
\]
By $\alpha >0$, we get $\Gamma \cdot C<0$. 
It follows from (1) that $-\Gamma$ is $\mu$-ample. 
Thus (2) holds.

We now prove (3) and (4) for the case when 
$k$ is of positive characteristic.  
Let us show (3). 
Fix a prime divisor $D_V$ on $V$ and set $D_W := \mu_*^{-1}D_V$. 
We may assume that $v \in D_V$. 
By (1) and (2), 
there is $\beta \in \Q$ such that $D_W + \beta \Gamma \equiv_{\mu} 0$. 
In particular, $D_W + \beta \Gamma$ is $\mu$-nef and $\mu$-big. 
By (1) and $v \in D_V$, $D_W$ is $\mu$-ample. 
It holds by (2) that $\beta > 0$. 
It follows from \cite[Lemma 2.18(1)]{CT17} that $\mathbb E_{\mu}(D_W + \beta \Gamma) \subseteq \Gamma$. 
By (iii) and \cite[Proposition 2.20]{CT17}, 
$D_W + \beta \Gamma$ is $\mu$-semi-ample. 
Hence, its pushforward $\mu_*(D_W + \beta \Gamma) = D_V$ is $\Q$-Cartier. Thus (3) holds. 
Let us show (4). 
By (3), there is $\gamma \in \Q$ such that 
\[
K_W + (1-\gamma) \Gamma  =\mu^*K_V.  
\]
For the different $\Delta_{\Gamma}$ defined by 
$(K_W + \Gamma)|_{\Gamma}  = (K_{\Gamma} +\Delta_{\Gamma})$, the following holds: 
\[
(K_W + \Gamma) \cdot C = (K_{\Gamma} +\Delta_{\Gamma}) \cdot C = (K_X+ \Delta) \cdot C_X <0, 
\]
where $C_X := \pi(C)$ and the second equality holds by Theorem \ref{t-Diff}. 
This inequality, together with $\Gamma \cdot C<0$, implies that $\gamma >0$. 
Then $(W, (1-\gamma) \Gamma))$ is sub-klt by Theorem \ref{t-toricity}, that is, 
$(1-\gamma)\Gamma$ is a (possibly non-effective) $\Q$-divisor and 
$a(E, W, (1-\gamma) \Gamma)) >-1$ for any prime divisor $E$ over $W$. 
Therefore, $(V, 0)$ is klt. 
This completes the proof for the case when 
$k$ is of positive characteristic.  

It suffices to show (3) and (4) for the case when 
$k$ is of characteristic zero. 
By the above argument, it holds that 
\[
K_W + \Gamma \equiv_{\mu} \gamma \Gamma 
\]
for some $\gamma >0$. 
By (2), $(W, B)$ is klt and $-(K_W+B)$ is $\mu$-ample 
for $B := (1-\delta )\Gamma$ and $0 < \delta \ll 1$. 
Hence $\mu : W \to V$ is a $(K_W+B)$-negative divisorial contraction. 
It follows from the same argument as in \cite[Corollary 3.18]{KM98} that $V$ is $\Q$-factorial. 
By the same argument as above, (4) holds. 
\end{proof}

%% file: main.bbl

%% file: main.bbl
\begin{bibdiv}
\begin{biblist}

\bib{ABL20}{article}{
      author={Arvidsson, Emelie},
      author={Bernasconi, Fabio},
      author={Lacini, Justin},
       title={On the {K}awamata-{V}iehweg vanishing theorem for log del {P}ezzo
  surfaces in positive characteristic},
        date={2022},
        ISSN={0010-437X},
     journal={Compos. Math.},
      volume={158},
      number={4},
       pages={750\ndash 763},
         url={https://doi.org/10.1112/S0010437X22007394},
      review={\MR{4438290}},
}

\bib{AM69}{book}{
      author={Atiyah, M.~F.},
      author={Macdonald, I.~G.},
       title={Introduction to commutative algebra},
   publisher={Addison-Wesley Publishing Co., Reading, Mass.-London-Don Mills,
  Ont.},
        date={1969},
      review={\MR{0242802}},
}

\bib{Ber}{article}{
      author={Bernasconi, Fabio},
       title={Kawamata-{V}iehweg vanishing fails for log del {P}ezzo surfaces
  in characteristic 3},
        date={2021},
        ISSN={0022-4049},
     journal={J. Pure Appl. Algebra},
      volume={225},
      number={11},
       pages={Paper No. 106727, 16},
         url={https://doi-org.utokyo.idm.oclc.org/10.1016/j.jpaa.2021.106727},
      review={\MR{4228436}},
}

\bib{BH93}{book}{
      author={Bruns, Winfried},
      author={Herzog, J\"{u}rgen},
       title={Cohen-{M}acaulay rings},
      series={Cambridge Studies in Advanced Mathematics},
   publisher={Cambridge University Press, Cambridge},
        date={1993},
      volume={39},
        ISBN={0-521-41068-1},
      review={\MR{1251956}},
}

\bib{CLS11}{book}{
      author={Cox, David~A.},
      author={Little, John~B.},
      author={Schenck, Henry~K.},
       title={Toric varieties},
      series={Graduate Studies in Mathematics},
   publisher={American Mathematical Society, Providence, RI},
        date={2011},
      volume={124},
        ISBN={978-0-8218-4819-7},
         url={https://doi.org/10.1090/gsm/124},
      review={\MR{2810322}},
}

\bib{CT19-2}{article}{
      author={Cascini, Paolo},
      author={Tanaka, Hiromu},
       title={Purely log terminal threefolds with non-normal centres in
  characteristic two},
        date={2019},
        ISSN={0002-9327},
     journal={Amer. J. Math.},
      volume={141},
      number={4},
       pages={941\ndash 979},
         url={https://doi.org/10.1353/ajm.2019.0025},
      review={\MR{3992570}},
}

\bib{CT17}{article}{
      author={Cascini, Paolo},
      author={Tanaka, Hiromu},
       title={Relative semi-ampleness in positive characteristic},
        date={2020},
        ISSN={0024-6115},
     journal={Proc. Lond. Math. Soc. (3)},
      volume={121},
      number={3},
       pages={617\ndash 655},
         url={https://doi.org/10.1112/plms.12323},
      review={\MR{4100119}},
}

\bib{CTW15b}{article}{
      author={Cascini, Paolo},
      author={Tanaka, Hiromu},
      author={Witaszek, Jakub},
       title={On log del {P}ezzo surfaces in large characteristic},
        date={2017},
        ISSN={0010-437X},
     journal={Compos. Math.},
      volume={153},
      number={4},
       pages={820\ndash 850},
         url={http://dx.doi.org/10.1112/S0010437X16008265},
}

\bib{CTW15a}{article}{
      author={Cascini, Paolo},
      author={Tanaka, Hiromu},
      author={Witaszek, Jakub},
       title={Klt del {P}ezzo surfaces which are not {G}lobally {F}-split},
        date={2018},
     journal={Int. Math. Res. Not. IMRN},
      number={7},
       pages={2135\ndash 2155},
}

\bib{Demazure}{article}{
      author={Demazure, Michel},
       title={Automorphismes et d\'eformations des vari\'et\'es de {B}orel},
        date={1977},
        ISSN={0020-9910},
     journal={Invent. Math.},
      volume={39},
      number={2},
       pages={179\ndash 186},
         url={http://dx.doi.org/10.1007/BF01390108},
}

\bib{Dem88}{incollection}{
      author={Demazure, Michel},
       title={Anneaux gradu\'{e}s normaux},
        date={1988},
   booktitle={Introduction \`a la th\'{e}orie des singularit\'{e}s, {II}},
      series={Travaux en Cours},
      volume={37},
   publisher={Hermann, Paris},
       pages={35\ndash 68},
      review={\MR{1074589}},
}

\bib{fga2005}{book}{
      author={Fantechi, B.},
      author={G{\"o}ttsche, L.},
      author={Illusie, L.},
      author={Kleiman, S.~L.},
      author={Nitsure, N.},
      author={Vistoli, A.},
       title={Fundamental algebraic geometry},
   publisher={American Mathematical Society, Providence, RI},
        date={2005},
}

\bib{gk14}{article}{
      author={Graf, Patrick},
      author={Kov\'{a}cs, S\'{a}ndor~J.},
       title={An optimal extension theorem for 1-forms and the
  {L}ipman-{Z}ariski conjecture},
        date={2014},
     journal={Doc. Math.},
      volume={19},
       pages={815\ndash 830},
}

\bib{GKKP}{article}{
      author={Greb, Daniel},
      author={Kebekus, Stefan},
      author={Kov\'{a}cs, S\'{a}ndor~J.},
      author={Peternell, Thomas},
       title={Differential forms on log canonical spaces},
        date={2011},
        ISSN={0073-8301},
     journal={Publ. Math. Inst. Hautes \'{E}tudes Sci.},
      number={114},
       pages={87\ndash 169},
         url={https://doi.org/10.1007/s10240-011-0036-0},
      review={\MR{2854859}},
}

\bib{GNT16}{article}{
      author={Gongyo, Yoshinori},
      author={Nakamura, Yusuke},
      author={Tanaka, Hiromu},
       title={Rational points on log {F}ano threefolds over a finite field},
        date={2019},
        ISSN={1435-9855},
     journal={J. Eur. Math. Soc. (JEMS)},
      volume={21},
      number={12},
       pages={3759\ndash 3795},
         url={https://doi.org/10.4171/JEMS/913},
      review={\MR{4022715}},
}

\bib{graf21}{article}{
      author={Graf, Patrick},
       title={Differential forms on log canonical spaces in positive
  characteristic},
        date={2021},
     journal={J. Lond. Math. Soc. (2)},
      volume={104},
      number={5},
       pages={2208\ndash 2239},
}

\bib{Got-Wat78}{article}{
      author={Goto, Shiro},
      author={Watanabe, Keiichi},
       title={On graded rings. {I}},
        date={1978},
        ISSN={0025-5645},
     journal={J. Math. Soc. Japan},
      volume={30},
      number={2},
       pages={179\ndash 213},
         url={https://doi.org/10.2969/jmsj/03020179},
      review={\MR{494707}},
}

\bib{hartshorne_deformation}{book}{
      author={Hartshorne, R.},
       title={Deformation theory},
      series={Graduate Texts in Mathematics},
   publisher={Springer, New York},
        date={2010},
      volume={257},
        ISBN={978-1-4419-1595-5},
         url={http://dx.doi.org/10.1007/978-1-4419-1596-2},
}

\bib{hartshorne77}{book}{
      author={Hartshorne, R.},
       title={Algebraic geometry},
   publisher={Springer-Verlag, New York},
        date={1977},
}

\bib{HMX14}{article}{
      author={Hacon, Christopher~D.},
      author={McKernan, James},
      author={Xu, Chenyang},
       title={A{CC} for log canonical thresholds},
        date={2014},
        ISSN={0003-486X,1939-8980},
     journal={Ann. of Math. (2)},
      volume={180},
      number={2},
       pages={523\ndash 571},
         url={https://doi.org/10.4007/annals.2014.180.2.3},
      review={\MR{3224718}},
}

\bib{illusie_de_rham_witt}{article}{
      author={Illusie, Luc},
       title={Complexe de de\thinspace {R}ham-{W}itt et cohomologie
  cristalline},
        date={1979},
        ISSN={0012-9593},
     journal={Ann. Sci. \'Ecole Norm. Sup. (4)},
      volume={12},
      number={4},
       pages={501\ndash 661},
         url={http://www.numdam.org/item?id=ASENS_1979_4_12_4_501_0},
}

\bib{Kaw3}{article}{
      author={Kawakami, Tatsuro},
       title={Bogomolov-{S}ommese vanishing and liftability for surface pairs
  in positive characteristic},
        date={2022},
        ISSN={0001-8708},
     journal={Adv. Math.},
      volume={409},
       pages={Paper No. 108640},
         url={https://doi.org/10.1016/j.aim.2022.108640},
      review={\MR{4473638}},
}

\bib{Kaw4}{article}{
      author={Kawakami, Tatsuro},
       title={Extendability of differential forms via {C}artier operators},
        date={2022},
     journal={arXiv preprint arXiv:2207.13967, to appear in
  J.~Eur.~Math.~Soc.~(JEMS)},
}

\bib{KM98}{book}{
      author={Koll{\'a}r, J.},
      author={Mori, S.},
       title={Birational {G}eometry of {A}lgebraic {V}arieties},
      series={Cambridge {T}racts in {M}athematics},
   publisher={Cambridge University Press},
        date={1998},
      volume={134},
}

\bib{Kawakami-Nagaoka22}{article}{
      author={Kawakami, Tatsuro},
      author={Nagaoka, Masaru},
       title={Pathologies and liftability of {D}u {V}al del {P}ezzo surfaces in
  positive characteristic},
        date={2022},
        ISSN={0025-5874},
     journal={Math. Z.},
      volume={301},
      number={3},
       pages={2975\ndash 3017},
         url={https://doi.org/10.1007/s00209-022-02998-6},
      review={\MR{4437346}},
}

\bib{kollar13}{book}{
      author={Koll{\'a}r, J{\'a}nos},
       title={Singularities of the minimal model program},
      series={Cambridge Tracts in Mathematics},
   publisher={Cambridge University Press, Cambridge},
        date={2013},
      volume={200},
        ISBN={978-1-107-03534-8},
         url={http://dx.doi.org/10.1017/CBO9781139547895},
}

\bib{Kol94}{incollection}{
      author={Koll\'{a}r, J\'{a}nos},
       title={Log surfaces of general type; some conjectures},
        date={1994},
   booktitle={Classification of algebraic varieties ({L}'{A}quila, 1992)},
      series={Contemp. Math.},
      volume={162},
   publisher={Amer. Math. Soc., Providence, RI},
       pages={261\ndash 275},
         url={https://doi-org.utokyo.idm.oclc.org/10.1090/conm/162/01538},
      review={\MR{1272703}},
}

\bib{KTTWYY1}{article}{
      author={Kawakami, Tatsuro},
      author={Takamatsu, Teppei},
      author={Tanaka, Hiromu},
      author={Witaszek, Jakub},
      author={Yobuko, Fuetaro},
      author={Yoshikawa, Shou},
       title={Quasi-${F}$-splittings in birational geometry},
        date={2022},
     journal={arXiv preprint arXiv:2208.08016},
}

\bib{Nag62}{book}{
      author={Nagata, Masayoshi},
       title={Local rings},
      series={Interscience Tracts in Pure and Applied Mathematics, No. 13},
   publisher={Interscience Publishers (a division of John Wiley \& Sons, Inc.),
  New York-London},
        date={1962},
      review={\MR{0155856}},
}

\bib{tanaka12}{article}{
      author={Tanaka, H.},
       title={Minimal models and abundance for positive characteristic log
  surfaces},
        date={2014},
        ISSN={0027-7630},
     journal={Nagoya Math. J.},
      volume={216},
       pages={1\ndash 70},
         url={http://dx.doi.org/10.1215/00277630-2801646},
}

\bib{tanakainseparable}{article}{
      author={Tanaka, H.},
       title={Behavior of canonical divisors under purely inseparable base
  changes},
        date={2016},
        ISSN={0075-4102},
     journal={J. Reine Angew. Math.},
}

\bib{tanaka16_excellent}{article}{
      author={Tanaka, Hiromu},
       title={Minimal model program for excellent surfaces},
        date={2018},
        ISSN={0373-0956},
     journal={Ann. Inst. Fourier (Grenoble)},
      volume={68},
      number={1},
       pages={345\ndash 376},
         url={http://aif.cedram.org/item?id=AIF_2018__68_1_345_0},
}

\bib{tanaka22}{article}{
      author={Tanaka, Hiromu},
       title={Vanishing theorems of {K}odaira type for {W}itt canonical
  sheaves},
        date={2022},
        ISSN={1022-1824},
     journal={Selecta Math. (N.S.)},
      volume={28},
      number={1},
       pages={Paper No. 12, 50},
}

\bib{wat81}{article}{
      author={Watanabe, Keiichi},
       title={Some remarks concerning {D}emazure's construction of normal
  graded rings},
        date={1981},
        ISSN={0027-7630},
     journal={Nagoya Math. J.},
      volume={83},
       pages={203\ndash 211},
         url={http://projecteuclid.org/euclid.nmj/1118786485},
      review={\MR{632654}},
}

\bib{Watanabe91}{article}{
      author={Watanabe, K.},
       title={{$F$}-regular and {$F$}-pure normal graded rings},
        date={1991},
        ISSN={0022-4049},
     journal={J. Pure Appl. Algebra},
      volume={71},
      number={2-3},
       pages={341\ndash 350},
         url={http://dx.doi.org/10.1016/0022-4049(91)90157-W},
}

\end{biblist}
\end{bibdiv}
